\documentclass[a4paper,reqno]{amsart}

\usepackage{amsmath}
\usepackage{amsthm}
\usepackage{amsfonts}
\usepackage{amssymb}
\usepackage{float}

\newcounter{rh}
\newcounter{lp}
\usepackage{graphicx}
\usepackage{subcaption}

\usepackage{pdfsync}

\usepackage{asymptote}
\usepackage{verbatim}

\usepackage{subcaption}

\usepackage{color}

\usepackage{hyperref}

\usepackage{graphicx}

\hypersetup{
    colorlinks=false,
    pdfborder={0 0 0},
    pdfstartview={XYZ null null 1.00}
}

\numberwithin{equation}{section}

  \def\<{\langle}
  \def\>{\rangle}

  \def\ve{\varepsilon}

  \def\R{\mathbb{R}}

  \def\C{\mathbb{C}}

  \def\AA{\mathcal{A}}
  \def\CC{\mathcal{C}}
  \def\UU{\mathcal{U}}
  
  \def\RR{\mathcal{R}}

  \def\KK{\mathcal{K}}
  \def\LL{\mathcal{L}}

  \def\t{\widetilde}

\theoremstyle{plain}
  \newtheorem{theorem}{Theorem}[section]

  \newtheorem{proposition}[theorem]{Proposition}
  
  \newtheorem{lemma}[theorem]{Lemma}

\theoremstyle{definition}

  \newtheorem{remark}[theorem]{Remark}

\begin{document}



\title[Total integrals of Ablowitz-Segur...]{Total integrals of Ablowitz-Segur solutions for the inhomogeneous Painlev\'e II equation}



\author{Piotr Kokocki}
\address{\noindent Faculty of Mathematics and Computer Science \newline Nicolaus Copernicus University \newline Chopina 12/18, 87-100 Toru\'n, Poland}
\email{pkokocki@mat.umk.pl}
\thanks{The researches supported by the MNiSW Iuventus Plus Grant no. 0338/IP3/2016/74}



\begin{abstract}
In this paper, we establish a formula determining the value of the Cauchy integrals of the real and purely imaginary Ablowitz-Segur solutions for the inhomogeneous second Painlev\'e equation. Our approach relies on the analysis of the corresponding Riemann-Hilbert problem and the construction of an appropriate parametrix in a neighborhood of the origin. Obtained integral formulas are consistent with already known analogous results for Ablowitz-Segur solutions of homogeneous Painlev\'e II equation.
\end{abstract}


\maketitle

\setcounter{tocdepth}{2}

\section{Introduction}

We are concerned with the inhomogeneous second Painlev\'e (PII) equation 
\begin{align}\label{PII}
u''(x) = x u(x) + 2u^{3}(x) - \alpha,\quad x\in\C,
\end{align}
where $\alpha\in\C$. The solutions of the equation \eqref{PII} are related with Riemann-Hilbert (RH) problems characterized by {\em Stokes multipliers}, that is, the triple of parameters $(s_{1},s_{2},s_{3})\in\C^{3}$ satisfying the constraint condition 
\begin{equation}\label{stokes2}
s_{1}-s_{2}+s_{3} + s_{1}s_{2}s_{3} = -2\sin(\pi\alpha).
\end{equation}
Roughly speaking, any choice of $(s_{1},s_{2},s_{3})\in\C^{3}$ satisfying the condition \eqref{stokes2}, gives us a solution $\Phi(\lambda,x)$ of the corresponding RH problem, which is a $2\times 2$ matrix valued function sectionally holomorphic in $\lambda$ and meromorphic with respect to the variable  $x$ (see \cite{MR0588248}, \cite{MR0723758}, \cite{MR2264522}, \cite{MR0851569}, \cite{jimbo}). If we assume that $\theta(\lambda,x) := i(4\lambda^{3}/3 + x\lambda)$ is a phase function and $\sigma_{3}$ is the third Pauli matrix, then the function $u(x)$ obtained by the following limit 
\begin{align*}
u(x)=\lim_{\lambda\to\infty}(2\lambda\Phi(\lambda,x)e^{\theta(\lambda,x)\sigma_{3}})_{12},
\end{align*}
is a solution of the PII equation \eqref{PII}. Proceeding in this way, we can define a map 
\begin{align*}
\{(s_{1},s_{2},s_{3})\in\C^{3} \text{ satisfying }\eqref{stokes2}\}\to \{\text{solutions of the PII}\},
\end{align*}
which is a bijection between the set of all Stokes multipliers and the set of solutions of the Painlev\'e II equation (see e.g. \cite{MR2264522}). Let us restrict our attention to the solutions of \eqref{PII} corresponding to the following choice of the Stokes data
\begin{equation}\label{stokes-11}
s_1 = -\sin(\pi\alpha) -i k, \quad s_2=0,\quad s_3 = -\sin(\pi\alpha) + i k,\quad k,\alpha\in\C,
\end{equation}
that, for the brevity, we denote by $u(\,\cdot\,;\alpha,k)$. Among them we can specify the real Ablowitz-Segur (AS) solutions that are determined by \eqref{stokes-11} with 
\begin{equation}
\begin{gathered}\label{stokes2-ras}
\alpha\in(-1/2,1/2),\quad k\in(-\cos(\pi\alpha),\cos(\pi\alpha)).
\end{gathered}
\end{equation}
and the Hasting-McLeod (HM) solutions that correspond to the borderline case 
\begin{equation*}
\alpha\in (-1/2,1/2),\quad k=\pm\cos(\pi\alpha).
\end{equation*}
The real Ablowitz-Segur and Hasting-McLeod solutions satisfy $\mathrm{Im}\, u(x;\alpha,k) = 0$ for $x\in\R$ (see e.g. \cite[Chapter 11]{MR2264522}) and the results of \cite{MR2434886} and \cite{MR3670014} show that they have no poles on the real line. 
Considering the multipliers \eqref{stokes-11} with $\alpha,k\in i\R$ we obtain the purely imaginary Ablowitz-Segur solutions satisfying the relation $u(x;\alpha,k)\in i\R$ for $x\in\R$ (see e.g. \cite[Chapter 11]{MR2264522}). In this case $u$ also does not admits poles on the real axis, which follows from the fact that the residues of the poles of $u$ are equal to $\pm 1$ (see e.g. \cite{MR1960811}). All above kinds of the PII transcendents are important due to their applications in mathematics and physics. For example the Hasting-McLeod solutions appear in the random matrix theory and they are related with the study of the asymptotics expansions of distribution functions (see \cite{MR2395479}, \cite{MR3781450} and \cite{dai}). Furthermore we refer the reader to \cite{MR3660311} and \cite{MR3740377} for the recent applications of these solutions in the theory of liquid crystals. See also \cite{MR0481656} and \cite{MR1207209}, where the Ablowitz-Segur solutions are considered as the self-similar profiles of the modified Korteweg-de Vries equation. The aim of this paper is to study the Cauchy integrals
\begin{align}\label{t-int}
\int_{-\infty}^{\infty}u(y)\,dy :=\lim_{x\to +\infty} \int_{-x}^{x}u(y)\,dy,
\end{align}
where $u$ is a solution of the second  Painlev\'e equation. The problem of finding the values of \eqref{t-int} was considered in \cite[Theorem 2.1 and Theorem 3.1]{MR2501035}, where the following formula was derived for the real and purely imaginary Ablowitz-Segur solutions of the homogeneous PII equation 
\begin{align}\label{c-int-1}
\int_{-\infty}^{+\infty} u(y;0,k)\,dy = \frac{1}{2}\ln\left(\frac{1+k}{1-k}\right) \quad \text{if either} \quad k\in(-1,1) \quad \text{or} \quad k\in i\R.
\end{align}
The equality \eqref{c-int-1} was later extended in \cite[Lemma 7.1]{MR3781450} to the following asymptotic relation, which is valid for the real Ablowitz-Segur solutions ($k\in(-1,1)$)
\begin{align}\label{tot-int-asym}
\int_{s}^{+\infty} u(y;0,k)\,dy = \frac{1}{2}\ln\left(\frac{1+k}{1-k}\right) + O((-s)^{-3/4}),\quad s\to-\infty.
\end{align}
In \cite[Theorem 2.2]{MR2501035} the following total integral formula was obtained for the Hasting-McLeod solutions of the homogeneous PII equation ($\alpha=0$ and $k=\pm1$)
\begin{align}\label{int-hm}
\int_{c}^{\infty} u(y;0,\pm 1)\,dy + \int_{-\infty}^{c}\left(u(y;0,\pm 1) \!\mp\!\sqrt{\frac{|y|}{2}}\right)\,dy = \mp\frac{\sqrt{2}}{3}c|c|^{\frac{1}{2}} \pm \frac{1}{2}\log(2).
\end{align}
In this case $c$ is an arbitrary real number and the form of the left hand side of \eqref{int-hm} is a consequence of the fact that the solution $u(x;0,\pm 1)$ decays exponentially to zero as $x\to+\infty$ and 
\begin{align*}
u(x;0,\pm 1) = \pm\sqrt{-x/2} + O((-x)^{-5/2}),\quad x\to-\infty
\end{align*}
See \cite{MR0555581} and \cite{MR1322812} for the formal derivation and the rigorous proof of these asymptotics, respectively. 
Recently, in \cite[Theorem 1.1]{dai}, the value of the total integral was established for the solutions of inhomogeneous Painlev\'e II equation that are related to the following choice of the monodromy data 
\begin{align}\label{tronq}
s_{1} = e^{-i\pi(\alpha+\frac{1}{2})}, \qquad s_{2} = \omega, \qquad s_{3} = e^{i\pi(\alpha+\frac{1}{2})}
\end{align}
where $\alpha>-1/2$ and $\omega\ge 0$ are parameters. These PII transcendents are examples the classical tronque\'e solutions introduced in \cite{boutr} (see e.g. \cite{MR2264522} and \cite{kapaev} for more details). Observe that, if $\omega = 0$ then the Stokes multipliers \eqref{tronq} determine the Hasting-McLeod solution of the inhomogeneous PII equation with $\alpha\in (-1/2,1/2)$ and $k=\cos(\pi\alpha)$. In particular \cite{dai} reproves the formula \eqref{int-hm} for the solution $u(\,\cdot\,;0, 1)$. 
Furthermore we refer the reader to \cite[Theorem 1.3]{miller}, where the total integral formula was established for the increasing tritronqu\'ee solutions (see \cite{Joshi}) of the PII equation.
In this paper we intend to prove the following theorem, which provides a formula expressing the value of the Cauchy integral \eqref{t-int} for the real Ablowitz-Segur solutions of the inhomogeneous second Painlev\'e equation in the terms of the parameters $\alpha$ and $k$. 
\begin{theorem}\label{th-total}
If $u(\,\cdot\,;\alpha,k)$ is a real Ablowitz-Segur solution for the inhomogeneous second Painlev\'e equation, then 
\begin{align}\label{wz1}
&\lim_{x\to+\infty} \int_{-x}^x u(y;\alpha,k)\,dy = \frac{1}{2}\ln\left(\frac{\cos(\pi\alpha)+k}{\cos(\pi\alpha)-k}\right).
\end{align}
\end{theorem}
Furthermore we will prove a corresponding result for the purely imaginary Ablowitz-Segur solutions in the following exponential form.
\begin{theorem}\label{th-total-2}
If $u(\,\cdot\,;\alpha,k)$ is a purely imaginary Ablowitz-Segur solution for the inhomogeneous second Painlev\'e equation, then 
\begin{align}\label{wz2}
&\lim_{x\to+\infty} \exp\left(\int_{-x}^x u(y;\alpha,k)\,dy\right)=\frac{\cos(\pi\alpha)+k}{\left(\cos^2(\pi\alpha)-k^2\right)^{1/2}}.
\end{align}
\end{theorem}
Let us observe that the obtained integral formulas \eqref{wz1} and \eqref{wz2} are consistent with \eqref{c-int-1} in the case of homogeneous PII equation ($\alpha=0$).
Before the proofs of Theorems \ref{th-total} and \ref{th-total-2} we recall the steepest descent deformations of the the RH problem contour of $\Phi(\lambda)$, that were used in \cite{MR3670014}, \cite{MR2264522} and \cite{MR1950792} to the study of the asymptotic behavior of the solution $u(x;\alpha,k)$ as $x\to\pm\infty$. Depending on the sign of the parameter $x\in\R$, we consider the equivalent Riemann-Hilbert problems defined on the deformed graphs and obtain representations of their solutions in the terms of appropriate parametrices and component functions. The main difficulty follows from the fact that the case of $x<0$ requires us to find a new parametrix describing singularity appearing at the origin of the complex plane. The  existence of the parametrix was recently proved in \cite[Section 3.5]{MR3670014} by the application of a vanishing lemma, however its exact form remains unknown. In this paper we provide an explicit formula for the local parametrix in the terms of the classical Bessel functions (see Theorem \ref{lok-aprox}) and use it to establish asymptotics of appropriate functions related with the solutions of the deformed RH problems. These in turn will provide the total integral formulas \eqref{wz1} and \eqref{wz2}. \\[5pt]
{\bf Outline.} In Section 2 we recall the Riemann-Hilbert problem for the inhomogeneous PII equation and analyze the related Flaschka-Newell Lax system to reduce the proof of the formulas \eqref{wz1} and \eqref{wz2} to finding appropriate asymptotics of functions related with the solution $\Phi(\lambda,x)$. In Section 3 we introduce an auxiliary RH problem that will be used in the construction of a local parametrix. Sections 4 and 5 are devoted to representations of solutions for deformed RH problems and asymptotic estimates of their component functions. Finally, in Section 6 we provide the proof of Theorems \ref{th-total} and \ref{th-total-2}. \\[5pt]
\noindent{\bf Notation and terminology.} We write $M_{2\times 2}(\C)$ for the linear space of $2\times 2$ matrices with complex coefficients, equipped with the Frobenius norm 
\begin{align*}
\|A\|:=\sqrt{|a_{11}|^{2} + |a_{12}|^{2} + |a_{21}|^{2} + |a_{22}|^{2}},\quad A=[a_{lm}]\in M_{2\times 2}(\C).
\end{align*}
It is known that the norm is sub-multiplicative, that is, 
\begin{equation*}
\|AB\|\le \|A\| \|B\|,\quad A,B\in M_{2\times 2}(\C).
\end{equation*}
Let $\Gamma$ be a contour contained in the complex $z$-plane, which is the sum of a finite number of possibly unbounded oriented curves that are smooth in the Riemann sphere. Let us assume that the set $S$, consisting of the intersection points of these curves, has a finite number of elements and the complement $\C\setminus\Gamma$ has a finite number of connected components. The graph $\Gamma$ has the natural orientation determined by the orientations of its component curves. Therefore, for any point of the set $\Gamma\setminus S$, we can naturally define the (+) and (--) sides of $\Gamma$. Given $1\le p<\infty$, we denote by $L^{p}(\Gamma)$ the space consisting of functions $f:\Gamma\to M_{2\times 2}(\C)$ with the property that
\begin{align*}
\|f\|^{p}_{L^{p}(\Gamma)}:= \int_{\Gamma}\|f(z)\|^{p}\,|dz|<\infty.
\end{align*}
Furthermore we write $L^{\infty}(\Gamma)$ for the space of functions $f$ such that 
\begin{align*}
\|f\|_{L^{\infty}(\Gamma)}:= \mathrm{ess\,sup}_{z\in \Gamma}\,\|f(z)\|<\infty.
\end{align*}
If $1\le p<\infty$ and the graph $\Gamma$ is unbounded, then we consider the space $L^{p}_{I}(\Gamma)$, whose elements are functions $f:\Gamma\to M_{2\times 2}(\C)$ with the property that there is $f_{\infty}\in M_{2\times 2}(\C)$ such that $f-f_{\infty}\in L^{p}(\Gamma)$ (see \cite{zhou}). It is not difficult to check that the matrix $f_{\infty}$ is uniquely determined by $f$, which allows us to set the norm 
\begin{align*}
\|f\|^{p}_{L^{p}_{I}(\Gamma)}:=\|f-f_{\infty}\|^{p}_{L^{p}(\Gamma)} + \|f_{\infty}\|^{p},\quad f\in L^{p}_{I}(\Gamma).
\end{align*}
Throughout  the paper we use the notation $\mathcal{C}_\pm$ for the Cauchy operator defined on the contour $\Gamma$, which, for any $f\in L^{p}(\Gamma)$ with $1\le p<\infty$, is given by
\begin{align*}
[\CC_\pm f](z) := \lim_{z'\to z^{\pm}}\frac{1}{2\pi i}\int_{\Gamma} \frac{f(\xi)}{\xi-z'}\,d\xi, \quad z\in \Gamma.
\end{align*}
In the above limit the argument $z'$ tends non-tangentially to $z$ from the $(\pm)$-side of the graph $\Gamma$. We will also frequently write $A\lesssim B$ to denote $A\le CB$, for some $C>0$. Then the notation $A\sim B$ means that there are constants $C_{1},C_{2}>0$ such that $C_{1} B \!\le\! A\le C_{2} B$. \\[5pt]
\noindent {\bf Acknowledgements.} The study of the author are supported by the MNiSW Iuventus Plus Grant no. 0338/IP3/2016/74.



\section{Solution $\Phi(\lambda,x)$ of the RH problem for the PII equation}\label{sec-abl-seg}

\begin{figure}[h]
\begin{subfigure}{0.49\textwidth}
\centering
\includegraphics[scale=0.65]{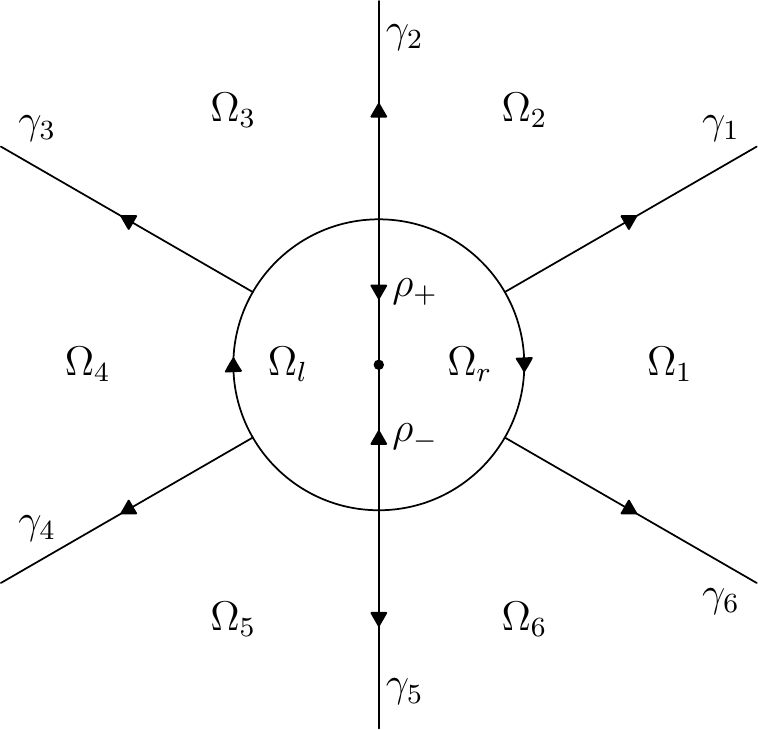}
\end{subfigure} 
\begin{subfigure}{0.49\textwidth}
\centering
\includegraphics[scale=0.7]{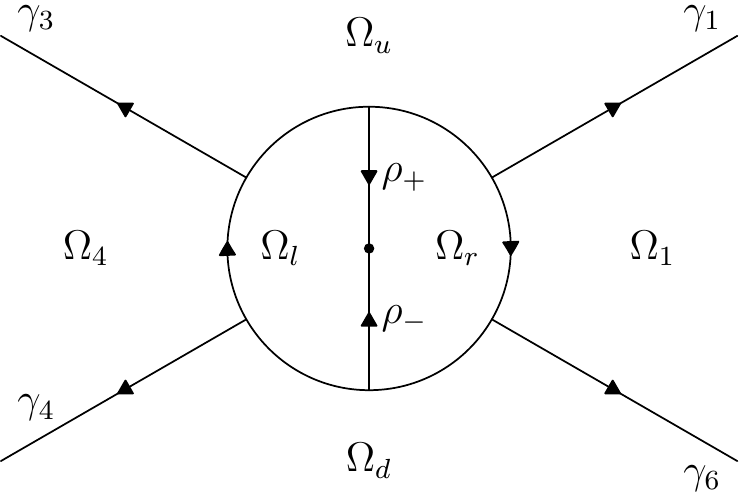}
\end{subfigure} 
\caption{Left: contour $\Sigma$ and regions of the set $\C\setminus\Sigma$. Right: contour $\Sigma$ in the case $s_{2}=0$.}
\label{fig:10}
\end{figure}

Assume that $\alpha\in\C$ and $\Sigma:=C\cup\rho_+\cup\rho_-\cup\gamma_{1}\cup\ldots\cup\gamma_{6}$ is the contour in the complex $\lambda$-plane, where $C:=\{\lambda\in\C \ | \ |\lambda|=r\}$ for $r>0$, is a clockwise oriented circle, $\rho_\pm:=\{\lambda\in\C \ | \ |\lambda|<r, \, \mathrm{arg}\,\lambda = \pm\frac{\pi}{2}\}$ are two radii oriented to the origin and $$\gamma_{k}:=\{\lambda\in\C \ | \ |\lambda|>r, \ \mathrm{arg}\, \lambda = \pi/6 + (k-1)\pi/3\},\quad k=1,2,\ldots,6$$ are six rays oriented to the infinity. The contour $\Sigma$ divides the complex plane into regions $\Omega_r$, $\Omega_l$ and $\Omega_k$ for $k=1,\ldots,6$, as it is depicted in the left diagram of Figure \ref{fig:10}.
Let us consider the triangular matrices
$$S_{k}:=\begin{pmatrix}1& 0\\ s_{k} & 1\end{pmatrix}, \ \ k=1,3,5 \quad\text{and}\quad
S_{k}:=\begin{pmatrix}1& s_{k}\\ 0 & 1\end{pmatrix}, \ \ k=2,4,6,$$
where constants $s_{k}$, for $k=1,\ldots,6$, satisfy the constrain condition
\begin{align}\label{stokes-mult}
s_{k+3}=-s_{k},\quad s_{1}-s_{2}+s_{3} + s_{1}s_{2}s_{3} = -2\sin(\pi\alpha).
\end{align}
Moreover we assume that $\sigma_1$, $\sigma_2$ and $\sigma_3$ are the usual Pauli matrices
\begin{align*}
\sigma_{1}:=\begin{pmatrix}
0& 1\\
1& 0
\end{pmatrix},\qquad
\sigma_{2}:=\begin{pmatrix}
0& -i\\
i& 0
\end{pmatrix},\qquad
\sigma_{3}:=\begin{pmatrix}
1& 0\\
0&-1
\end{pmatrix}
\end{align*}
and $E$ is a unimodular connection matrix such that the following equality holds
\begin{align}\label{equa-monodromy}
ES_1S_2S_3 = \sigma_2 M^{-1}E\sigma_2,\quad\text{where} \ \ M:=-ie^{i\pi(\alpha-\frac{1}{2})\sigma_3}\sigma_2.
\end{align}
The Riemann-Hilbert problem associated with the second Painlev\'e equation is to find a function $\Phi(\lambda)=\Phi(\lambda,x)$ with values in the space $M_{2\times 2}(\C)$ such that the following conditions are satisfied. \\[2pt]
\noindent\makebox[8mm][l]{(A1)}\parbox[t][][t]{119mm}{For any $k=1,2,\ldots,6$, the restriction $\Phi_k:=\Phi_{|\Omega_k}$ is holomorphic on $\Omega_k$ and continuous up to the boundaries of $\Omega_k$.}\\[3pt]
\noindent\makebox[8.5mm][l]{$(A2)$}\parbox[t][][t]{119mm}{Given $k\in\{l,r\}$, the restrictions $\Phi_{k}:=\Phi_{|\Omega_{k}}$ is holomorphic on $\Omega_{k}$ and $\Phi_{k}\in C(\overline{\Omega_{k}^\ve})$ for sufficiently small $\ve>0$, where $\Omega^\ve_{k}:=\Omega_k\setminus\{\lambda\in\C \ | \ |\lambda|<\ve\}$.}\\[3pt]
\noindent\makebox[8mm][l]{(A3)}\parbox[t][][t]{119mm}{Given $\lambda\in\Sigma\setminus\{0\}$, if we denote by $\Phi^+(\lambda)$ and $\Phi^-(\lambda)$ the limits of $\Phi(\lambda')$ as $\lambda'\to\lambda$ from the left (+) and right (--) side of the contour $\Sigma$, respectively, then the following {\em jump condition} is satisfied
\begin{equation*}
\Phi^+(\lambda) = \Phi^-(\lambda)S(\lambda),\quad \lambda\in\Sigma\setminus\{0\},
\end{equation*}}\\
\noindent\makebox[8mm][l]{}\parbox[t][][t]{119mm}{where the jump matrix $S(\lambda)$ is constructed as follows. On the rays $\gamma_k$ the matrix $S(\lambda)$ is given by the equation
\begin{equation*}
S(\lambda) = S_k, \quad \lambda\in\gamma_k, \quad k=1,2,\ldots,6,
\end{equation*}}\\
\noindent\makebox[8mm][l]{}\parbox[t][][t]{119mm}{while on the circle $C$ the matrix $S(\lambda)$ is obtained by the following relations
\begin{equation*}
\begin{aligned}
&\Phi_1^+(\lambda)=\Phi_r^-(\lambda)E, && \Phi_2^+(\lambda)=\Phi_r^-(\lambda)ES_1, \\
&\Phi_3^+(\lambda)=\Phi_l^-(\lambda)\sigma_2E\sigma_2S_3^{-1}, && \Phi_4^+(\lambda)=\Phi_l^-(\lambda)\sigma_2E\sigma_2,\\
&\Phi_5^+(\lambda)=\Phi_l^-(\lambda)\sigma_2ES_1\sigma_2, && \Phi_6^-(\lambda)=\Phi_r^-(\lambda)ES_6^{-1}.
\end{aligned}
\end{equation*}}\\
\noindent\makebox[8mm][l]{}\parbox[t][][t]{119mm}{Furthermore, on the radii $\rho_\pm$, the matrix $S(\lambda)$ is determined by the equations
\begin{equation}
\begin{aligned}\label{jump-circ}
&\Phi_l^-(\lambda)=\Phi_r^+(\lambda)M, && \quad\lambda\in \rho_-\setminus\{0\},\\
&\Phi_r^+(\lambda)=\Phi_l^-(\lambda)\sigma_2M\sigma_2, && \quad\lambda\in \rho_+\setminus\{0\}.
\end{aligned}
\end{equation}}\\
\noindent\makebox[8mm][l]{(A4)}\parbox[t][][t]{119mm}{The function $\Phi_r(\lambda)\lambda^{-\alpha\sigma_3}$ is bounded for $\lambda$ sufficiently close to zero, where the branch $\lambda^{-\alpha}$ is chosen arbitrarily.}\\[3pt]
\noindent\makebox[8mm][l]{(A5)}\parbox[t][][t]{119mm}{The function $\Phi$ has the following asymptotic behavior
$$\Phi(\lambda) = (I + O(\lambda^{-1}))e^{-\theta(\lambda)\sigma_{3}},\quad\lambda\to\infty,$$
where the phase function is given by $\theta(\lambda,x) := i(\frac{4}{3}\lambda^{3} + x\lambda)$.}\\[5pt]

In the following lemma we obtain useful asymptotic relations of the function $\Phi(\lambda)$ at the origin of the complex plane.
\begin{lemma}\label{lem-asy-phi}
If $0<\mathrm{Re}\,\alpha<1/2$ then the function $\Phi(\lambda)$ has the following behavior
\begin{align}\label{aas-1}
\Phi(\lambda)= O\begin{pmatrix}|\lambda|^{-\alpha}& |\lambda|^{-\alpha}\\[5pt] |\lambda|^{-\alpha}& |\lambda|^{-\alpha}\end{pmatrix},\quad \lambda\to 0
\end{align}
and furthermore, if $1/2<\mathrm{Re}\,\alpha\le 0$ then 
\begin{align}\label{aas-2}
\Phi(\lambda)= O\begin{pmatrix}|\lambda|^{\alpha}&|\lambda|^{\alpha}\\[5pt]|\lambda|^{\alpha}& |\lambda|^{\alpha}\end{pmatrix},\quad\lambda\to 0.
\end{align}
\end{lemma}
\proof If we define the following functions
\begin{equation*}
\begin{aligned}
A_{1}(\lambda):=\Phi(\lambda)\sigma_{2}M\sigma_{2}\lambda^{-\alpha\sigma_3},\quad \lambda\in\Omega_{l},\quad 
A_{2}(\lambda):=\Phi(\lambda)\lambda^{-\alpha\sigma_3},\quad \lambda\in\Omega_{r}.\\
\end{aligned}	
\end{equation*}
then we have the representations
\begin{equation*}
\begin{aligned}
\Phi(\lambda) = A_{1}(\lambda)\lambda^{\alpha\sigma_{3}}\sigma_{2}M^{-1}\sigma_{2},\quad \lambda\in\Omega_{l}, \quad 
\Phi(\lambda) = A_{2}(\lambda)\lambda^{\alpha\sigma_{3}},\quad \lambda\in\Omega_{r}. 
\end{aligned}
\end{equation*}
Since the condition (A5) is satisfied, it follows that the function $A_{2}(\lambda)$ is bounded whenever $\lambda\in\Omega_{r}$ is sufficiently close to zero. Using the jump relations \eqref{jump-circ}, we deduce that the function $A_{1}(\lambda)$ is a holomorphic extension of $A_{2}(\lambda)$ over the set $\Omega_{l}$ and therefore, using the condition (A5) once again, we infer that the function $A_{1}(z)$ is also bounded provided $\lambda\in\Omega_{r}$ is sufficiently close to the origin. This in turn gives us the asymptotics \eqref{aas-1} and \eqref{aas-2}. \hfill $\square$ 
\begin{remark}
In view of the choice of the monodromy data \eqref{stokes-11}, we have $S_{2}=S_{5}=I$ and therefore the diagram $\Sigma$ takes the form depicted on the right diagram of Figure \ref{fig:10}, where we denote 
$\Omega_{u}:=\Omega_{2}\cup\Omega_{3}$ and $\Omega_{d}:=\Omega_{5}\cup\Omega_{6}$. By the condition \eqref{stokes-mult} the remaining Stokes multipliers $s_{1}$, $s_{3}$ satisfy the equality 
\begin{equation}\label{cc1}
s_{1}+s_{3} = -2\sin(\pi\alpha).
\end{equation}
Furthermore direct calculations using the equation \eqref{equa-monodromy} show that the connection matrix $E$ takes the following form (see \cite[Chapter 11.6]{MR2264522})
\begin{align}\label{eq-e}
E=\begin{pmatrix}p&0\\[5pt] 0&q\end{pmatrix}\begin{pmatrix}1&ie^{-i\pi\alpha}\\[5pt] 1&-ie^{i\pi\alpha}\end{pmatrix},
\end{align}
where $p=1$ and $q=-(2i\cos(\pi\alpha))^{-1}$. \hfill $\square$
\end{remark}
From \cite{MR1677884}, \cite{MR0588248}, \cite[Chapter 11]{MR2264522} and \cite{jimbo} we know that the RH problem associated with the PII equation is uniquely meromorphically (with respect to $x$) solvable for any choice of the Stokes multipliers $(s_1,s_2,s_3)$. Therefore there is a function $\Phi(\lambda,x)$ and the countable set of complex numbers $X=\{x_k\}_{k\ge 1}$ such that, for any $x\in \C\setminus X$, the function $\Phi(\lambda,x)$ satisfies the conditions (A1)\,--\,(A5) and $\Phi(\lambda,x)$ is holomorphic on $(\C\setminus\Sigma)\times(\C\setminus X)$. Furthermore, for any $k\ge 1$, the point $x_{k}$ is a pole of the function $\Phi(\lambda,x)$ such that the coefficients of the corresponding Laurent series are functions of $\lambda$.
It is also known that the function $\Phi(\lambda,x)$ satisfies the following Flaschka-Newell Lax pair
\begin{equation}\label{lax-sys}
\left\{\begin{aligned}
\partial_{\lambda}\Phi(\lambda,x) = \AA(\lambda,x)\Phi(\lambda,x),\\[1pt]
\partial_{x}\Phi(\lambda,x) = \UU(\lambda,x)\Phi(\lambda,x).
\end{aligned}\right.
\end{equation}
In the above system $\AA$ and $\UU$ are $2\times 2$ matrix functions that are given by
\begin{align*}
\AA(\lambda,x)&:= -i(4\lambda^{2} + x + 2u^{2}(x))\sigma_{3} - (4\lambda u(x) + \alpha\lambda^{-1})\sigma_{2} - 2u_{x}(x)\sigma_{1},\\
\UU(\lambda,x)&:= -i\lambda\sigma_{3} - u(x)\sigma_{2}.
\end{align*}
By \cite[Chapter 11.3]{MR2264522} we know that the solution $\Phi$ has the following representation
\begin{align}\label{eq-zet-phi}
\Phi(\lambda) = Z(\lambda)e^{-\theta(\lambda)\sigma_{3}}\lambda^{\alpha\sigma_3},\quad \lambda\in \Omega_{r},
\end{align}
where the function $Z(\lambda)$ holomorphic on the $B(0,r)$ and the branch cut of the logarithm is taken such that $\mathrm{arg}\,\lambda\in(-\pi/2,3\pi/2)$. Using the second equation of the system \eqref{lax-sys}, we infer that
\begin{align*}
\frac{\partial Z(\lambda)}{\partial x}&= \frac{\partial \Phi(\lambda)}{\partial x}e^{\theta(\lambda)\sigma_{3}}\lambda^{-\alpha\sigma_3} +
i\lambda\Phi(\lambda)\sigma_3e^{\theta(\lambda)\sigma_{3}}\lambda^{-\alpha\sigma_3}  \\
&= (-i\lambda\sigma_{3} - u(x)\sigma_{2})\Phi(\lambda)e^{\theta(\lambda)\sigma_{3}}\lambda^{-\alpha\sigma_3} +
i\lambda\Phi(\lambda)\sigma_3e^{\theta(\lambda)\sigma_{3}}\lambda^{-\alpha\sigma_3} \\
&= (-i\lambda\sigma_{3} - u(x)\sigma_{2})Z(\lambda) +  i\lambda Z(\lambda)\sigma_3 = -i\lambda [\sigma_3,Z(\lambda)] -u(x)\sigma_{2}Z(\lambda).
\end{align*}
Therefore, if we define $P(x):=Z(0,x)$, then, passing to the limit with $\lambda\to 0$, gives the following linear equation
\begin{align}\label{eq-p}
\frac{\partial P(x)}{\partial x}= -u(x)\sigma_{2}P(x),\quad x\in\C\setminus X.
\end{align}
As it was proved in \cite[Section 2.2]{MR3670014}, if $u$ is either a real or purely imaginary AS solution, then the set of poles $X$ does not contain any real number and the function $\Phi(\lambda,x)$ is defined for all $x\in\R$. Therefore, if we denote
\begin{align*}
v(x_{1},x_{2}):=\int_{x_{1}}^{x_{2}}u(y)\,dy, \quad x_{1}<x_{2},
\end{align*}
then the solution of the equation \eqref{eq-p} satisfies the following formula
\begin{align*}
P(x_{2}) = \exp\left(-v(x_{1},x_{2})\sigma_{2}\right)P(x_{1}),\quad x_{1}<x_{2}
\end{align*}
and therefore, for any $x>0$, we have
\begin{align}\label{eq-p-2}
P(x)P(-x)^{-1} = \begin{pmatrix}\frac{1}{2}(e^{v(-x,x)}+e^{-v(-x,x)}) & \frac{i}{2}(e^{v(-x,x)}-e^{-v(-x,x)}) \\[5pt] -\frac{i}{2}(e^{v(-x,x)}-e^{-v(-x,x)})& \frac{1}{2}(e^{v(-x,x)}+e^{-v(-x,x)})\end{pmatrix}.
\end{align}
Consequently, the proof of the integral formula \eqref{wz1} reduces to studying asymptotic behavior of the function $P(x)$ as $x\to\pm\infty$. 



\section{Auxiliary RH problem}

Let us consider the function $\hat\Phi^{0}(z)$, given by the formula
\begin{equation}\label{funct-def-2}
\hat\Phi^{0}(z) \!:=\! B(z)\!\begin{pmatrix}v_1(z) \hspace{-2mm}& v_2(z)\\[5pt] v_1'(z) \hspace{-2mm}& v_2'(z)\end{pmatrix}\!,\ \ \text{where}\ \
B(z) := \frac{1}{2}e^{-i\frac{\pi}{4}\sigma_3}\begin{pmatrix}1&1\\[5pt]-1&1\end{pmatrix}
\begin{pmatrix}1&0\\[5pt] -\frac{\alpha}{z} &1\end{pmatrix}
\end{equation}
and the functions $v_1$, $v_2$ are defined by 
\begin{align}\label{form-v1}
v_1(z)= \sum_{k=0}^\infty \frac{\Gamma(\alpha+\frac{1}{2})z^{\alpha+2k}}{4^k k!\Gamma(\alpha+\frac{1}{2}+k)}=
2^{\alpha-\frac{1}{2}}\Gamma(\alpha+\frac{1}{2})e^{i\frac{\pi}{2}(\alpha-\frac{1}{2})}z^{\frac{1}{2}}J_{\alpha-\frac{1}{2}}(e^{-i\frac{\pi}{2}}z)
\end{align}
and
\begin{align}\label{form-v2}
v_2(z) := \sum_{k=0}^\infty \frac{\Gamma(\frac{3}{2}-\alpha)z^{1-\alpha+2k}}{4^k k!\Gamma(\frac{3}{2}-\alpha+k)}=
2^{\frac{1}{2}-\alpha}\Gamma(\frac{3}{2}-\alpha)e^{i\frac{\pi}{2}(\frac{1}{2}-\alpha)}z^{\frac{1}{2}}J_{\frac{1}{2}-\alpha}(e^{-i\frac{\pi}{2}}z),
\end{align}
where $J_\nu(z)$ is the classical Bessel function defined on the universal covering of the punctured complex plane $\C\setminus\{0\}$. Then, it is known that the function $\hat\Phi^{0}(z) z^{-\alpha\sigma_3}$ is holomorphic on the complex plane (see e.g. \cite{ba-er}, \cite{MR2264522}). Let us define the matrices
\begin{equation}
\begin{gathered}\label{matrix-s}
\hat E = \frac{\sqrt{\pi}}{2\cos\pi\alpha} \begin{pmatrix}\frac{2^{1-\alpha}}{\Gamma(\frac{1}{2}+\alpha)}&0\\[5pt] 0 &
\frac{2^{\alpha}}{\Gamma(\frac{3}{2}-\alpha)}\end{pmatrix} e^{i\frac{\pi}{4}\sigma_3}
\begin{pmatrix}e^{-i\pi\alpha}&i\\[5pt] ie^{i\pi\alpha}&1\end{pmatrix},\\[5pt]
\hat S_1 = \begin{pmatrix}1 & 2\sin(\pi\alpha) \\[5pt] 0&1\end{pmatrix},\qquad
\hat S_2 = \begin{pmatrix}1 & 0 \\[5pt] -2\sin(\pi\alpha)&1\end{pmatrix}
\end{gathered}
\end{equation}
and consider the functions 
\begin{gather}\label{funct-def}
\hat\Phi^{1}(z):=\hat\Phi^{0}(z)\hat E,\quad \hat\Phi^{2}(z):=\hat\Phi^{1}(z)\hat S_1, \quad \hat\Phi^{3}(z):=\hat\Phi^{2}(z)\hat S_2.
\end{gather}
Then, for any $1\le k\le 3$, we have the asymptotic behavior
\begin{gather*}
\hat\Phi^k(z)=(I - \frac{i\alpha}{2z}\sigma_1 + O(\frac{1}{z^2}))e^{z\sigma_3},\quad z\to\infty, 
\end{gather*}
with $\mathrm{arg}\,z\in(\pi(k-3/2),\pi(k+1/2))$ and the following equality holds
\begin{align}\label{eq-dd}
\sigma_{2}\hat\Phi^{k+1}(e^{i\pi}z)\sigma_{2} = \hat\Phi^{k}(z),\quad k=1,2
\end{align}
(see e.g. \cite{ba-er}, \cite{MR2264522}). It is also not difficult to check the useful equality 
\begin{equation}\label{eekk-1}
\hat E\hat S_{1} = DE,
\end{equation}
where the matrix $E$ is defined by \eqref{eq-e} and 
\begin{align*}
D:= \frac{\sqrt{\pi}e^{i\frac{\pi}{4}}}{\cos\pi\alpha} \begin{pmatrix}\frac{2^{-\alpha}e^{-i\pi\alpha}}{\Gamma(1/2+\alpha)}&0\\[5pt] 0 &
\frac{-i2^{\alpha}\cos(\pi\alpha)e^{i\pi\alpha}}{\Gamma(3/2-\alpha)}\end{pmatrix}.
\end{align*}

In the following lemma we obtain the value of the holomorphic function $\hat\Phi^{0}(z) z^{-\alpha\sigma_3}$ at the origin of the complex plane.
\begin{lemma}\label{lem-lim-1}
The following convergence holds
\begin{align}\label{eq-ff1}
\lim_{z\to 0}\hat\Phi^{0}(z)z^{-\alpha\sigma_3}=\frac{1}{2}e^{-i\frac{\pi}{4}\sigma_3}\begin{pmatrix}1&1\\[5pt]-1&1\end{pmatrix}\begin{pmatrix}1&0\\[5pt] 0 &1-2\alpha\end{pmatrix}.
\end{align}
\end{lemma}
\proof From the definition of $\hat\Phi^{0}(z)$ it follows that
\begin{align*}
\hat\Phi^{0}(z)z^{-\alpha\sigma_3} = \frac{1}{2}e^{-i\frac{\pi}{4}\sigma_3}\begin{pmatrix}1&1\\[5pt]-1&1\end{pmatrix}
\begin{pmatrix}1&0\\[5pt] -\alpha/z&1\end{pmatrix}\begin{pmatrix}v_1(z) & v_2(z)\\[5pt] v_1'(z) & v_2'(z)\end{pmatrix}
\begin{pmatrix}z^{-\alpha}&0\\[5pt] 0&z^\alpha\end{pmatrix}.
\end{align*}
On the other hand, we have
\begin{equation*}
\begin{aligned}
C(z)&:=\begin{pmatrix}1&0\\[5pt] -\alpha/z&1\end{pmatrix}\begin{pmatrix}v_1(z) & v_2(z)\\[5pt] v_1'(z) & v_2'(z)\end{pmatrix}
\begin{pmatrix}z^{-\alpha}&0\\[5pt] 0&z^\alpha\end{pmatrix}\\
&=\begin{pmatrix}1&0\\[5pt] -\alpha/z&1\end{pmatrix} \begin{pmatrix}v_1(z)z^{-\alpha} & v_2(z)z^{\alpha}\\[5pt] v_1'(z)z^{-\alpha} & v_2'(z)z^{\alpha}\end{pmatrix}\\
& = \begin{pmatrix}v_1(z)z^{-\alpha} & v_2(z)z^{\alpha}\\[5pt] -\alpha v_1(z)z^{-\alpha-1}+ v_1'(z)z^{-\alpha} & -\alpha v_2(z)z^{\alpha-1}+v_2'(z)z^{\alpha}\end{pmatrix}.
\end{aligned}
\end{equation*}
Using the formula \eqref{form-v1}, we have
\begin{align*}
v_1'(z)z^{-\alpha}= \alpha z^{-1}\sum_{k=0}^\infty \frac{\Gamma(\alpha+\frac{1}{2})z^{2k}}{4^k k!\Gamma(\alpha+\frac{1}{2}+k)}+
\sum_{k=1}^\infty 2k\frac{\Gamma(\alpha+\frac{1}{2})z^{2k-1}}{4^k k!\Gamma(\alpha+\frac{1}{2}+k)}
\end{align*}
and consequently
\begin{align}\label{form-v3}
-\alpha v_1(z)z^{-\alpha-1}+ v_1'(z)z^{-\alpha}= \sum_{k=1}^\infty 2k\frac{\Gamma(\alpha+\frac{1}{2})z^{2k-1}}{4^k k!\Gamma(\alpha+\frac{1}{2}+k)}.
\end{align}
On the other hand, from the formula \eqref{form-v2}, it follows that 
\begin{align*}
v_2'(z)z^{\alpha}&=\sum_{k=0}^\infty (2k+1-\alpha)\frac{\Gamma(\frac{3}{2}-\alpha)z^{2k}}{4^k k!\Gamma(\frac{3}{2}-\alpha+k)},
\end{align*}
which gives
\begin{align}\label{form-v4}
-\alpha v_2(z)z^{\alpha-1}+v_2'(z)z^{\alpha}= \sum_{k=0}^\infty (2k+1-2\alpha)
\frac{\Gamma(\frac{3}{2}-\alpha)z^{2k}}{4^k k!\Gamma(\frac{3}{2}-\alpha+k)}.
\end{align}
Combining \eqref{form-v1}, \eqref{form-v2}, \eqref{form-v3} and \eqref{form-v4}, we infer that
\begin{align*}
\lim_{z\to 0} C(z) = \begin{pmatrix}1&0\\[5pt] 0 &1-2\alpha\end{pmatrix},
\end{align*}
which finally gives the convergence \eqref{eq-ff1} and the proof of lemma is completed. \hfill $\square$\\

\begin{figure}[h]
\includegraphics[scale=0.65]{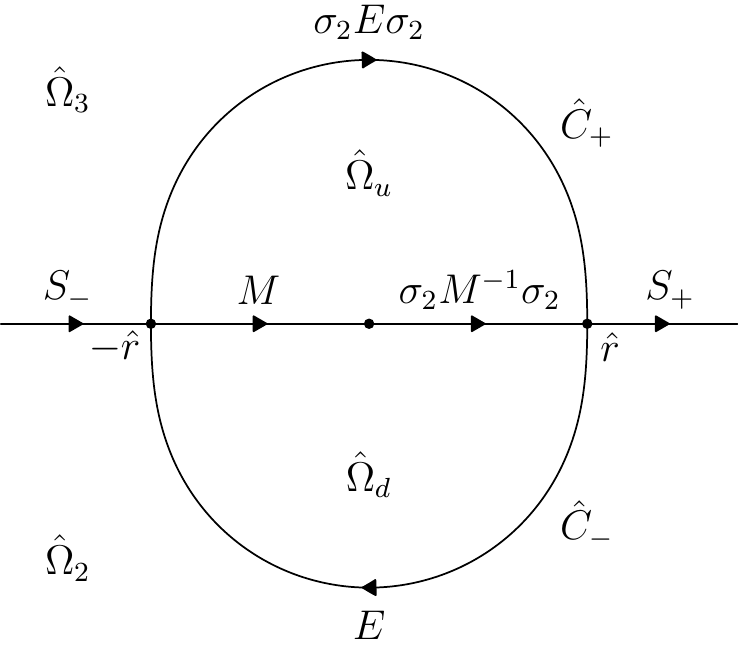}
\caption{The contour $\hat\Sigma$ for the auxiliary RH problem.}
\label{d6b}
\end{figure}

Given $x>0$, let us consider the function $z:\C\to\C$ given by the formula
\begin{equation}\label{map-z}
z(\lambda):=-\theta(\lambda) = -i(4\lambda^{3}/3 + x\lambda),\quad \lambda\in\C
\end{equation}
Let us assume that $r,R>0$ are taken such that 
\begin{align}\label{radius-r}
0<r<R:=\frac{1}{4}x^{\frac{1}{2}}.
\end{align}
It is not difficult to check that $z$ is an injective map on the open ball $B(0,R)$ and therefore, by the open mapping theorem for holomorphic functions, we infer that the set $V:=z(B(0,R))$ is open and the inverse $z^{-1}:V\to B(0,R)$ is also a holomorphic function. In the complex $z$-plane, we consider the contour $\hat\Sigma:=\R\cup \hat C_{+}\cup \hat C_{-}$, where $\hat C_{\pm}$ is the image of the set $\{\lambda\in \C \ | \ |\lambda| = r, \ \pm\mathrm{Re}\,\lambda\ge 0\}$ under the map $z$ (see Figure \ref{d6b}). If we put $\hat r := z(ir)$, then $\pm\hat r$ are the intersection points of $\hat C:=\hat C_{+}\cup \hat C_{-}$ with the real axis. 
The contour $\hat\Sigma$ divides the complex $z$-plane into four sets $\hat\Omega_d$, $\hat\Omega_u$, $\hat\Omega_2$, $\hat\Omega_3$ such that the sets $\hat\Omega_d$, $\hat\Omega_u$ lie 
in the interior of the circle $\hat C$ and the regions $\hat\Omega_2$, $\hat\Omega_3$ are located outside $\hat C$. We define the function $\hat\Phi(z)$ as follows
\begin{equation*}
\begin{gathered}
\hat\Phi(z)=\hat\Phi^{2}(e^{2\pi i}z), \ z\in\hat\Omega_2,\quad \hat\Phi(z)=\hat\Phi^{3}(e^{2\pi i}z), \ z\in\hat\Omega_3,\\
\hat\Phi(z)=\hat\Phi^{0}(e^{2\pi i}z)D, \ z\in\hat\Omega_d,\quad \hat\Phi(z)=\sigma_2\hat\Phi^{0}(e^{\pi i}z)D\sigma_2, \ z\in\hat\Omega_u,
\end{gathered}
\end{equation*}
where we recall that the maps $\hat\Phi^k(z)$, for $0\le k\le 3$, are defined on the universal covering of the punctured complex plane $\C\setminus\{0\}$ and the branch cut is chosen such that $\mathrm{arg}\,z\in(-\pi,\pi)$ (see \eqref{funct-def-2} and \eqref{funct-def}). 
\begin{lemma}\label{lem-eq}
For any $z\in\hat\Omega_{u}$ then following equality holds
\begin{align*}
\hat\Phi(z)\sigma_{2}E\sigma_{2}=\hat\Phi^{0}(e^{2\pi i}z)DE\hat S_{2}.
\end{align*}
\end{lemma}
\proof Combining \eqref{funct-def}, \eqref{eekk-1} and \eqref{eq-dd}, for any $z\in\hat\Omega_{u}$, gives
\begin{align*}
\hat\Phi(z)\sigma_{2}E\sigma_{2} & = \sigma_2\hat\Phi^{0}(e^{\pi i}z)DE\sigma_{2} = 
\sigma_2\hat\Phi^{0}(e^{\pi i}z)\hat E\hat S_{1}\sigma_{2} \\
& = \sigma_2\hat\Phi^{1}(e^{\pi i}z)\hat S_{1}\sigma_{2}=\sigma_2\hat\Phi^{2}(e^{\pi i}z)\sigma_{2}=\hat\Phi^{3}(e^{2\pi i}z)\\
& =\hat\Phi^{0}(e^{2\pi i}z)\hat E\hat S_{1}\hat S_{2}=\hat\Phi^{0}(e^{2\pi i}z)DE\hat S_{2}
\end{align*}
and the proof of lemma is completed. \hfill $\square$\\

The Riemann-Hilbert problem formulated in the following proposition will be used in the construction of a local parametrix for the steepest descent contour around the origin. For the proof we refer the reader to \cite[Section 11.6]{MR2264522}.

\begin{figure}[h]
\begin{subfigure}{0.49\textwidth}
\centering
\includegraphics[scale=0.65]{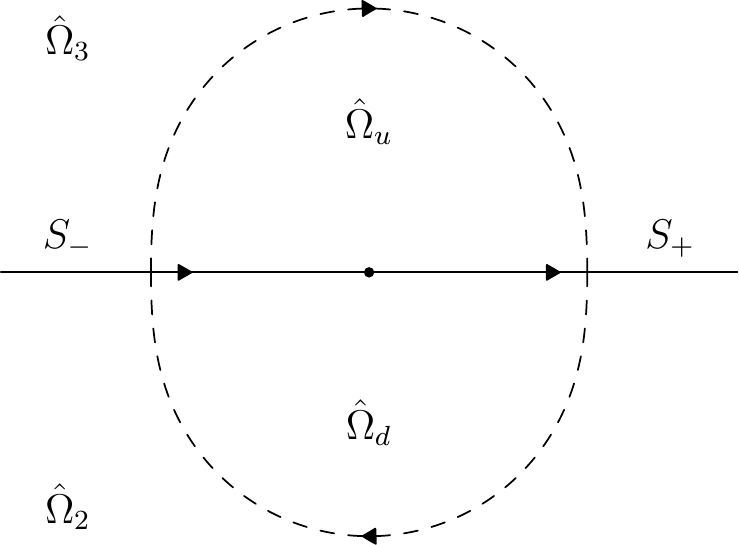}
\end{subfigure} 
\begin{subfigure}{0.49\textwidth}
\centering
\includegraphics[scale=0.65]{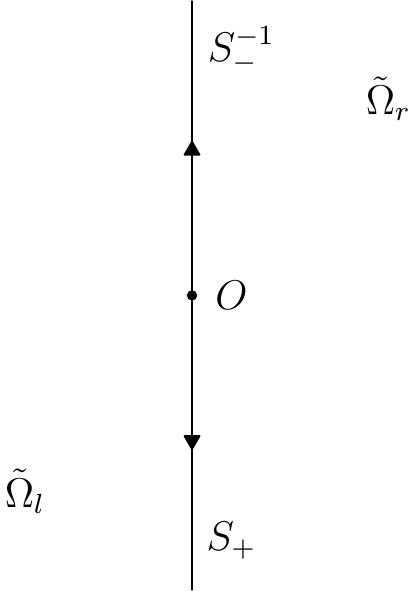}
\end{subfigure}
\caption{Left: the contour deformation between $\hat\Sigma$ and $\check\Sigma$. Right: the graph $\tilde\Sigma$ together 
with the jump matrices.}
\label{d17b}
\end{figure}

\begin{proposition}\label{prop-bessel}
The function $\hat\Phi(z)$ solves the following RH problem. \\
\noindent\makebox[6mm][l]{$(a)$}\parbox[t][][t]{121mm}{The function $\hat\Phi_{|\hat\Omega_d}(z)z^{-\alpha\sigma_3}$ is analytic on the open set confined by the curve $\hat C$.}\\
\noindent\makebox[6mm][l]{$(b)$}\parbox[t][][t]{121mm}{We have the jump relation $\hat\Phi_+(z) = \hat\Phi_-(z)\hat S(z)$ for $z\in\hat\Sigma$, where 
\begin{equation*}
\hat S(z):=\left\{\begin{aligned}
&S_+:=\hat S_2,&& \text{for}\quad z\in\R, \ z>\hat r,\\
&S_-:=\hat S_1^{-1},&& \text{for}\quad z\in\R, \ z<-\hat r\\
\end{aligned}\right.
\end{equation*}}
\noindent\makebox[6mm][l]{}\parbox[t][][t]{121mm}{and furthermore
\begin{gather*}
\hat S(z):= M, \ \ z\in\R, \ -\hat r<z<0,\quad \hat S(z):= \sigma_2 M^{-1}\sigma_2, \ \ z\in\R, \ 0<z<\hat r, \\
\hat S(z):= E, \ \ z\in\hat C_{-},\quad \hat S(z):=\sigma_2 E\sigma_2, \ \ z\in\hat C_{+}.
\end{gather*}}\\
\noindent\makebox[6mm][l]{$(c)$}\parbox[t][][t]{121mm}{The function $\hat\Phi(z)$ has the asymptotic behavior $$\hat\Phi(z) = (I + O(z^{-1}))e^{z\sigma_3},\quad z\to\infty.$$}
\end{proposition}
Let us assume that $\check\Sigma$ is an oriented contour in the complex plane consisting of two rays $\mathrm{arg}\,\lambda = 0$ and $\mathrm{arg}\,\lambda = \pi$ as it is shown on the left diagram of Figure \ref{d17b}. Let us consider the function $\check\Phi(z)$, given by the formulas
\begin{gather*}
\check\Phi(z) = \hat\Phi(z), \ \  z\in\hat\Omega_2\cup\hat\Omega_3,\quad \check\Phi(z)=\hat\Phi(z)\sigma_2 E\sigma_2, \ \ z\in\hat\Omega_u,\\
\check\Phi(z) = \hat\Phi(z)E, \ \ z\in\hat\Omega_d,
\end{gather*}
where the branch cut is chosen such that $\mathrm{arg}\,z\in(-\pi,\pi)$. 
\begin{proposition}\label{prop-l-h}
The function $\check \Phi(z)$ is a solution of the following RH problem. \\[2pt]
\noindent\makebox[5.5mm][l]{$(a)$}\parbox[t][][t]{121mm}{The function $\check\Phi(z)$ is holomorphic on $\C\setminus\check\Sigma$.}\\[3pt]
\noindent\makebox[5.5mm][l]{$(b)$}\parbox[t][][t]{121mm}{On the contour $\check\Sigma$, the function $\check\Phi(z)$ satisfies 
the following jump relations 
\begin{equation*}
\hspace{-30pt}\check\Phi_+(z) = \check\Phi_-(z)S_{+}, \ z\in\R, \, z>0 \ \text{ and } \ \check\Phi_+(z) = \check\Phi_-(z)S_{-}, \ z\in\R, \, z<0.
\end{equation*}}\\
\noindent\makebox[5.5mm][l]{$(c)$}\parbox[t][][t]{121mm}{If $0<\mathrm{Re}\,\alpha<1/2$ then the function $\check\Phi(z)$ satisfies the asymptotic relation
\begin{align}\label{asym-orig-1}
\check\Phi(z) = O\begin{pmatrix}|z|^{-\alpha} & |z|^{-\alpha} \\[5pt]
|z|^{-\alpha} & |z|^{-\alpha}\end{pmatrix},\quad z\to 0
\end{align}}\\
\noindent\makebox[5.5mm][l]{}\parbox[t][][t]{121mm}{and furthermore if $-1/2<\mathrm{Re}\,\alpha\le 0$ then
\begin{align}\label{asym-orig-2}
\check\Phi(z) = O\begin{pmatrix}|z|^{\alpha} & |z|^{\alpha} \\[5pt]
|z|^{\alpha} & |z|^{\alpha}\end{pmatrix},\quad z\to 0.
\end{align}}\\
\noindent\makebox[5.5mm][l]{$(d)$}\parbox[t][][t]{121mm}{We have the following asymptotic behavior 
\begin{equation*}
\check\Phi(z) =(I+O(z^{-1}))e^{z\sigma_3}, \quad z\to\infty.
\end{equation*}}
\end{proposition}
\noindent {\em Proof.} The conditions $(a)$, $(b)$ and $(d)$ are straightforward consequences of Proposition \ref{prop-bessel}. To show that $\check\Phi(z)$ satisfies also the point $(c)$, we define the functions
\begin{equation*}
\begin{aligned}
&A_{1}(z):=\hat\Phi(z)\sigma_{2}M\sigma_{2}z^{-\alpha\sigma_3},&& z\in \hat\Omega_{u},\\
&A_{2}(z):=\hat\Phi(z)z^{-\alpha\sigma_3},&& z\in \hat\Omega_{d}.\\
\end{aligned}	
\end{equation*}
If we write $B_{1}:=\sigma_{2}M^{-1}E\sigma_{2}$ and $B_{2}:= E$, then we have the representation
\begin{equation}\label{aa-mm-1}
\check\Phi(z) = A_{1}(z)z^{\alpha\sigma_{3}}B_{1}, \ z\in\hat\Omega_{u} \ \ \text{ and } \ \ \check\Phi(z) = A_{2}(z)z^{\alpha\sigma_{3}}B_{2}, \ z\in\hat\Omega_{d}. 
\end{equation}
By the point $(a)$ of Proposition \ref{prop-bessel}, the function $A(z)$ given by the formula $$A(z):=A_{1}(z), \ z\in\hat\Omega_{u},\quad 
A(z):=A_{2}(z), \ z\in\hat\Omega_{d}$$ is holomorphic in a neighborhood of the origin, which implies that 
\begin{equation}\label{aa-mm-2}
A_{1}(z) = O(1), \ z\to 0, \ z\in\hat\Omega_{u} \ \ \text{ and } \ \ A_{2}(z) = O(1), \ z\to 0, \ z\in\hat\Omega_{d}.
\end{equation}
Combining \eqref{aa-mm-1} and \eqref{aa-mm-2} gives the asymptotic relations \eqref{asym-orig-1} and \eqref{asym-orig-2}. Thus the proof of proposition is completed. \hfill $\square$ \\

\begin{figure}[h]
\begin{subfigure}{0.49\textwidth}
\centering
\hspace{10pt}\includegraphics[scale=0.6]{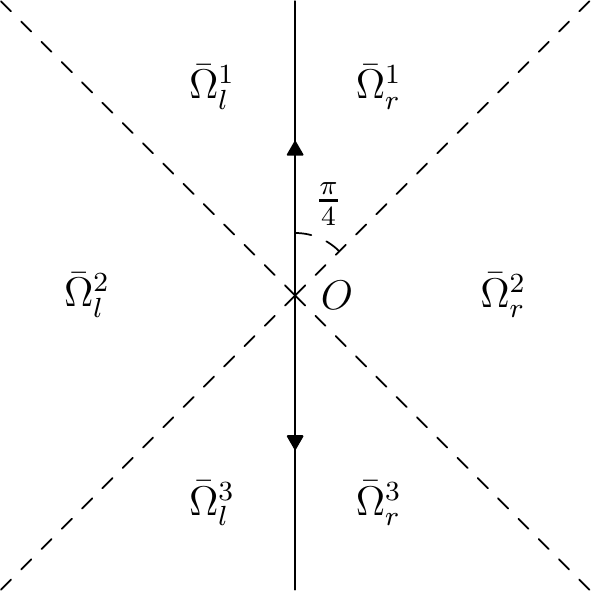}
\end{subfigure} 
\begin{subfigure}{0.49\textwidth}
\centering
\hspace{-10pt}\includegraphics[scale=0.65]{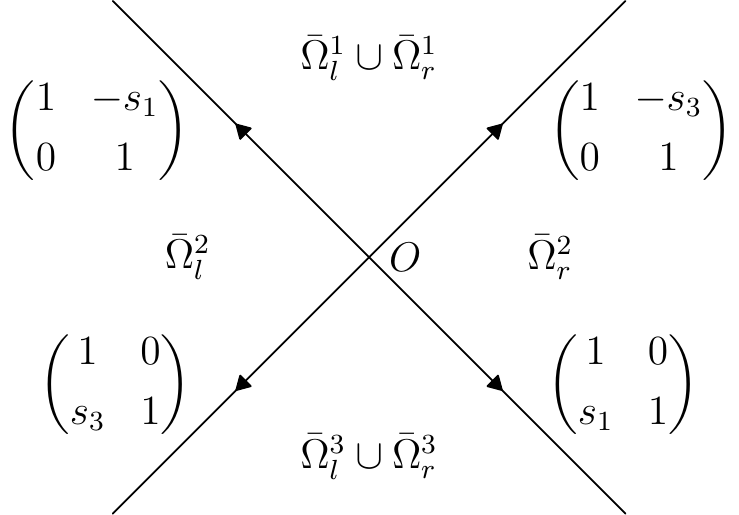}
\end{subfigure}
\caption{Left: a contour deformation between $\tilde\Sigma$ and $\bar\Sigma$. Right: the graph $\bar\Sigma$ with the associated jump matrices.}
\label{d5b}
\end{figure}

Let us assume that $\tilde\Sigma$ is an oriented contour in the complex plane consisting of two rays $\mathrm{arg}\,\lambda = \pm\pi/2$. The contour 
divides the complex plane into the sets $\tilde\Omega_{r}:=\{\mathrm{Re}\, z >0\}$ and $\tilde\Omega_{l}:=\{\mathrm{Re}\, z <0\}$ as it is depicted on the right diagram of Figure \ref{d17b}. We define the rotated function $\tilde\Phi(z)$ given by the formula
\begin{equation*}
\tilde\Phi(z):=\check\Phi(iz),\quad z\in \tilde\Omega_{r}\cup\tilde\Omega_{l},
\end{equation*}
where the branch cut is chosen such that $\mathrm{arg}\,z\in(-3\pi/2,\pi/2)$. Then we have the following direct consequence of Proposition \ref{prop-l-h}.
\begin{proposition}\label{prop-tilde-l}
The function $\tilde\Phi(z)$ is a solution of the following RH problem. \\[2pt]
\noindent\makebox[5.5mm][l]{$(a)$}\parbox[t][][t]{121mm}{The function $\tilde\Phi(z)$ is an analytic function on $\C\setminus\tilde\Sigma$.}\\[2pt]
\noindent\makebox[5.5mm][l]{$(b)$}\parbox[t][][t]{121mm}{We have the following jump relation 
$$\begin{aligned}
&\tilde\Phi_+(z)=\tilde\Phi_-(z)S_{-}^{-1}, && z\in\tilde\Sigma, \ \mathrm{Im}\,z>0, \\  
&\tilde\Phi_+(z)=\tilde\Phi_-(z)S_{+}, && z\in\tilde\Sigma, \ \mathrm{Im}\,z<0.
\end{aligned}$$}\\
\noindent\makebox[5.5mm][l]{$(c)$}\parbox[t][][t]{121mm}{At $z=0$, the function $\tilde\Phi(z)$ has the same behavior as $\check\Phi(z)$ in \eqref{asym-orig-1} and \eqref{asym-orig-2}.}\\[5pt]
\noindent\makebox[5.5mm][l]{$(d)$}\parbox[t][][t]{121mm}{We have the following asymptotic relation
\begin{equation}\label{c-asy}
\tilde\Phi(z) =(I+O(z^{-1}))e^{iz\sigma_3}, \quad z\to\infty.
\end{equation}}
\end{proposition}
In view of \eqref{stokes2} and \eqref{stokes-11}, we have $s_{1} + s_{3} = -2\sin(\pi\alpha)$,
which implies that
\begin{align}\label{factor-1}
S_+ = \begin{pmatrix}1&0\\[5pt] -2sin(\pi\alpha)& 1\end{pmatrix}=
\begin{pmatrix}1& 0\\[5pt] s_1& 1\end{pmatrix}\begin{pmatrix}1& 0\\[5pt] s_3& 1\end{pmatrix}
\end{align}
and furthermore 
\begin{align}\label{factor-2}
S_-^{-1}=\begin{pmatrix}1&2sin(\pi\alpha)\\[5pt] 0 & 1\end{pmatrix} = \begin{pmatrix}1& -s_1\\[5pt] 0& 1\end{pmatrix}\begin{pmatrix}1& -s_3\\[5pt] 0& 1\end{pmatrix}.
\end{align}
The contour $\tilde\Sigma$ together with the four rays $\mathrm{arg}\,z = \pm\frac{\pi}{4}$ and $\mathrm{arg}\,z = \pm\frac{3\pi}{4}$ divide the complex plane on six regions as it is shown on the left diagram of Figure \ref{d5b}. Then we can represent the sets $\tilde\Omega_{l}$ and $\tilde\Omega_{r}$ in the form of the following sums $$\tilde\Omega_{r}=\bar\Omega_{r}^{1}\cup\bar\Omega_{r}^{2}\cup\bar\Omega_{r}^{3},\quad\tilde\Omega_{l}=\bar\Omega_{l}^{1}\cup\bar\Omega_{l}^{2}\cup\bar\Omega_{l}^{3}.$$
Assume that $\bar\Sigma$ is the contour consisting of four rays $\mathrm{arg}\,z = \pm\frac{\pi}{4}$ and $\mathrm{arg}\,z = \pm\frac{3\pi}{4}$ as it is shown on the right diagram of Figure \ref{d5b}. We define the function $\bar\Phi(z)$ by
\begin{gather*}
\bar\Phi(z) := \tilde\Phi(z)\begin{pmatrix}1& -s_3\\[3pt] 0& 1\end{pmatrix}, \ z\in \bar\Omega_r^1, \quad
\bar\Phi(z) := \tilde\Phi(z)\begin{pmatrix}1& 0\\[3pt] s_1& 1\end{pmatrix}^{-1}, \ z\in \bar\Omega_r^3  \\
\bar\Phi(z) := \tilde\Phi(z)\begin{pmatrix}1& -s_1\\[3pt] 0& 1\end{pmatrix}^{-1}, \ z\in \bar\Omega_l^1, \quad
\bar\Phi(z) := \tilde\Phi(z)\begin{pmatrix}1& 0\\[3pt] s_3& 1\end{pmatrix}, \ z\in \bar\Omega_l^3  \\
\bar\Phi(z)=\tilde\Phi(z), \ z\in \bar\Omega_r^2\cup \bar\Omega_l^2.
\end{gather*}
In the following theorem we formulate the auxiliary RH problem that will be applied in the construction of the local parametrix for the deformed RH problem. 
\begin{theorem}\label{th-l}
The function $\bar\Phi(z)$ satisfies the following RH problem.\\[2pt]
\noindent\makebox[5.5mm][l]{$(a)$}\parbox[t][][t]{121mm}{The function $\bar\Phi(z)$ is an analytic function on $\C\setminus\bar\Sigma$;}\\[2pt]
\noindent\makebox[5.5mm][l]{$(b)$}\parbox[t][][t]{121mm}{On the contour $\bar\Sigma$, the following jump relation is satisfied
\begin{equation*}
\bar\Phi_+(z) = \bar\Phi_-(z)\bar S(z),\quad z\in \bar\Sigma,
\end{equation*}
where the jump matrix $\bar S(z)$ is given on the right diagram of Figure \ref{d5b}.}\\[2pt]
\noindent\makebox[5.5mm][l]{$(c)$}\parbox[t][][t]{121mm}{If $0<\mathrm{Re}\,\alpha<1/2$, then the function $\bar\Phi(z)$ has the following asymptotic behavior 
\begin{align*}
\bar\Phi(z) = O\begin{pmatrix}|z|^{-\alpha} & |z|^{-\alpha} \\[5pt]
|z|^{-\alpha} & |z|^{-\alpha}\end{pmatrix},\quad z\to 0,
\end{align*}
and furthermore, for $-1/2<\mathrm{Re}\,\alpha\le 0$, we have}\\
\noindent\makebox[5.5mm][l]{}\parbox[t][][t]{121mm}{
\begin{align*}
\bar\Phi(z) = O\begin{pmatrix}|z|^{\alpha} & |z|^{\alpha} \\[5pt]
|z|^{\alpha} & |z|^{\alpha}\end{pmatrix},\quad z\to 0.
\end{align*}}\\
\noindent\makebox[6mm][l]{$(d)$}\parbox[t][][t]{121mm}{The function $\bar\Phi(z)$ has the following behavior at infinity
\begin{equation*}
\bar\Phi(z) = \left(I+O\left(z^{-1}\right)\right)e^{iz\sigma_3}, \quad z\to\infty.
\end{equation*}}
\end{theorem}
\proof Using the factorizations \eqref{factor-1} and \eqref{factor-2}, it is not difficult to check that the conditions $(a)$ and $(b)$ are satisfied. On the other hand, using the point $(c)$ of Proposition \ref{prop-tilde-l}, we infer that the point $(c)$ holds true. To check that $(d)$ is valid, let us observe that the asymptotic relation \eqref{c-asy}, implies that
\begin{equation*}
\bar\Phi(z) = \left(I+O\left(z^{-1}\right)\right)e^{iz\sigma_3}, \quad z\to\infty, \ z\in\bar\Omega_r^2\cup \bar\Omega_l^2.
\end{equation*}
Furthermore, for $z\in\bar\Omega_{r}^{1}\cup\bar\Omega_{l}^{1}$, we have
\begin{align*}
\bar\Phi(z) e^{-iz\sigma_3} = \tilde\Phi(z)e^{-iz\sigma_3} \begin{pmatrix}1& c e^{2iz}\\[3pt]0&1\end{pmatrix},
\end{align*}
where the parameter $c$ is either $s_{1}$ or $-s_{3}$. Consequently we can write
\begin{align}\label{c-asy-1}
\bar\Phi(z) e^{-iz\sigma_3} - I = (\tilde\Phi(z) e^{-iz\sigma_3} - I)\begin{pmatrix}1& c e^{2iz}\\[3pt]0&1\end{pmatrix}+
\begin{pmatrix}1& c e^{2iz}\\[3pt]0&1\end{pmatrix} -I.
\end{align}
Considering the polar coordinates $z=|z|e^{i\varphi}$, we infer that the argument $\varphi$ is an element of the interval $(\frac{\pi}{4}, \frac{3\pi}{4})$, whenever 
$z\in\bar\Omega_{r}^{1}\cup\bar\Omega_{l}^{1}$, which implies that
\begin{align*}
|e^{2iz}| = |e^{2i|z|\cos\varphi}e^{-2|z|\sin\varphi}|\le e^{-2|z|\sin\varphi}\le e^{-\sqrt{2}|z|},\quad z\in\bar\Omega_{r}^{1}\cup \bar\Omega_{l}^{1}.
\end{align*}
Combining this with \eqref{c-asy-1}, we deduce that 
$$\bar\Phi(z) e^{-iz\sigma_3} - I = O(z^{-1}),\quad z\to\infty, \ z\in\bar\Omega_{r}^{1}\cup\bar\Omega_{l}^{1}.$$
Arguing in the similar way we can write
\begin{align*}
\bar\Phi(z) e^{-iz\sigma_3} = \tilde\Phi(z)e^{-iz\sigma_3} \begin{pmatrix}1& 0 \\[5pt] d e^{-2iz} &1\end{pmatrix},\quad 
\lambda\in\bar\Omega_{r}^{3}\cup \bar\Omega_{l}^{3},
\end{align*}
where the parameter $d$ is equal to either $-s_{1}$ or $s_{3}$. Then we can write
\begin{align}\label{c-asy-2}
\bar\Phi(z) e^{-iz\sigma_3} - I = (\tilde\Phi(z) e^{-iz\sigma_3} - I)\begin{pmatrix}1& 0\\[5pt] s_1 e^{-2iz}&1\end{pmatrix}+
\begin{pmatrix}1& 0 \\[5pt] s_1 e^{-2iz}&1\end{pmatrix} -I.
\end{align}
Considering the parameter $z$ in polar coordinates once again, we have $\varphi\in (-\frac{3\pi}{4}, -\frac{\pi}{4})$ for $z\in\bar\Omega_{r}^{3}\cup\bar\Omega_{l}^{3}$, and hence
\begin{align}\label{c-asy-kk}
|e^{-2iz}| = |e^{-2i|z|\cos\varphi}e^{2|z|\sin\varphi}|\le e^{2|z|\sin\varphi}\le e^{-\sqrt{2}|z|},\quad z\in\bar\Omega_{r}^{3}\cup \bar\Omega_{l}^{3}
\end{align}
By \eqref{c-asy-kk} and \eqref{c-asy-2}, we obtain
$$\bar\Phi(z) e^{-iz\sigma_3} - I = O(z^{-1}),\quad z\to\infty, \ z\in\bar\Omega_{r}^{3}\cup\bar\Omega_{l}^{3}$$
and the proof of theorem is completed. \hfill $\square$\\

\section{Asymptotic analysis of the function $\Phi(\lambda,x)$ for $x>0$}

\subsection{Contour deformation} 
\begin{figure}[h]
\includegraphics[scale=0.6]{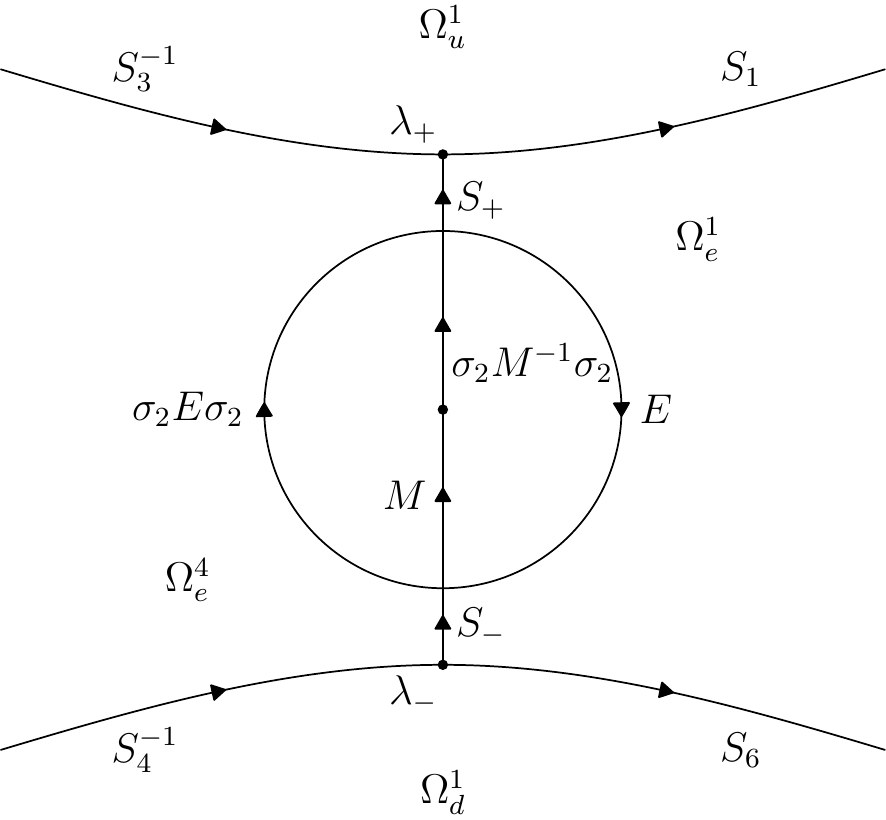}
\caption{The graph $\Sigma_1$ with the jump matrices.}
\label{d11bb}
\end{figure}

\begin{figure}[h]
\begin{subfigure}{0.49\textwidth}
\centering
\includegraphics[scale=0.55]{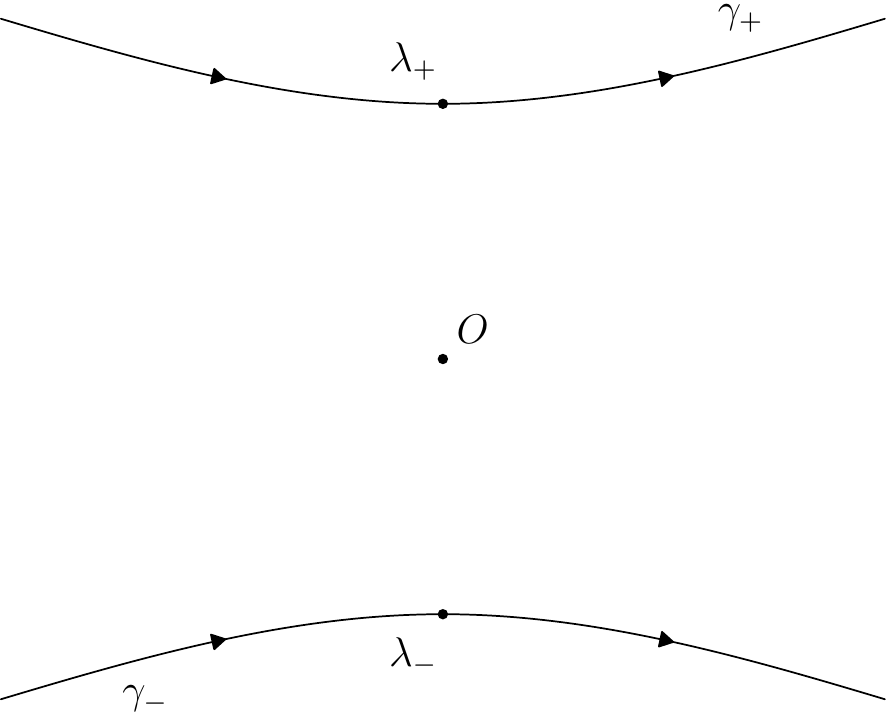}
\end{subfigure}
\begin{subfigure}{0.49\textwidth}
\centering
\includegraphics[scale=0.6]{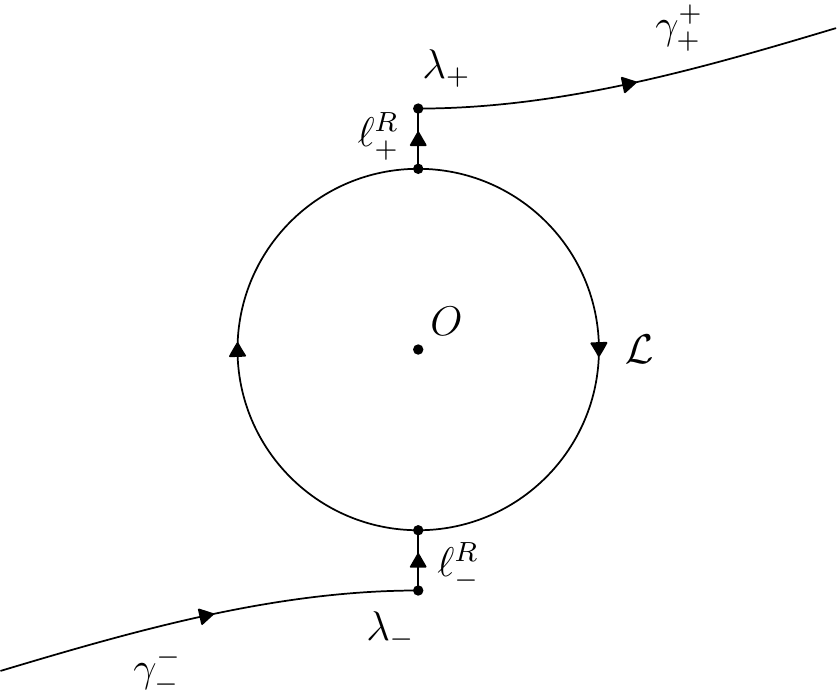}
\end{subfigure}
\caption{the contours $\Gamma_{1}$ and $\Gamma_{2}$ for the RH problems with the function $\chi^{1}(\lambda)$ and $\chi^{2}(\lambda)$, respectively.}
\label{reduced-cont2}
\end{figure}

Let us write $\lambda_{\pm}:=\pm\frac{i}{2}\sqrt{x}$ for the stationary points of the phase function $\theta(\lambda,x)$. From the equation $\mathrm{Im}\,\theta(z)=\mathrm{Im}\,\theta(\lambda_{\pm})=0$, it follows that the steepest descent paths passing through the points $\lambda_{\pm}$ are either the line $\mathrm{Re}\,\lambda = 0$ or the curves 
\begin{equation}\label{curve-gamma}
\gamma_\pm(t):=t\pm i\left(t^{2}/3 + x/4\right)^{1/2},\quad t\in\R,
\end{equation}
that are asymptotic to the rays $\mathrm{arg}\,\lambda = \pm\frac{\pi}{6},\pm\frac{5\pi}{6}$. Assume that $\ell_\pm$ is a vertical segment connecting the origin with the stationary point $\lambda_\pm$. Let us consider the contour $\Sigma_1$ consisting of the circle $C$ of radius $r>0$ (see \eqref{radius-r}), steepest descent paths $\gamma_{\pm}$ and vertical segments $\ell_{\pm}$, as it is depicted on the Figure \ref{d11bb}. 
We define the sets $\Omega^{1}_{e}:=\Omega^{3}_{u}\cup\Omega_{1}\cup\Omega^{3}_{d}$ and $\Omega^{4}_{e}:=\Omega^{2}_{u}\cup\Omega_{4}\cup\Omega^{2}_{d}$.
In view of \cite[Chapter 11]{MR2264522}, we infer that the function $\Phi(z)$ satisfying the RH problem (A1)\,--\,(A5) can be deformed to the function $\Phi^{1}(\lambda)$ with values in the space $M_{2\times 2}(\C)$, which satisfies the following deformed RH problem on the graph $\Sigma_1$. \\[3pt]
\noindent\makebox[8.5mm][l]{$(B1)$}\parbox[t][][t]{118mm}{Given $k\in \{u,d\}$ and $l\in \{1,4\}$, the restriction of the function $\Phi^{1}(\lambda)$ to the sets $\Omega_k^{1}$ and $\Omega_e^{l}$ are holomorphic and continuous up to $\overline{\Omega_k^{1}}$ and $\overline{\Omega_e^{l}}$, respectively.}\\[3pt]
\noindent\makebox[8.5mm][l]{$(B2)$}\parbox[t][][t]{118mm}{Given $k\in\{l,r\}$, the restrictions $\Phi^{1}_{|\Omega_{k}}$ is holomorphic and
$\Phi^{1}_{|\Omega_{k}}\in C(\overline{\Omega_{k}^\ve})$, where $\Omega^\ve_{k}:=\Omega_k\setminus\{\lambda\in\C \ | \ |\lambda|<\ve\}$ and $\ve>0$ is sufficiently small.}\\[5pt]
\noindent\makebox[8.5mm][l]{$(B3)$}\parbox[t][][t]{118mm}{The following jump relation is satisfied
\begin{equation*}
\Phi^{1}_{+}(\lambda) = \Phi^{1}_{-}(\lambda)S_{1}(\lambda),\quad \lambda\in\Sigma_1,
\end{equation*}}\\
\noindent\makebox[8.5mm][l]{}\parbox[t][][t]{118mm}{where the jump matrix is such that 
\begin{equation*}
S_{1}(\lambda)= S(\lambda), \ \ \lambda\in \rho_{+}\cup\rho_{-}\cup C,\quad 
S_{1}(\lambda)=S_\pm, \ \ \lambda\in\ell_\pm, \ |\lambda|>r
\end{equation*}}\\
\noindent\makebox[8.5mm][l]{}\parbox[t][][t]{118mm}{and furthermore, if we define 
\begin{equation*}
\begin{aligned}
\gamma_\pm^+:=\{\lambda\in\gamma_\pm \ | \ \mathrm{Re}\,\gamma_\pm>0\}, \quad \gamma_\pm^-:=\{\lambda\in\gamma_\pm \ | \ \mathrm{Re}\,\gamma_\pm<0\},
\end{aligned}
\end{equation*}}\\
\noindent\makebox[8.5mm][l]{}\parbox[t][][t]{118mm}{then the matrix $S_{1}$ has the following form:
\begin{gather*}
S_{1}(\lambda)= S_3^{-1}, \ \ \lambda\in\gamma_{+}^{-},\quad S_{1}(\lambda)= S_1, \ \ \lambda\in\gamma_{+}^{+},\\
S_{1}(\lambda)= S_4^{-1}, \ \ \lambda\in\gamma_{-}^{-},\quad S_{1}(\lambda)= S_6, \ \ \lambda\in\gamma_{-}^{+}.
\end{gather*}}\\
\noindent\makebox[8.5mm][l]{$(B4)$}\parbox[t][][t]{118mm}{The function $\Phi^{1}_r(\lambda)\lambda^{-\alpha\sigma_3}$ is bounded for $\lambda$ sufficiently close to zero, where the branch of the multifunction $\lambda^{-\alpha}$ is chosen arbitrarily.}\\[3pt]
\noindent\makebox[8.5mm][l]{$(B5)$}\parbox[t][][t]{118mm}{As $\lambda\to\infty$, the function $\Phi^{1}(\lambda)$ has the following asymptotic behavior
$$\Phi^{1}(\lambda) = (I + O(\lambda^{-1}))e^{-\theta(\lambda)\sigma_{3}}.$$}\\[5pt]
By the construction of the function $\Phi^{1}(\lambda)$, we have
\begin{align}\label{funct-phi}
\Phi^{1}(\lambda) =\Phi(\lambda), \quad  \lambda\in \Omega_{1}\cup\Omega_{4}\cup\Omega_{r}\cup\Omega_{l}\cup\Omega^{1}_{u}\cup\Omega^{1}_{d}.
\end{align}
Let us consider two oriented contours $\Gamma_{1}$ and $\Gamma_{2}$ (see Figure \ref{reduced-cont2}), where $\Gamma_{1}:=\gamma_{+}\cup\gamma_{-}$, which is the union of the steepest descent paths as is shown on the left diagram of Figure \ref{reduced-cont2} and furthermore $\Gamma_{2}$ is a union of the circle $\LL$ of the radius $R>0$ (see \eqref{radius-r}), the segment $\ell_{\pm}^{R}$ is the part of the curve $\ell_{\pm}$ lying outside the circle $\LL$ and the steepest descent paths $\gamma_{\pm}^{\pm}$.
In \cite[p. 418--421]{MR2264522} it is proved that, for sufficiently large $x>0$, the function $\Phi^{1}(\lambda)$ has the following representation 
\begin{equation}\label{representation}
\Phi^{1}(\lambda,x):=\left\{\begin{aligned}
&\chi^{1}(\lambda,x)\chi^{2}(\lambda,x)\hat\Phi(z(\lambda),x),&&\quad |\lambda|<R,\\
&\chi^{1}(\lambda,x)\chi^{2}(\lambda,x)e^{-\theta(\lambda)\sigma_3},&&\quad |\lambda|>R,
\end{aligned}\right.
\end{equation}
where the component functions satisfy the following conditions.\\[2pt]
\noindent\makebox[8.5mm][l]{$(C1)$}\parbox[t][][t]{119mm}{The function $z(\lambda)$ is defined by equation \eqref{map-z} and $\hat\Phi(\lambda)$ is a solutions of the RH problem from Proposition \ref{prop-bessel}.}\\[5pt]
\noindent\makebox[8.5mm][l]{$(C2)$}\parbox[t][][t]{119mm}{The function $\chi^{1}(\lambda)$ is holomorphic on the sets $\C\setminus\Gamma_{1}$ and has the following asymptotic behavior $\chi^{1}(\lambda)\to I$ as $\lambda\to\infty$.}\\[2pt]
\noindent\makebox[8.5mm][l]{$(C3)$}\parbox[t][][t]{119mm}{We have the following jump relation
$$\chi^{1}_+(\lambda) = \chi^{1}_-(\lambda)H_{1}(\lambda), \ \ \lambda\in\Gamma_{1},$$ where the matrix function $H_{1}$ (defined precisely in 
\cite[p. 419]{MR2264522}) satisfies for some $c>0$, the inequality
\begin{align}\label{est-t}
\|H_{1}(\lambda) - I\|&\le c e^{-\frac{2}{3}|x|^{3/2}}e^{-c|x|^{1/2}|\lambda - \lambda_\pm|^2}, \quad \lambda\in\gamma_\pm, \ x>0.
\end{align}}\\
\noindent\makebox[8.5mm][l]{$(C4)$}\parbox[t][][t]{119mm}{The function $\chi^{2}(\lambda)$ is holomorphic on the set $\C\setminus\Gamma_{2}$ and has the following asymptotic behavior $\chi^{2}(\lambda)\to I$ as $\lambda\to\infty$.}\\[2pt]
\noindent\makebox[8.5mm][l]{$(C5)$}\parbox[t][][t]{119mm}{We have the following jump relations 
$$\chi^{2}_+(\lambda) = \chi^{2}_-(\lambda)H_{2}(\lambda), \ \ \lambda\in\Gamma_{2},$$ where the matrix function $H_{2}$ (defined precisely in 
\cite[p. 421]{MR2264522}) satisfies for some $c>0$, the inequality
\begin{align}\label{est-t-2}
\|H_{2}(\lambda) - I\|&\le 
\left\{\begin{aligned}
&ce^{-c|x||\lambda|}&&\text{for}\ \ \lambda\in \ell_{\pm}^{R}\cup\gamma_{\pm}^{\pm},\\
&cR^{-3}&&\text{for}\ \ \lambda\in\mathcal{L}.
\end{aligned}\right.
\end{align}}
In the following two lemmata we provide simple estimates of the jump matrices $H_{1}$ and $H_{2}$ in the spaces $L^{p}(\Gamma_{1})$ and $L^{p}(\Gamma_{2})$, respectively.
\begin{lemma}\label{lem-est-2}
Given $1\le p<+\infty$, there is constant $c_{p}>0$ such that 
\begin{align}\label{lp-est1b}
\|H_{1}-I\|_{L^{p}(\Gamma_{1})}\le c_{p}|x|^{-\frac{1}{4p}}e^{-\frac{2}{3}|x|^{3/2}},\quad x>0.
\end{align}
Furthermore there is a constant $c_{\infty}>0$ such that 
\begin{align}\label{lp-est2b}
\|H_{1}-I\|_{L^{\infty}(\Gamma_{1})}\le c_{\infty}e^{-\frac{2}{3}|x|^{3/2}},\quad x>0.
\end{align}
\end{lemma}
\proof As a consequence of the inequality \eqref{est-t}, we obtain 
\begin{align*}
\|H_{1}(\lambda) - I\|\le c e^{-\frac{2}{3}|x|^{3/2}}, \quad \lambda\in\gamma_\pm, \ x>0,
\end{align*}
which immediately gives \eqref{lp-est2b}. For the proof of \eqref{lp-est1b}, observe that from the definition \eqref{curve-gamma} of the curves $\gamma_{\pm}$, we have
\begin{align*}
|\gamma_\pm(t) - \lambda_\pm| = |t\pm i\left(t^{2}/3 + x/4\right)^{1/2}\mp i\sqrt{x}/2 |\ge t ,\quad t\in\R.
\end{align*}
Combining this with \eqref{est-t} yields
\begin{align*}
&\|H_{1}-I\|_{L^p(\gamma_{\pm})}^p  = \int_{\gamma_\pm} \|H_{1}(\lambda) - I\|^p\,|d\lambda|
\lesssim e^{-\frac{2p}{3}|x|^{3/2}}\int_{\gamma_\pm}e^{-pc|x|^{1/2}|\lambda - \lambda_\pm|^2}\,|d\lambda|\\
& = e^{-\frac{2p}{3}|x|^{\frac{3}{2}}}\!\!\int_\R e^{-pc|x|^{\frac{1}{2}}|\gamma_\pm(t) - \lambda_\pm|^2}|\dot\gamma_\pm(t)|\,dt
 \lesssim e^{-\frac{2p}{3}|x|^{\frac{3}{2}}}\!\!\int_\R e^{-pc|x|^{\frac{1}{2}}t^2}dt  \sim \frac{1}{|x|^{\frac{1}{4}}}e^{-\frac{2p}{3}|x|^{\frac{3}{2}}}
\end{align*}
and completes the proof of the inequality \eqref{lp-est1b}. \hfill $\square$

\begin{lemma}\label{lem-est-1}
Given $1\le p<+\infty$, we have the following asymptotic behavior 
\begin{align}\label{lp-est1}
\|H_{2}-I\|_{L^{p}(\Gamma_{2})}=O(|x|^{-\frac{3p-1}{2p}}),\quad x\to+\infty.
\end{align}
Furthermore we have
\begin{align}\label{lp-est2}
\|H_{2}-I\|_{L^{\infty}(\Gamma_{2})}=O(|x|^{-\frac{3}{2}}),\quad x\to+\infty.
\end{align}
\end{lemma}
\proof In view of the inequality \eqref{est-t-2} and the definition \eqref{radius-r}, we find that
\begin{equation}\label{ineq-rr-11}
\|H_{2}(\lambda) - I\|\le 
\left\{\begin{aligned}
&ce^{-c|x|^{\frac{3}{2}}/4}&&\text{for}\ \ \lambda\in \ell_{\pm}^{R}\cup\gamma_{\pm}^{\pm},\\
&64 c|x|^{-\frac{3}{2}}&&\text{for}\ \ \lambda\in\mathcal{L},
\end{aligned}\right.
\end{equation}
and hence \eqref{lp-est2} holds. To show \eqref{lp-est1}, we observe that \eqref{ineq-rr-11} gives
\begin{align}\label{ee1}
\int_{\LL}\|H_{2}(\lambda) - I\|^{p}\,|d\lambda| \lesssim |x|^{-(3p-1 )/2}.
\end{align}
On the other hand simple calculations show that
\begin{align}\label{est-t2-bb}
|\gamma_{\pm}(t)| \ge\sqrt{3}|t|/3+x^{\frac{1}{2}}/4\quad\text{and}\quad |\dot\gamma_\pm(t)|\le 2\quad\text{for}\quad t\in\R. 
\end{align}
For the estimate on the steepest descent path $\gamma_{+}^{+}$, we use \eqref{est-t-2} and \eqref{est-t2-bb} to obtain
\begin{equation}
\begin{aligned}\label{ee2}
&\int_{\gamma_{+}^{+}}\|H_{2}(\lambda) - I\|^{p}\,|d\lambda| \lesssim \int_{\gamma_{+}^{+}}e^{-pc|x||\lambda|}\,|d\lambda| 
= \int_{0}^{\infty}e^{-pc|x||\gamma_{+}(t)|}|\dot\gamma_{+}(t)|\,dt\\
&\ \lesssim \int_{0}^{\infty}e^{-\frac{1}{4}pc|x|(\sqrt{3}t/3+x^{\frac{1}{2}})}\,dt 
 = e^{-\frac{1}{4}pc|x|^{\frac{3}{2}}}\int_{0}^{\infty}e^{-\frac{\sqrt{3}}{3}pc|x|t}\,dt
\sim \frac{1}{|x|}e^{-\frac{1}{4}pc|x|^{\frac{3}{2}}}
\end{aligned}
\end{equation}
and furthermore, for the curve $\ell_{+}^{R}$, we have
\begin{equation}
\begin{aligned}\label{est-t2}
&\int_{\ell_{+}^{R}}\|H_{2}(\lambda) - I\|^{p}\,|d\lambda| \lesssim \int_{\ell_{+}^{R}}e^{-pc|x||\lambda|}\,|d\lambda| 
= \int_{\frac{1}{4}|x|^{1/2}}^{\frac{1}{2}|x|^{1/2}}e^{-pc|x||t|}\,dt\\
& \qquad \sim |x|^{-1}\left(e^{-pc|x|^{3/2}/2}-e^{-pc|x|^{3/2}/4}\right)
\le |x|^{-1}e^{-pc|x|^{3/2}/2}.
\end{aligned}
\end{equation}
Proceeding in the similar way we obtain analogous estimates for the steepest descent paths $\gamma_{-}^{-}$ and $\ell_{-}^{R}$, that read as follows
\begin{align}\label{est-t3}
\hspace{-2pt}\int_{\gamma_{-}^{-}\cup \ell_{-}^{R}}\|H_{2}(\lambda)-I\|^{p}\,|d\lambda|\lesssim |x|^{-1} e^{-\frac{pc}{4}|x|^{\frac{3}{2}}}.
\end{align}
Combining \eqref{ee1}, \eqref{ee2}, \eqref{est-t2} and \eqref{est-t3} we obtain the inequality \eqref{lp-est1} 
and the proof of lemma is completed. \hfill $\square$


\subsection{Representation of solutions of the deformed RH problem} We intend to prove Proposition \ref{prop-p-g}, which gives representation of the function $P(x)$ and provides information about its asymptotic behavior as $x\to+\infty$. We start with the following result concerning asymptotics of the function $\chi^{1}(0,x)$ as $x\to+\infty$.
\begin{proposition}\label{prop-as-1}
There is $x_{1}>0$ such that, for any $x>x_{1}$, the RH problem $(C2)$, $(C3)$ has a unique solution $\chi^{1}(\lambda,x)$ with the property that 
\begin{align}\label{est-kk}
\|\chi^{1}(0,x) - I\|= O(|x|^{-3/4}),\quad x\to+\infty.
\end{align}
\end{proposition}
\proof Applying Lemma \ref{lem-est-2}, we obtain the existence of $c_{0},x_{0}>0$ such that 
\begin{align}\label{est-t9-bb}
\|H_{1} - I\|_{L^{1}(\Gamma_{1})}\le c_{0}|x|^{-\frac{1}{4}},\ \ \|H_{1} - I\|_{(L^{2}\cap L^{\infty})(\Gamma_{1})}\le c_{0}|x|^{-\frac{1}{8}}, \quad x>x_{0}.
\end{align}
Let us consider a complex linear map $\KK_{1}:L^{2}_{I}(\Gamma_{1})\to L^{2}_{I}(\Gamma_{1})$ given by 
\begin{equation*}
\KK_{1}(\rho):= \CC_{-}(\rho(H_{1} - I)),\quad \rho\in L^{2}_{I}(\Gamma_{1}),
\end{equation*}
where $\mathcal{C}_-$ is the Cauchy operator on the contour $\Gamma_{1}$. If $\rho\in L^{2}_{I}(\Gamma_{1})$ is such that $\rho = \rho_{0}+\rho_{\infty}$, where $\rho_{0}\in L^{2}(\Gamma_{1})$ and $\rho_{\infty}\in M_{2\times 2}(\C)$, then  
\begin{equation*}
\KK_{1}(\rho)= \CC_{-}(\rho_{0}(H_{1} - I))+ \CC_{-}(\rho_{\infty}(H_{1} - I)),
\end{equation*}
which implies that $\KK_{1}(\rho)\in L^{2}(\Gamma_{1})$ and the following estimates hold
\begin{equation}
\begin{aligned}\label{est-t5b}
\|\KK_{1}(\rho)\|_{L^{2}(\Gamma_{1})}&\le \|\CC_{-}\|_{L^{2}(\Gamma_{1})}\left(\|\rho_{0}(H_{1} - I)\|_{L^{2}(\Gamma_{1})}+
\|\rho_{\infty}\|\|(H_{1} - I)\|_{L^{2}(\Gamma_{1})}\right)\\
&\le \|\CC_{-}\|_{L^{2}(\Gamma_{1})}\|H_{1}-I\|_{(L^{2}\cap L^{\infty})(\Gamma_{1})}\left(\|\rho_{0}\|_{L^{2}(\Gamma_{1})}+\|\rho_{\infty}\|\right)\\
&\le \|\CC_{-}\|_{L^{2}(\Gamma_{1})}\|\rho\|_{L^{2}_{I}(\Gamma_{1})}\|H_{1}-I\|_{(L^{2}\cap L^{\infty})(\Gamma_{1})},\quad x>0.
\end{aligned}
\end{equation}
By \eqref{est-t2-bb} we infer that $\gamma_{\pm}$ are Lipschitz curves and their Lipschitz constants are not dependent from the parameter $x>0$. Then, from \cite[Section 2.5.4]{MR3450072} it follows that the norm of the Cauchy operator $\CC_{-}$ satisfies the following inequality 
\begin{equation}\label{ee9b}
\|\CC_{-}\|_{L^{2}(\Gamma_{1})}\le m,\quad x>0,
\end{equation}
where $m>0$ is a constant independent from the parameter $x>0$. Using \eqref{est-t5b} together with \eqref{est-t9-bb} and \eqref{ee9b}, we find that
\begin{align}\label{ee10b}
\|\KK_{1}(\rho)\|_{L^{2}_{I}(\Gamma_{1})}\lesssim |x|^{-\frac{1}{8}}\|\rho\|_{L^{2}_{I}(\Gamma_{1})}, \quad x>x_{0},
\end{align}
and consequently
\begin{align}\label{sum-s}
\|\KK_{1} I\|_{L^{2}_{I}(\Gamma_{1})} \lesssim |x|^{-\frac{1}{8}}, \quad x>x_{0}.
\end{align}
Furthermore \eqref{ee10b} implies that, there is $x_{1}>x_{0}$ such that 
\begin{align}\label{sum-series}
\|\KK_{1}\|_{L^{2}_{I}(\Gamma_{1})}\le 1/2, \quad x>x_{1}.
\end{align}
Hence the element $\rho\in L^{2}_{I}(\Gamma_{1})$, given by the series $\rho:=\sum_{i=0}^{\infty} \KK_{1}^{i}I$ 
is the unique solution of the equation $\rho - \KK_{1}(\rho) = I$ and the inequalities \eqref{sum-s} and \eqref{sum-series} yield
\begin{equation}\label{ee15b}
\|\rho - I\|_{L^{2}_{I}(\Gamma_{1})}\!\le\! \|\KK_{1} I\|_{L^{2}(\Gamma_{1})}\sum_{i=0}^{\infty} \|\KK_{1}\|^{i}_{L^{2}_{I}(\Gamma_{1})}\!\lesssim\! \|\KK_{1} I\|_{L^{2}(\Gamma_{1})} \!\lesssim\! |x|^{-\frac{1}{8}},  \quad x>x_{1}.
\end{equation}
Observe that the solution $\chi^{1}(\lambda)$ of the RH problem $(C2)$, $(C3)$ is given by 
\begin{align*}
\chi^{1}(\lambda) = I+\frac{1}{2\pi i}\int_{\Gamma}\frac{\rho(\xi)(H_{1}(\xi)-I)}{\xi-\lambda}\,d\xi,\quad \lambda\not\in\Gamma_{1}.
\end{align*}
Combining this with \eqref{sum-s} and the H\"older inequality implies that 
\begin{align*}
&\|\chi^{1}(0) - I\| \le \frac{1}{2\pi}\int_{\Gamma_{1}} \frac{\|\rho(\xi)(H_{1}(\xi)-I)\|}{|\xi|}\,|d\xi|\lesssim
|x|^{-\frac{1}{2}}\int_{\Gamma_{1}}\|\rho(\xi)(H_{1}(\xi)-I)\|\,|d\xi|\\
&\quad\le |x|^{-\frac{1}{2}}\int_{\Gamma_{1}}\|(\rho(\xi)-I)(H_{1}(\xi)-I)\|\,|d\xi|+
|x|^{-\frac{1}{2}}\int_{\Gamma_{1}}\|H_{1}(\xi)-I\|\,|d\xi| \\
&\quad\le |x|^{-\frac{1}{2}}\|\rho-I\|_{L^{2}(\Gamma_{1})}\|H_{1}-I\|_{L^{2}(\Gamma_{1})}
+|x|^{-\frac{1}{2}}\|H_{1}-I\|_{L^{1}(\Gamma_{1})},\quad x> 0,
\end{align*}
which together with \eqref{est-t9-bb} and \eqref{ee15b} yield
\begin{align*}
\|\chi^{1}(0) - I\|\lesssim |x|^{-\frac{1}{2}}|x|^{-\frac{1}{8}}|x|^{-\frac{1}{8}} + |x|^{-\frac{1}{2}}|x|^{-\frac{1}{4}} \sim |x|^{-\frac{3}{4}},\quad x>x_{1}.
\end{align*} 
This gives the inequality \eqref{est-kk} and completes the proof of proposition. \hfill $\square$\\

We proceed to the proof of the following proposition concerning asymptotic behavior of the function $\chi^{2}(0,x)$ as $x\to+\infty$.
\begin{proposition}\label{prop-as-2}
There is $x_{2}>0$ such that, for any $x>x_{2}$, the RH problem $(C4)$, $(C5)$ has a unique solution $\chi^{2}(\lambda,x)$ with the property that 
\begin{align}\label{asy-11bb}
\|\chi^{2}(0,x) - I\|= O(|x|^{-3/2}),\quad x\to+\infty.
\end{align}
\end{proposition}
\proof Observe that by Lemma \ref{lem-est-1}, there are $c_{0}>0$ and $x_{0}>0$ such that 
\begin{align}\label{est-t9}
\|H_{2} - I\|_{L^{1}(\Gamma_{2})}\le c_{0}|x|^{-1},\quad \|H_{2} - I\|_{(L^{2}\cap L^{\infty})(\Gamma_{2})}\le c_{0}|x|^{-\frac{5}{4}}, \quad x>x_{0}.
\end{align}
Let us assume that $\KK_{2}:L^{2}_{I}(\Gamma_{2})\to L^{2}_{I}(\Gamma_{2})$ is a complex linear map given by
\begin{equation*}
\KK_{2}(\rho):= \CC_{-}(\rho(H_{2} - I)),\quad \rho\in L^{2}_{I}(\Gamma_{2}),
\end{equation*}
where $\mathcal{C}_-$ is the Cauchy operator on the contour $\Gamma_{2}$. If we take $\rho\in L^{2}_{I}(\Gamma_{2})$ with $\rho = \rho_{0}+\rho_{\infty}$, where $\rho_{0}\in L^{2}(\Gamma_{2})$ and $\rho_{\infty}\in M_{2\times 2}(\C)$, then 
\begin{equation*}
\KK_{2}(\rho)= \CC_{-}(\rho_{0}(H_{2} - I))+ \CC_{-}(\rho_{\infty}(H_{2} - I)).
\end{equation*}
Therefore $\KK_{2}(\rho)\in L^{2}(\Gamma_{2})$ and we have the following estimates
\begin{equation}
\begin{aligned}\label{est-t5}
\|\KK_{2}(\rho)\|_{L^{2}(\Gamma_{2})}&\le \|\CC_{-}\|_{L^{2}(\Gamma_{2})}\left(\|\rho_{0}(H_{2} - I)\|_{L^{2}(\Gamma_{2})}+
\|\rho_{\infty}\|\|(H_{2} - I)\|_{L^{2}(\Gamma_{2})}\right)\\
&\le \|\CC_{-}\|_{L^{2}(\Gamma_{2})}\|H_{2}-I\|_{(L^{2}\cap L^{\infty})(\Gamma_{2})}\left(\|\rho_{0}\|_{L^{2}(\Gamma_{2})}+\|\rho_{\infty}\|\right)\\
&= \|\CC_{-}\|_{L^{2}(\Gamma_{2})}\|\rho\|_{L^{2}_{I}(\Gamma_{2})}\|H_{2}-I\|_{(L^{2}\cap L^{\infty})(\Gamma_{2})}.
\end{aligned}
\end{equation}
Although the contour $\Gamma_{2}$ depends on the parameter $x>0$, from \cite[Section 2.5.4]{MR3450072} we know that the norm of the operator $\CC_{-}$ satisfies the inequality 
\begin{equation}\label{ee9}
\|\CC_{-}\|_{L^{2}(\Gamma_{2})}\le m,\quad x>0,
\end{equation}
where $m>0$ is a constant. By the inequalities \eqref{est-t9} and \eqref{est-t5}, we have
\begin{align}\label{ee10}
\|\KK_{2}(\rho)\|_{L^{2}_{I}(\Gamma_{2})}\lesssim |x|^{-5/4}\|\rho\|_{L^{2}_{I}(\Gamma_{2})}, \quad x>x_{0},
\end{align}
which, in particular, implies that 
\begin{align}\label{ee11}
\|\KK_{2}I\|_{L^{2}_{I}(\Gamma_{2})} \lesssim |x|^{-5/4}, \quad x>x_{0}.
\end{align}
Furthermore \eqref{ee10} shows that there is a constant $x_{1}>x_{0}$ such that 
\begin{equation*}
\|\KK_{2}\|_{L^{2}_{I}(\Gamma_{2})}<1/2,\quad x>x_{1},
\end{equation*}
which implies that the equation $\rho - \KK_{2}(\rho) = I$ has a unique solution, given by the series 
$\rho = \sum_{i=0}^{\infty} \KK_{2}^{i}I$, which is convergent in the space $L^{2}_{I}(\Gamma_{2})$. Therefore the inequalities \eqref{ee9}, \eqref{ee10} and \eqref{ee11} provide
\begin{equation}
\begin{aligned}\label{ee15}
\|\rho - I\|_{L^{2}_{I}(\Gamma_{2})}\!\le\!\|\KK_{2}I\|_{L^{2}(\Gamma_{2})}\!\sum_{i=0}^{\infty} \|\KK_{2}\|^{i}_{L^{2}_{I}(\Gamma_{2})}
\!\lesssim\! \|\KK_{2}I\|_{L^{2}(\Gamma_{2})} \!\lesssim\! |x|^{-\frac{5}{4}}, \ \ x>x_{1}.
\end{aligned}
\end{equation}
Using the representation formula for the solutions of the RH problem, we obtain
\begin{align*}
\chi^{2}(\lambda)= I +\frac{1}{2\pi i}\int_{\Gamma_{2}} \frac{\rho(\xi)(H_{2}(\xi)-I)}{\xi-\lambda}\,d\xi, \quad \lambda\not\in\Gamma_{2},
\end{align*}
which together with the H\"older inequality imply that 
\begin{align*}
&\|\chi^{2}(0) - I\| \le \frac{1}{2\pi}\int_{\Gamma_{2}} \frac{\|\rho(\xi)(H_{2}(\xi)-I)\|}{|\xi|}\,|d\xi|\lesssim
|x|^{-\frac{1}{2}}\int_{\Gamma_{2}}\|\rho(\xi)(H_{2}(\xi)-I)\|\,|d\xi|\\
&\quad\le |x|^{-\frac{1}{2}}\int_{\Gamma_{2}}\|(\rho(\xi)-I)(H_{2}(\xi)-I)\|\,|d\xi|+
|x|^{-\frac{1}{2}}\int_{\Gamma_{2}}\|H_{2}(\xi)-I\|\,|d\xi| \\
&\quad\le |x|^{-\frac{1}{2}}\|\rho-I\|_{L^{2}(\Gamma_{2})}\|H_{2}-I\|_{L^{2}(\Gamma_{2})}
+|x|^{-\frac{1}{2}}\|H_{2}-I\|_{L^{1}(\Gamma_{2})},\quad x> x_{0}.
\end{align*}
Combining this with \eqref{est-t9} and \eqref{ee15}, gives
\begin{align*}
&\|\chi^{2}(0) - I\|\lesssim |x|^{-\frac{1}{2}}|x|^{-\frac{5}{4}}|x|^{-\frac{5}{4}} + |x|^{-\frac{1}{2}}|x|^{-1} = |x|^{-3}+|x|^{-\frac{3}{2}}\lesssim |x|^{-\frac{3}{2}}, \ \ x>x_{1},
\end{align*}
and hence the proof of the relation \eqref{asy-11bb} is completed. \hfill $\square$\\

In the following proposition we obtain the representation formula for the matrix $P(x)$ provided $x>0$ is sufficiently large.
\begin{proposition}\label{prop-p-g}
There is $x_{+}>0$ such that, for any $x>x_{+}$, we have
\begin{align*}
P(x)= \frac{1}{2}\chi^{1}(0,x)\chi^{2}(0,x) e^{-i\frac{\pi}{4}\sigma_3}\begin{pmatrix}1&1-2\alpha \\[5pt]-1&1-2\alpha\end{pmatrix}e^{2\pi i\alpha\sigma_{3}}
(-ix)^{\alpha\sigma_3} D.
\end{align*}
\end{proposition}
\proof According to the notation of Section \ref{sec-abl-seg} and the equality \eqref{eq-zet-phi}, we have
\begin{align*}
P(x) = \lim_{\lambda\to 0} Z(\lambda,x) = \lim_{\lambda\to 0}\Phi(\lambda,x)e^{\theta(\lambda,x)\sigma_{3}}\lambda^{-\alpha\sigma_3},
\end{align*}
where in the above limit the parameter $\lambda$ belongs to the set $\Omega_{r}$. By Propositions \ref{prop-as-1} and \ref{prop-as-2}, there is $x_{+}>0$ such that, for any $x>x_{+}$, the representation formula \eqref{representation} holds and the component functions $\chi^{1}(\lambda)$ and $\chi^{2}(\lambda)$ have the asymptotic behaviors \eqref{est-kk} and \eqref{asy-11bb}, respectively. Therefore, if we confine our attention to the ray $\mathrm{arg}\,\lambda = 0$, then $z(\lambda)\in\hat\Omega_d$ for sufficiently small $|\lambda|>0$ (see the left diagram of Figure \ref{d17b}), which together with \eqref{funct-phi} gives
\begin{equation}
\begin{aligned}\label{passing-lim}
&\Phi(\lambda,x)e^{\theta(\lambda,x)\sigma_{3}}\lambda^{-\alpha\sigma_3} = \Phi^{1}(\lambda,x)e^{\theta(\lambda,x)\sigma_{3}}\lambda^{-\alpha\sigma_3}\\
&\qquad=\chi^{1}(\lambda,x)\chi^{2}(\lambda,x)\hat\Phi(z(\lambda))e^{\theta(\lambda,x)\sigma_{3}}\lambda^{-\alpha\sigma_3} \\
&\qquad=\chi^{1}(\lambda,x)\chi^{2}(\lambda,x)\hat\Phi^{0}(e^{2\pi i}z(\lambda,x))De^{\theta(\lambda,x)\sigma_{3}}\lambda^{-\alpha\sigma_3}\\
&\qquad=\chi^{1}(\lambda,x)\chi^{2}(\lambda,x)\hat\Phi^{0}(e^{2\pi i}z(\lambda,x)) \lambda^{-\alpha\sigma_3}De^{\theta(\lambda,x)\sigma_{3}},
\end{aligned}
\end{equation}
where the last equality follows from the fact that $D$ is a diagonal matrix. Then
\begin{equation*}
\begin{aligned}
&\lim_{\lambda\in\R_{+},\lambda\to 0}(e^{2\pi i}z(\lambda,x))^{\alpha\sigma_3} \lambda^{-\alpha\sigma_3} = \lim_{\lambda\in\R_{+},\lambda\to 0} e^{2\pi i\alpha\sigma_{3}}[z(\lambda,x)/\lambda]^{\alpha\sigma_3} \\
&\qquad = \lim_{\lambda\in\R_{+},\lambda\to 0}e^{2\pi i\alpha\sigma_{3}} [-i(4\lambda^{2}/3 + x)]^{\alpha\sigma_3} = e^{2\pi i\alpha\sigma_{3}}(-ix)^{\alpha\sigma_3}.
\end{aligned}
\end{equation*}
which combined with Lemma \ref{lem-lim-1} and the equality 
\begin{align*}
\hat\Phi^{0}(e^{2\pi i}z(\lambda,x)) \lambda^{-\alpha\sigma_3} 
= \hat\Phi^{0}(e^{2\pi i}z(\lambda,x))(e^{2\pi i}z(\lambda,x))^{-\alpha\sigma_3} (e^{2\pi i}z(\lambda,x))^{\alpha\sigma_3} \lambda^{-\alpha\sigma_3}
\end{align*}
yields the following limit 
\begin{align}\label{lim-kk-11}
\lim_{\lambda\in\R_{+},\lambda\to 0}\hat\Phi^{0}(e^{2\pi i}z(\lambda,x)) \lambda^{-\alpha\sigma_3}  = \frac{1}{2}e^{-i\frac{\pi}{4}\sigma_3}\begin{pmatrix}1\hspace{-3pt}&1-2\alpha \\[5pt]-1\hspace{-3pt}&1-2\alpha\end{pmatrix}e^{2\pi i\alpha\sigma_{3}}(-ix)^{\alpha\sigma_3}.
\end{align}
Using \eqref{lim-kk-11} and the fact that the functions $\chi^{1}(\lambda)$ and $\chi^{2}(\lambda)$ are holomorphic in a neighborhood of the origin, we pass in \eqref{passing-lim} to the limit with $\lambda\to 0$ along the ray $\mathrm{arg}\,\lambda = 0$ and obtain
\begin{align*}
P(x)&=\lim_{\lambda\in\R_{+},\lambda\to 0}\chi^{1}(\lambda,x)\chi^{2}(\lambda,x)\hat\Phi^{0}(e^{2\pi i}z(\lambda,x)) \lambda^{-\alpha\sigma_3}De^{\theta(\lambda,x)\sigma_{3}}\\
&=\frac{1}{2}\chi^{1}(0,x)\chi^{2}(0,x) e^{-i\frac{\pi}{4}\sigma_3}\begin{pmatrix}1&1-2\alpha \\[5pt]-1&1-2\alpha\end{pmatrix} e^{2\pi i\alpha\sigma_{3}}
(-ix)^{\alpha\sigma_3} D. 
\end{align*}
Thus the proof of proposition is completed. \hfill $\square$

\section{Asymptotic analysis of the function $\Phi(\lambda,x)$ for $x<0$}

\subsection{Contour deformation}\label{sect-constr-sigma-t}
\begin{figure}[h]
\begin{subfigure}{0.49\textwidth}
\centering
\includegraphics[scale=0.7]{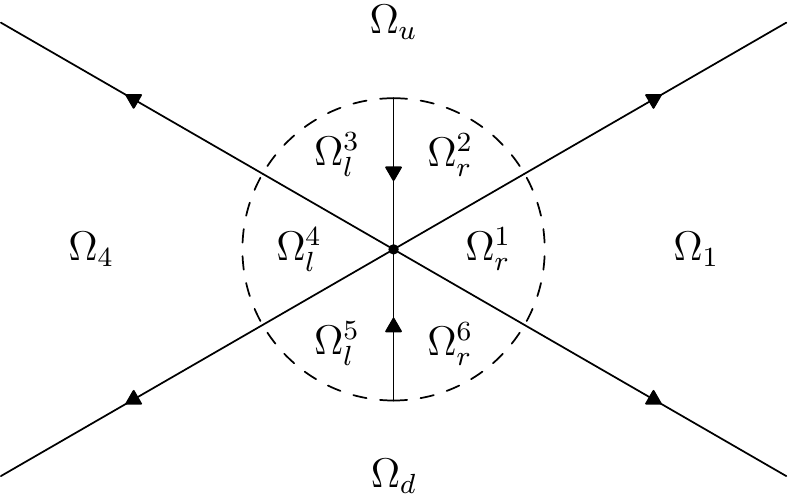}
\end{subfigure}
\begin{subfigure}{0.49\textwidth}
\centering
\includegraphics[scale=0.75]{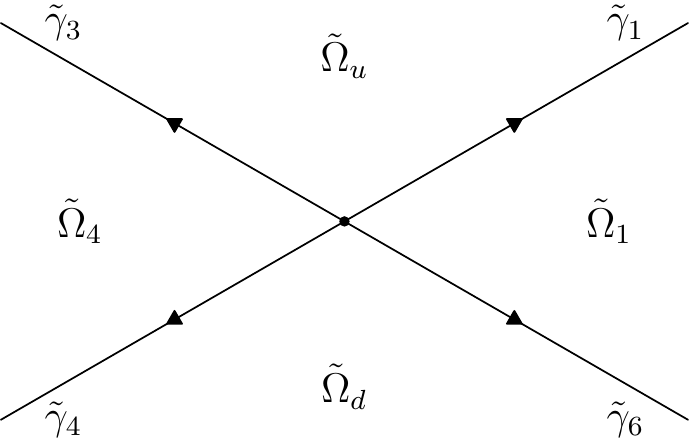}
\end{subfigure}
\caption{Contour deformation between $\Sigma$ and $\Sigma_{2}$.}
\label{fig:rh1}
\end{figure}

Let us assume that $\Sigma_{2}$ is a contour in the complex $\lambda$-plane, consisting of the four rays oriented from zero to infinity: $$\tilde\gamma_{k}:\quad\mathrm{arg}\, \lambda = \pi/6 + (k-1)\pi/3, \qquad k=1,3,4,6.$$ The contour divides the complex plane into four regions as it is shown on the right diagram of Figure \ref{fig:rh1}. Let us define the sets 
\begin{gather*}
\Omega_{r}^{1}:=\Omega_{r}\cap\tilde\Omega_{1}, \ \Omega_{r}^{2}:=\Omega_{r}\cap\tilde\Omega_{u}, \ \Omega_{r}^{6}:=\Omega_{r}\cap\tilde\Omega_{d}, \
\Omega_{l}^{4} := \Omega_{l}\cap \tilde\Omega_{4}, \\ \Omega_{l}^{5} := \Omega_{l}\cap \tilde\Omega_{d}, \ \Omega_{l}^{3} := \Omega_{l}\cap \tilde\Omega_{u},
\end{gather*}
that are shown on the left diagram of Figure \ref{fig:rh1}. Let $\Phi^{2}(\lambda)$ be a $2\times 2$ matrix valued function, defined as follows:
\begin{gather*}
\Phi^{2}(\lambda)\!:=\!\Phi(\lambda)E, \, \lambda\in \Omega_r^1,\ \Phi^{2}(\lambda)\!:=\!\Phi(\lambda)ES_1, \, \lambda\in \Omega_r^2,\ \Phi^{2}(\lambda)\!:=\!\Phi(\lambda)ES_6^{-1},\, \lambda\in \Omega_r^6\\
\Phi^{2}(\lambda):=\Phi(\lambda)\sigma_2E\sigma_2S_3^{-1}, \ \lambda\in \Omega_l^3,\quad \Phi^{2}(\lambda):=\Phi(\lambda)\sigma_2E\sigma_2, \ \lambda\in \Omega_l^4,\\
\Phi^{2}(\lambda):=\Phi(\lambda)\sigma_2ES_1\sigma_2,\ \lambda\in \Omega_l^5,\quad
\Phi^{2}(\lambda):=\Phi(\lambda), \ \lambda\in \Omega_1\cup\Omega_4\cup\Omega_u\cup\Omega_d.
\end{gather*}
At the beginning we prove the following proposition.
\begin{proposition}
The function $\Phi^{2}(\lambda)$ is a solution of the following RH problem.\\[3pt]
\noindent\makebox[5.5mm][l]{$(a)$}\parbox[t][][t]{122mm}{The function $\Phi^{2}(\lambda)$ is analytic 
for $\lambda\in\C\setminus\Sigma_{2}$.}\\[3pt]
\noindent\makebox[5.5mm][l]{$(b)$}\parbox[t][][t]{122mm}{We have the following jump relation 
$$\Phi^{2}_{+}(\lambda) = \Phi^{2}_{-}(\lambda) S_{2}(\lambda), \quad\lambda\in\Sigma_{2},$$ where the jump matrix is given by 
$$S_{2}(\lambda):=S_{k},\quad \lambda\in\tilde\gamma_{k}, \ \ k=1,3,4,6.$$}\\
\noindent\makebox[5.5mm][l]{$(c)$}\parbox[t][][t]{122mm}{The function $\Phi^{2}(\lambda)$ has the following asymptotic behavior
$$\Phi^{2}(\lambda) = (I + O(\lambda^{-1}))e^{-\theta(\lambda)\sigma_{3}}, \quad \lambda\to\infty.$$}\\
\noindent\makebox[5.5mm][l]{$(d)$}\parbox[t][][t]{122mm}{If $0<\mathrm{Re}\,\alpha<1/2$ then the function $\Phi^{2}(\lambda)$ satisfies the asymptotic relation
\begin{align*}
\Phi^{2}(\lambda)=O\begin{pmatrix}|\lambda|^{-\alpha}&|\lambda|^{-\alpha}\\[5pt]|\lambda|^{-\alpha}&|\lambda|^{-\alpha}\end{pmatrix},\quad\lambda\to0
\end{align*}}
\noindent\makebox[5.5mm][l]{}\parbox[t][][t]{121mm}{and furthermore, if $1/2<\mathrm{Re}\,\alpha\le 0$ then 
\begin{align*}
\Phi^{2}(\lambda)= O\begin{pmatrix}|\lambda|^{\alpha}&|\lambda|^{\alpha}\\[5pt]|\lambda|^{\alpha}& |\lambda|^{\alpha}\end{pmatrix},\quad\lambda\to 0.
\end{align*}}
\end{proposition}
\noindent {\em Proof.} It is not difficult to check that the function $\Phi^{2}(\lambda)$ satisfies conditions $(a)$, $(b)$ and $(c)$. The condition $(d)$ is a consequence of Lemma \ref{lem-asy-phi}. \hfill $\square$ \\

\begin{figure}[h]
\includegraphics[scale=0.65]{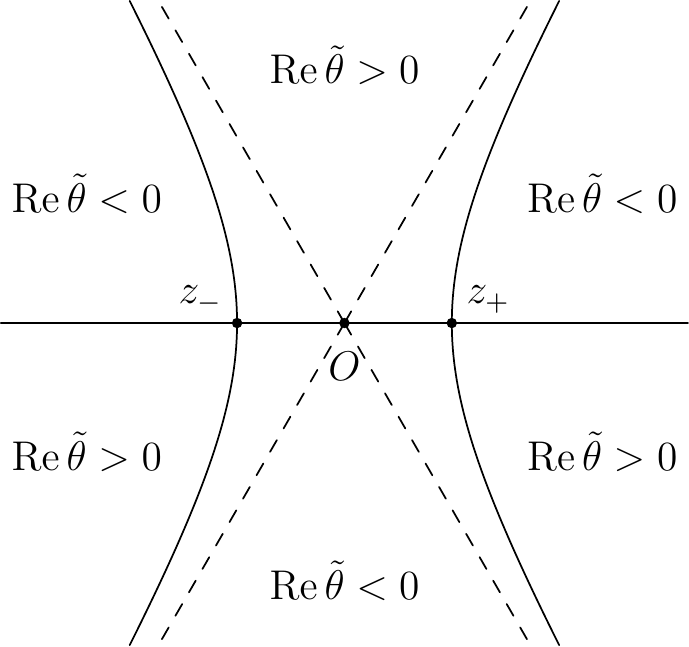}
\caption{The regions of sign changing of the function $\mathrm{Re}\,\tilde\theta(z)$. The dashed rays have directions $\exp(ik\pi/3)$ for $k=1,2,4,5$.}
\label{d26}
\end{figure}

\begin{figure}[h]
\begin{subfigure}{0.49\textwidth}
\centering
\includegraphics[scale=0.59]{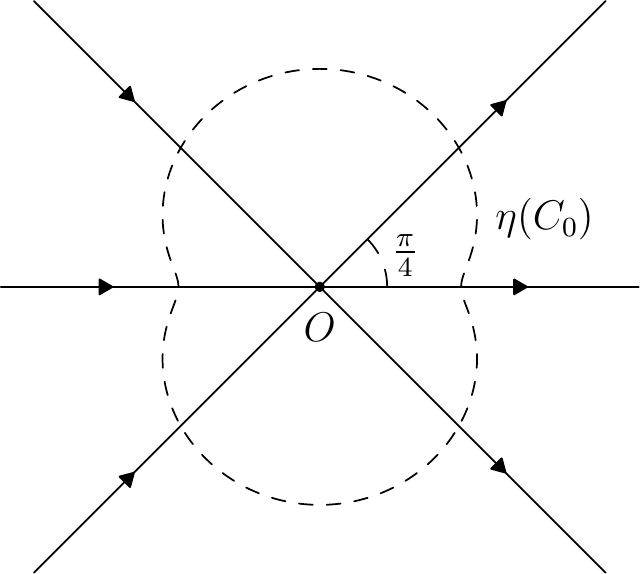}
\end{subfigure}
\begin{subfigure}{0.49\textwidth}
\centering
\includegraphics[scale=0.59]{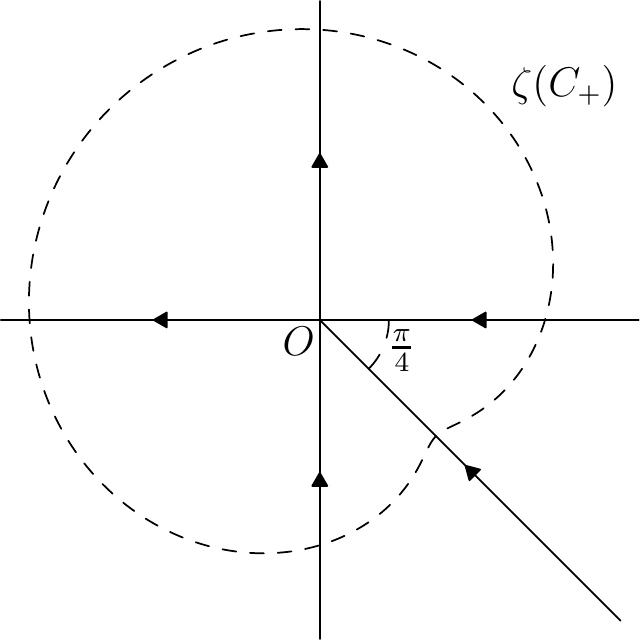}
\end{subfigure}
\caption{Left: the contour $\Sigma_{3}^{0}$ and the closed curve $\eta(C_{0})$. 
Right: the contour $\Sigma_{3}^{+}$ with $\zeta(C_{+})$.}
\label{d4bkk}
\end{figure}

\begin{figure}[h]
\includegraphics[scale=1]{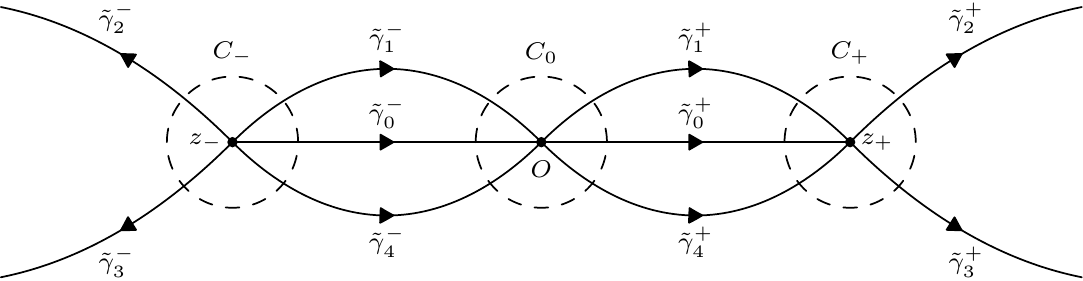}
\caption{The contour $\Sigma_{3}$ and the circles $C_{0}$, $C_{\pm}$ that are depicted by dashed lines.}
\label{fig:rh2mm}
\end{figure}

\begin{figure}[h]
\includegraphics[scale=0.85]{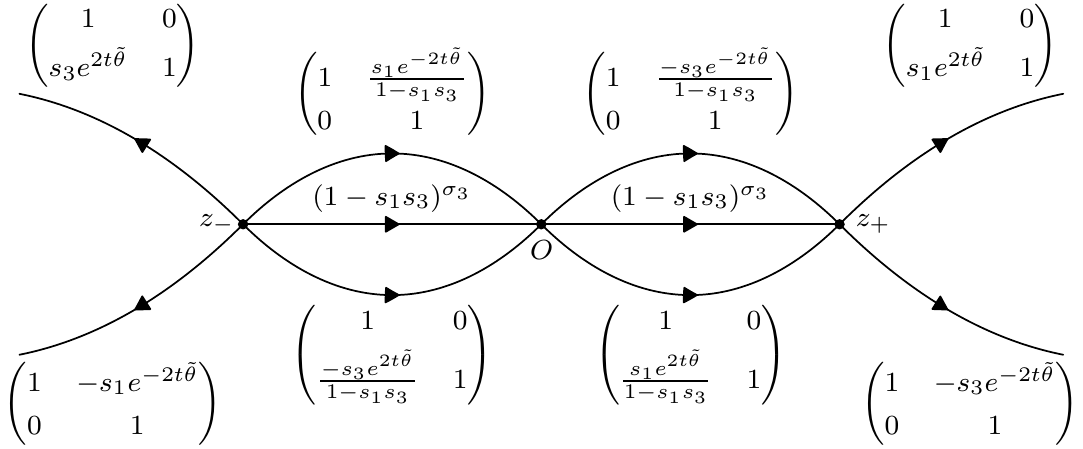}
\caption{The contour $\Sigma_{3}$ and the associated jump matrix $S_{3}$.}
\label{d4b}
\end{figure}

Let us consider the following scaling of variables 
\begin{equation*}
\lambda(z) = (-x)^{1/2}z, \ \ t(x)=(-x)^{3/2},\quad z\in \C, \ x<0
\end{equation*}
and define the function $\tilde\theta(z) := i(\frac{4}{3}z^{3} - z)$. Then it is not difficult to check that $\theta(\lambda(z),x)=t(x)\tilde\theta(z)$. 
Let us observe that $z_{\pm} = \pm 1/2$ are stationary points for the phase function $\tilde\theta$ and $\tilde\theta(\pm 1/2)=\mp i/3$. It is not difficult to check that the set of solutions of $\mathrm{Re}\,\tilde\theta(z)=0$ consists of the real axis and the curves 
$$h_{\pm}(t):=it \pm\left(t^{2}/3 +1/4\right)^{1/2},\quad t\in\R.$$
It is clear that $h_{+}$ and $h_{-}$ are asymptotic to the rays $\mathrm{arg}\,\lambda = \pm\frac{\pi}{3}$ and $\mathrm{arg}\,\lambda = \pm\frac{2\pi}{3}$, respectively (see Figure \ref{d26}).
Let us assume that $U$ is a function with values in the space $M_{2\times 2}(\C)$, which is given by the following formula
\begin{align*}
U(z,t) :=\Phi^{2}(\lambda(z),-t^{2/3})\exp(t\tilde\theta(z)\sigma_{3}).
\end{align*}
Let us consider the triangular matrices defined by
\begin{align*}
G_{2k}:= \begin{pmatrix} 1 & e^{-2t\tilde\theta(z)}s_{2k}\\[5pt] 0 & 1\end{pmatrix}, \  k=2,3 \ \ \text{and} \ \ 
G_{2k+1}:=\begin{pmatrix} 1 & 0\\[5pt] e^{2t\tilde\theta(z)}s_{2k+1} & 1\end{pmatrix}, \  k=0,1.
\end{align*}
Then it is not difficult to check that the function $U(z):=U(z,t)$ is a solution of the following scaled Riemann-Hilbert problem. \\[2pt]
\noindent\makebox[8.5mm][l]{$(D1)$}\parbox[t][][t]{119mm}{The function $U(z)$ is holomorphic for $z\in \C\setminus\Sigma_{2}$.}\\[2pt]
\noindent\makebox[8.5mm][l]{$(D2)$}\parbox[t][][t]{119mm}{For any $k=1,3,4,6$, we have the following relation $$U_{+}(z) = U_{-}(z)G_{k},\quad z\in \gamma_{k}.$$}\\
\noindent\makebox[8.5mm][l]{$(D3)$}\parbox[t][][t]{119mm}{The function $U(z)$ has the following asymptotic relation 
$$U(z) = I + O(z^{-1}),\quad z\to\infty.$$}\\
\noindent\makebox[8.5mm][l]{$(D4)$}\parbox[t][][t]{119mm}{If $0<\mathrm{Re}\,\alpha <1/2$ then the function $U(z)$ has the following asymptotic behavior
\begin{align}\label{cond-u-1}
U(z)= O\begin{pmatrix}|z|^{-\alpha}& |z|^{-\alpha}\\[5pt] |z|^{-\alpha}& |z|^{-\alpha}\end{pmatrix},\quad z\to 0
\end{align}}
\noindent\makebox[8.5mm][l]{}\parbox[t][][t]{121mm}{and furthermore, if $1/2<\mathrm{Re}\,\alpha\le 0$, then 
\begin{align}\label{cond-u-2}
U(z)= O\begin{pmatrix}|z|^{\alpha}& |z|^{\alpha}\\[5pt] |z|^{\alpha}& |z|^{\alpha}\end{pmatrix}, \quad z\to 0.
\end{align}}\\
Let us consider two auxiliary oriented graphs that are shown on Figure \ref{d4bkk}. One of them is $\Sigma_{3}^{0}$ consisting of the six rays 
\begin{align*}
\mathrm{arg}\,\lambda = 0, \quad \mathrm{arg}\,\lambda = \pi, \quad \mathrm{arg}\,\lambda = \pi/4 + k\pi/2, \quad 0\le k\le 3
\end{align*}
and the other one is $\Sigma_{3}^{+}$, formed by the five curves
\begin{align*}
\mathrm{arg}\,\lambda = 7\pi/4 \quad\text{and}\quad \mathrm{arg}\,\lambda = k\pi/2,\quad 0\le k\le 3.
\end{align*}
Let us consider the maps $\eta(z)$ and $\zeta(z)$, given by the formulas
\begin{equation}\label{et-ze}
\begin{gathered}
\eta(z):=i\tilde\theta(z) = z - 4z^3/3,\\
\zeta(z):= 2\sqrt{-\tilde\theta(z)+\tilde\theta(z_{+})}=4\sqrt{3}e^{\frac{3}{4}\pi i}\left(z-1/2\right)(z+1)^{\frac{1}{2}}/3,
\end{gathered}
\end{equation}
where the branch cut of the square root is taken such that $\mathrm{arg}\,(z-1/2)\in(-\pi,\pi)$.
Let us observe that $\eta(z)$ and $\zeta(z)$ are holomorphic functions in a neighborhood of the origin and $z_{+}$, respectively. Since $\eta'(0)\neq 0$ and $\zeta'(z_{+})\neq 0$, there is a small $\delta>0$ \label{page-11} with the property that the functions $\eta(z)$ and $\zeta(z)$ are biholomorphic on the balls $B(0,2\delta)$ and $B(z_{+},2\delta)$, respectively. \label{holo-eta} If we denote $C_{0}:=\partial B(0,\delta)$ and $C_{+}:= \partial B(z_{+},\delta)$, then both images $\eta(C_{0})$ and $\zeta(C_{+})$ are closed curves surrounding the origin (see Figure \ref{d4bkk}). 
We consider the contour $\Sigma_{3}$ (see Figure \ref{fig:rh2mm}), which consist of the curves $\tilde\gamma^{\pm}_{k}$, where $0\le k\le 4$, such that $\tilde\gamma^{\pm}_{0}$ are straight lines joining the origin with the stationary points $z_{\pm}$ and furthermore, its part contained in the ball $B(0,\delta)$ is an inverse image of the set $\Sigma_{3}^{0}\cap\eta(B(0,\delta))$ under the map $\eta$ restricted to the ball $B(0,2\delta)$. Since $\eta'(0) = 1$ it follows that the angle between the curves $\t\gamma^{+}_{1}$ and $\t\gamma^{+}_{0}$ is equal to $\pi/4$. We require also that the part of the contour $\Sigma_{3}$ contained in the ball $B(z_{+},\delta)$ is an inverse image of the set $\Sigma_{3}^{+}\cap\zeta(B(z_{+},\delta))$ under the map $\zeta$, restricted to the ball $B(z_{+},2\delta)$. Furthermore the part of the contour $\Sigma_{3}$ contained in the ball $B(z_{-},\delta)$ is taken such that it is the reflection across the origin of the set $\Sigma_{3}\cap B(z_{+},\delta)$. We also choose the unbounded components $\tilde \gamma_{2}^{\pm}$ and $\tilde \gamma_3^{\pm}$, emanating from the stationary point $z_{\pm}$, to be asymptotic to the rays $\{\mathrm{arg}\, \lambda = \pi/2\mp\pi/3\}$ and $\{\mathrm{arg}\, \lambda = 3\pi/2 \pm \pi/3\}$, respectively. 
In view of \cite[Section 3.1 and 3.2]{MR3670014}, the function $U(z)$ satisfying the RH problem (D1)\,--\,(D4) can be deformed to the function $T(z)$, which satisfies the following RH problem on a graph $\Sigma_3$. \\[2pt]
\noindent\makebox[8.5mm][l]{$(E1)$}\parbox[t][][t]{119mm}{The function $T(z)$ is holomorphic for $z\in \C\setminus\Sigma_{3}$.}\\[2pt]
\noindent\makebox[8.5mm][l]{$(E2)$}\parbox[t][][t]{119mm}{We have the jump relation 
$$T_{+}(z)= T_{-}(z)S_{3}(z),\quad z\in \Sigma_{3},$$
where the the jump matrices $S_{3}$ are presented on the Figure \ref{d4b}.}\\[2pt]
\noindent\makebox[8.5mm][l]{$(E3)$}\parbox[t][][t]{119mm}{The function $T(z)$ has the following asymptotic behavior
$$T(z) = I + O(z^{-1}),\quad z\to\infty.$$}\\
\noindent\makebox[8.5mm][l]{$(E4)$}\parbox[t][][t]{119mm}{At $z=0$, the function $T(z)$ has the same behavior as $U(z)$ in \eqref{cond-u-1} and \eqref{cond-u-2}.}\\[5pt]
From the construction of the function $T(z,t)$ we deduce that, for any $z\in\C$ such that $\mathrm{arg}\,z\in(0,\frac{\pi}{6})$ and $|z|$ is sufficiently small, the following equality holds
\begin{align}\label{rep-1}
U(z,t)=T(z,t)\begin{pmatrix}1& -\frac{s_3e^{-2t\tilde\theta(z)}}{1-s_1s_3}\\[5pt] 0& 1\end{pmatrix}
\begin{pmatrix}1& 0\\[5pt] -s_1 e^{2t\tilde\theta(z)}&1\end{pmatrix},\quad t>0.
\end{align}

\subsection{Parametrix near the origin} 
\begin{figure}[h]
\includegraphics[scale=0.8]{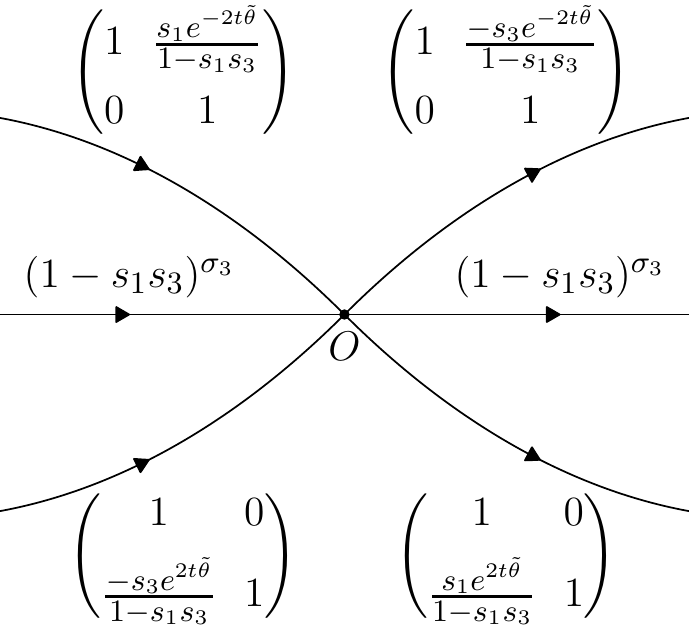}
\caption{The graph $D(0,\delta)\cap \Sigma_{3}$ for the RH problem (G1)\,--\,(G4).}
\label{d5bkk}
\end{figure}

Let us define $\nu := -(2\pi i)^{-1}\ln(1-s_1s_3)$ and consider the following function 
\begin{equation*}
N(z):=\left(\frac{z+1/2}{z-1/2}\right)^{\nu\sigma_{3}},\quad z\in \C\setminus[z_{-},z_{+}],
\end{equation*}
where the branch cut is taken such that $\mathrm{arg}\,(z\pm1/2)\in(-\pi,\pi)$. Then it is known that the function $N(z)$ satisfies the following conditions. \\[3pt]
\noindent\makebox[8.5mm][l]{$(F1)$}\parbox[t][][t]{119mm}{The function $N(z)$ is analytic on $\C\setminus[z_{-},z_{+}]$.}\\[3pt]
\noindent\makebox[8.5mm][l]{$(F2)$}\parbox[t][][t]{119mm}{If we denote $S_{D}=(1-s_{1}s_{3})^{\sigma_{3}}$, then the following jump relation holds
\begin{equation*}
N_{+}(z)=N_{-}(z)S_{D},\quad z\in [z_{-},z_{+}].
\end{equation*}}\\
\noindent\makebox[8.5mm][l]{$(F3)$}\parbox[t][][t]{119mm}{We have the asymptotic behavior $N(z) = I+O(1/z)$ as $z\to\infty$.}
\begin{remark}\label{rem-real}
A simple calculations show that $1-s_{1}s_{3}$ is a positive number. Indeed, let us observe that under condition \eqref{stokes-11}, we have
\begin{align*}
1-s_1s_3 =\cos^2(\pi\alpha)-k^2,
\end{align*}
which gives $1-s_1s_3>0$ provided \eqref{stokes2-ras} holds. On the other hand, if $\alpha,k\in i\R$, then $\alpha = i\eta$ and $k=ik_{0}$, for some $\eta,k_{0}\in\R$, and consequently
\begin{align*}
1-s_1s_3 = \cosh^{2}(\pi\eta) + k_{0}^{2}>0
\end{align*}
as desired. As a consequence we obtain that $\nu=-(2\pi i)^{-1}\ln(1-s_1s_3)$ is a purely imaginary complex number.  \hfill $\square$
\end{remark}
We are looking for the function $T_0(z)$ defined on the closed ball $D(0,\delta)$, with values in the space $M_{2\times 2}(\C)$, satisfying the following RH problem.\\[3pt]
\noindent\makebox[8.5mm][l]{$(G1)$}\parbox[t][][t]{119mm}{The function $T_0(z)$ is analytic in $D(0,\delta)\setminus\Sigma_3$.}\\[3pt]
\noindent\makebox[8.5mm][l]{$(G2)$}\parbox[t][][t]{119mm}{On the contour $D(0,\delta)\cap \Sigma_{3}$ the function $T_0(z)$ satisfies the jump conditions depicted on Figure \ref{d5bkk}.}\\[3pt]
\noindent\makebox[8.5mm][l]{$(G3)$}\parbox[t][][t]{119mm}{The function $T_0(z)$ satisfies the asymptotic relation
\begin{equation}\label{ineq-kk-ll-22}
T_0(z)N(z)^{-1} = I + O(t^{-1}), \quad t\to+\infty,
\end{equation}
uniformly for $z\in \partial D(0,\delta)$.}\\[3pt]
\noindent\makebox[8.5mm][l]{$(G4)$}\parbox[t][][t]{119mm}{At $z=0$, the function $T_0(z)$ has the same behavior as $U(z)$ in \eqref{cond-u-1} and \eqref{cond-u-2}.}
\begin{theorem}\label{lok-aprox}
The solution of the RH problem (G1)\,--\,(G4) is given by 
\begin{equation}\label{equa-p}
T_0(z):=\left\{\begin{aligned}
&E(z)\bar\Phi(t\eta(z))e^{-it\eta(z)\sigma_3} e^{-i\pi\nu\sigma_3},&&\mathrm{Im}\,z>0,\\
&E(z)\bar\Phi(t\eta(z))e^{-it\eta(z)\sigma_3} e^{i\pi\nu\sigma_3},&&\mathrm{Im}\,z<0,
\end{aligned}\right.
\end{equation}
where the function $E(z)$ is defined as follows
\begin{equation*}
E(z):=\left\{\begin{aligned}
&N(z)e^{i\pi\nu\sigma_3},&& \mathrm{Im}\, z>0,\\
&N(z)e^{-i\pi\nu\sigma_3},&& \mathrm{Im}\, z<0.
\end{aligned}\right.
\end{equation*}
\end{theorem}
\proof The fact that the function $T_0(z)$ satisfies the conditions (G1) and (G2) follows directly from the the formula \eqref{equa-p}. We show that $T_{0}(z)$ satisfies the condition $(G4)$. To this end, let us fix $t>0$ and assume that $0\le \mathrm{Re}\,\alpha<1/2$. The argument in the case of $1/2<\mathrm{Re}\,\alpha<0$ is analogous. Using the point $(c)$ of Theorem \ref{th-l} we obtain the existence of
$C_{0}>0$ and $\ve_{0}\in(0,\delta)$ such that
\begin{align}\label{eq-l-1}
|\bar\Phi(\eta)_{kl}|\le C_{0} |\eta|^{-\mathrm{Re}\,\alpha},\quad \eta\in B(0,\ve_{0}), \ \ 1\le k,l\le 2
\end{align}
Since $\eta(0)=0$ and $\eta'(0)=1$, there is $0<\ve_{1}<\ve_{0}$ such that 
\begin{equation}\label{eq-l-2}
|\eta(z)|/|z|\ge 1/2\quad\text{and}\quad |t\eta(z)|\le \ve_{0},\quad |z|\le\ve_{1}.
\end{equation}
Then, by \eqref{eq-l-1} and \eqref{eq-l-2}, for any $1\le k,l\le 2$, we have
\begin{align*}
|\bar\Phi(t\eta(z))_{kl}| \lesssim |t\eta(z)|^{-\mathrm{Re}\,\alpha} \lesssim t^{-\mathrm{Re}\,\alpha} |z|^{-\mathrm{Re}\,\alpha},\quad |z|\le\ve_{1},
\end{align*}
which implies that the following inequality holds
\begin{align}\label{ll-33}
\|\bar\Phi(t\eta(z))\|\lesssim t^{-\mathrm{Re}\,\alpha} |z|^{-\mathrm{Re}\,\alpha},\quad |z|\le\ve_{1},
\end{align}
Since the functions $E(z)$ and $e^{-it\eta(z)\sigma_3}$ are holomorphic in $D(0,2\delta)$, there is a constant $c>0$ such that
\begin{equation}\label{eq-form-3}
\|E(z)^{-1}\|\le c, \ \ \|E(z)\|\le c \ \ \text{and} \ \ \|e^{-it\eta(z)\sigma_3}\| \le c \ \ \text{for} \ \ |z|\le\delta.
\end{equation}
Combining the equality $\|e^{i\pi\nu\sigma_3}\| = \|e^{-i\pi\nu\sigma_3}\|$ with \eqref{equa-p}, \eqref{ll-33} and \eqref{eq-form-3}, for any $|z|\le\ve_{1}$, we obtain
\begin{align*}
\|T_{0}(z)\| &\le \|E(z)\|\|\bar\Phi(t\eta(z))\|\|e^{-it\eta(z)\sigma_3}\|\|e^{-i\pi\nu\sigma_3}\|\\
&\lesssim \|\bar\Phi(t\eta(z))\|\|e^{-i\pi\nu\sigma_3}\| \lesssim t^{-\mathrm{Re}\,\alpha}|z|^{-\mathrm{Re}\,\alpha},
\end{align*}
which proves that the condition $(G4)$ is satisfied. To show that the condition $(G3)$ also holds true, we observe that 
\begin{equation*}
\begin{aligned}
T^{0}(z)N(z)^{-1}&=\left\{\begin{aligned}
&E(z)\bar\Phi(t\eta(z))e^{-it\eta(z)\sigma_3} e^{-i\pi\nu\sigma_3}N(z)^{-1},&&\mathrm{Im}\,z>0,\\
&E(z)\bar\Phi(t\eta(z))e^{-it\eta(z)\sigma_3} e^{i\pi\nu\sigma_3}N(z)^{-1},&&\mathrm{Im}\,z<0,
\end{aligned}\right. \\[5pt]
& = E(z)\bar\Phi(t\eta(z))e^{-it\eta(z)\sigma_3}E(z)^{-1},
\end{aligned}
\end{equation*}
which implies that
\begin{align}\label{eq-form-1}
T_0(z)N(z)^{-1}-I = E(z)(\bar\Phi(t\eta(z))e^{-it\eta(z)\sigma_3}-I)E(z)^{-1}.
\end{align}
On the other hand, by the point $(d)$ of Theorem \ref{th-l}, there are $R,K>0$ such that  
\begin{align}\label{eq-ff-jj}
\|\bar\Phi(z)e^{-iz\sigma_3}-I\|\le K|z|^{-1},\quad |z|\ge R.
\end{align}
Since the radius $\delta>0$ is chosen so that the function $\eta(z)$ is biholomorphic on $B(0,2\delta)$ (see page \pageref{holo-eta}), we have 
$|\eta(z)|>c_{0}>0$ for $|z|=\delta$. Hence we can find $t_{0}>0$ such that $|t\eta(z)|\ge R$ for $t>t_{0}$ and $|z|=\delta$. Therefore, by the equations \eqref{eq-form-1}, \eqref{eq-ff-jj} and \eqref{eq-form-3}, for any $|z|= \delta$ and $t>t_{0}$, we have
\begin{align*}
\|T_0(z)N(z)^{-1}\!-\!I\| & \le \|E(z)\|\|(\bar\Phi(t\eta(z))e^{-it\eta(z)\sigma_3}\!-\!I)\|\|E(z)^{-1}\|\lesssim |t\eta(z)|^{-1}\lesssim t^{-1},
\end{align*}
which gives the relation \eqref{ineq-kk-ll-22} and the proof of theorem is completed. \hfill $\square$ 

\begin{remark}
The existence of solution for the problem (G1)\,--\,(G4) was proved in \cite[Section 3.5]{MR3670014} by the use of a vanishing lemma. Theorem \ref{lok-aprox} provides an explicit formula for the solution, which will be crucial in the proof of Theorem \ref{th-total}.
\end{remark}

\begin{remark}\label{l-par}
The explicit form of local parametrix around stationary points $z_{\pm}$ was constructed in \cite[Section 9.4]{MR2264522} using parabolic cylinder functions. In a consequence there is a function $T_{1}(z)$ defined on the union of the closed balls $D(z_{+},\delta)\cup D(z_{-},\delta)$ with values in the space $M_{2\times 2}(\C)$ such that the following RH problem is satisfied.\\[2pt]
\noindent\makebox[5.5mm][l]{$(a)$}\parbox[t][][t]{121mm}{The function $T_{1}(z)$ is analytic in $[D(z_{+},\delta)\cup D(z_{-},\delta)]\setminus\Sigma_{3}$.}\\[2pt]
\noindent\makebox[5.5mm][l]{$(b)$}\parbox[t][][t]{121mm}{On the contour $[D(z_{+},\delta)\cup D(z_{-},\delta)]\cap\Sigma_{3}$ the function $T_{1}(z)$ satisfies the same jump conditions as $T(z)$ (see Figure \ref{d4b}).}\\[2pt]
\noindent\makebox[5.5mm][l]{$(c)$}\parbox[t][][t]{121mm}{The function $T_{1}(z)$ satisfies the asymptotic relation
\begin{equation*}
T_{1}(z)N(z)^{-1} = I + O(t^{-1/2}),\quad  t\to+\infty,
\end{equation*}
uniformly for $z\in \partial D(z_{+},\delta)\cup\partial D(z_{-},\delta)$.}\\
\end{remark}

\subsection{Representation of solutions of the deformed RH problem} 
\begin{figure}[h]
\includegraphics[scale=0.7]{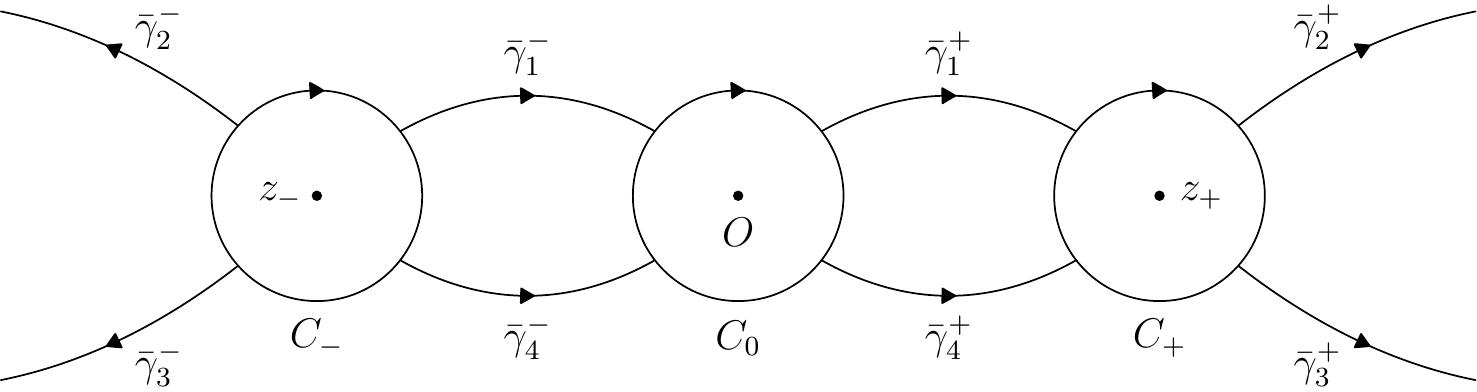}
\caption{The graph $\Sigma_{4}$ for the RH problem (H1)\,--\,(H3).}
\label{d5bkkll}
\end{figure}

We intend to use Theorem \ref{lok-aprox} to derive Proposition \ref{prop-p-l}, which provides information about the asymptotic behavior of the function $P(x)$ as $x\to-\infty$. Let us assume that $R(z)$ is a function given by
\begin{gather*}
R(z):=T(z)T_{1}(z)^{-1}, \ z\in [D(z_{+},\delta)\cup D(z_{-},\delta)]\setminus\Sigma_{3}\\
R(z):=T(z)T_{0}(z)^{-1}, \ z\in D(z_{0},\delta)\setminus\Sigma_{3}\\
R(z):=T(z)N(z)^{-1}, \ z\in\C\setminus(D(z_{+},\delta)\cup D(z_{-},\delta)\cup D(0,\delta)\cup\Sigma_{3})
\end{gather*}
and let $\Sigma_{4}$ be the contour depicted on the Figure \ref{d5bkkll}, which consists of circles $C_{\pm}$ and $C_{0}$ (see page \pageref{page-11}) and the parts $\bar\gamma^{\pm}_{k}$ of the curves $\tilde\gamma^{\pm}_{k}$ lying outside the set $C_{+}\cup C_{-}\cup C_{0}$ (see Figure \ref{fig:rh2mm}). Then it is clear that the function $R(z)$ is the solution of the following Riemann-Hilbert problem.\\[2pt]
\noindent\makebox[8.5mm][l]{$(H1)$}\parbox[t][][t]{119mm}{The function $R(z)$ is analytic in $\C\setminus\Sigma_{4}$.}\\[2pt]
\noindent\makebox[8.5mm][l]{$(H2)$}\parbox[t][][t]{119mm}{The following jump condition holds 
\begin{align*}
R_{+}(z) = R_{-}(z)S_{4}(z),\quad z\in\Sigma_{4},
\end{align*}
where the jump matrix is given by}\\[2pt]
\noindent\makebox[8.5mm][l]{}\parbox[t][][t]{119mm}{
\begin{gather*}
\hspace{-30pt}S_{4}(z):=T_{1}(z)N(z)^{-1},\ z\in\partial D(z_{+},\delta)\cup \partial D(z_{-},\delta),\\
\hspace{-30pt}S_{4}(z):=T_{0}(z)N(z)^{-1}, \ z\in \partial D(0,\delta)\\
\hspace{-30pt}S_{4}(z):=N(z)S_{3}(z)N(z)^{-1}, \ z\in\Sigma_{4}\setminus(\partial D(z_{+},\delta)\cup \partial D(z_{-},\delta)\cup \partial D(0,\delta))
\end{gather*}}\\
\noindent\makebox[8.5mm][l]{$(H3)$}\parbox[t][][t]{119mm}{We have the following asymptotic behavior 
\begin{align*}
R(z)=I + O(1/z),\quad z\to\infty.
\end{align*}}\\

In the following lemma we provide useful estimates for the jump matrix $S_{4}$ in the space $L^{p}(\Gamma_{1})$, where $1\le p\le +\infty$.
\begin{lemma}\label{lem-est-3}
Given $1\le p<+\infty$, we have the following asymptotic relation 
\begin{align}\label{ineq-11aa}
\|S_{4}-I\|_{(L^{p}\cap L^{\infty})(\Sigma_{4})} = O(t^{-1/2}),\quad t\to\infty.
\end{align}
\end{lemma}
\proof By Theorem \ref{lok-aprox} and Remark \ref{l-par}, we infer that 
\begin{align*}
\|S_{4} - I\|_{L^{\infty}(C_{+}\cup C_{-})} & = \|T_{1}N^{-1} - I\|_{L^{\infty}(C_{+}\cup C_{-})} = O(t^{-1/2}), \quad t\to\infty,\\
\|S_{4} - I\|_{L^{\infty}(C_{0})} & = \|T_{0}N^{-1} - I\|_{L^{\infty}(C_{0})} = O(t^{-1}), \quad t\to\infty.
\end{align*}
Furthermore, by Theorem \ref{lok-aprox}, we have the following asymptotic
\begin{equation*}
\|S_{4} - I\|_{L^{\infty}(C_{0})} = \|T_{0}N^{-1} - I\|_{L^{\infty}(C_{0})} = O(t^{-1}), \quad t\to\infty.
\end{equation*}
In particular, there is $t_{0}>0$ such that
\begin{equation}
\begin{aligned}\label{ee-m}
&\|S_{4} - I\|_{L^{\infty}(C_{+}\cup C_{-})} \lesssim t^{-1/2}, \quad \|S_{4} - I\|_{L^{\infty}(C_{0})} \lesssim t^{-1},\quad t\ge t_{0},
\end{aligned}
\end{equation}
which implies that, for any $t>t_{0}$, the following estimates hold
\begin{align}\label{ee-4-gg}
\|S_{4} - I\|^{p}_{L^{p}(C_{0})} = \int_{C_{0}}\|S_{4}(z) - I\|^{p}\,|dz|\lesssim \int_{C_{0}}t^{-p}\,|dz|\sim t^{-p}
\end{align}
and furthermore, for any $t>t_{0}$, we have
\begin{align}\label{ee-4}
\|S_{4} - I\|^{p}_{L^{p}(C_{+}\cup C_{-})} = \int_{C_{\pm}}\|S_{4}(z) - I\|^{p}\,|dz|\lesssim \int_{C_{+}\cup C_{-}}t^{-p/2}\,|dz|\sim t^{-p/2}.
\end{align}
Let us denote $\Sigma_{4}':=\Sigma_{4}\setminus [C_{+}\cup C_{-}\cup C_{0}]$. By the definition of $N(z)$ and the choice of the component curves of $\Sigma_{3}$, there is a constant $C>0$ such that
\begin{align*}
\|N(z)\|\le C \quad \text{and} \quad \|N(z)^{-1}\|\le C \quad \text{for} \ \ z\in \Sigma_{4}'.
\end{align*}
Consequently, for any $z\in \Sigma_{4}'$, we have
\begin{equation}
\begin{aligned}\label{eq-r-t}
&\|S_{4}(z) - I\| = \|N(z)[S_{3}(z) - I]N(z)^{-1}\|\\
& \qquad\le \|N(z)\|\|S_{3}(z) - I\|\|N(z)^{-1}\|\le C^{2}\|S_{3}(z) - I\|.
\end{aligned}
\end{equation}
Since the curve $\bar\gamma^{+}_{2}$ is asymptotic to the ray $\{se^{i\pi/6} \ | \ s>0\}$, there is a smooth function $h:[a,+\infty)\to\R$, where $a>0$, 
satisfying asymptotic condition
\begin{equation*}
h(s)/s\to \sqrt{3}/3,\quad s\to+\infty
\end{equation*}
such that the map $\bar\gamma^{+}_{2}:[a,+\infty)\to\C$ given by the formula $$\bar\gamma^{+}_{2}(s):=s + ih(s),\quad s\ge a,$$ 
is a parametrization of the curve $\bar\gamma^{+}_{2}$. Let us take sufficiently small $\ve_{0}>0$ such that 
\begin{equation}\label{eq-ineq-55}
4(\sqrt{3}/3+\ve_{0})^{3}/3- 4(\sqrt{3}/3-\ve_{0})<0
\end{equation}
and observe that, there is $a_{0}>a$ with the property that, for any $s>a_{0}$, we have
\begin{align*}
&\mathrm{Re}\,\tilde\theta(s+ih(s))= 4 h(s)^{3}/3- 4s^{2}h(s)+h(s)\\
&\quad\le (4(\sqrt{3}/3+\ve_{0})^{3}/3- 4(\sqrt{3}/3-\ve_{0}))s^{3} + (\sqrt{3}/3+\ve_{0})s. 
\end{align*}
Hence, in view of \eqref{eq-ineq-55}, there is $a_{1}>a_{0}>0$ such that 
\begin{align}\label{eq-hh-22}
\mathrm{Re}\,\tilde\theta(\bar\gamma^{+}_{2}(s)) = \mathrm{Re}\,\tilde\theta(s+ ih(s))\le -s,\quad s\ge a_{1}.
\end{align}
On the other hand, using the sign changing properties for the function $\mathrm{Re}\,\t\theta(z)$ (see Figure \ref{d26}), we obtain the existence of constants $c_{0}>0$ such that 
\begin{align*}
\mathrm{Re}\,\tilde\theta(\bar\gamma^{+}_{2}(s))=\mathrm{Re}\,\tilde\theta(s+ ih(s))\le -c_{0},\quad s\in [a,a_{1}].
\end{align*}
In particular, combining this inequality with \eqref{eq-hh-22} we infer that 
\begin{align}\label{kk-mm-11-22}
\mathrm{Re}\,\tilde\theta(\bar\gamma^{+}_{2}(s))\le -m,\quad s\ge a,
\end{align}
where we define $m:=\min(c_{0},a_{1})$. Let us observe that the form of the jump matrix $S_{3}$ on the curve $\tilde\gamma_{2}^{+}$ (see Figure \ref{d4b}), gives us 
\begin{align}\label{form-s-3}
\|S_{3}(\bar\gamma^{+}_{2}(s)) - I\| = |s_{1}|e^{2t\mathrm{Re}\,\tilde\theta(\bar\gamma^{+}_{2}(s))},\quad s>a, \ t>0,
\end{align}
which together with \eqref{eq-r-t} and \eqref{kk-mm-11-22}, give
\begin{align}\label{kk-ll-mm}
\|I - S_{4}\|_{L^{\infty}(\bar\gamma_{2}^{+})}\lesssim \|I - S_{3}\|_{L^{\infty}(\bar\gamma_{2}^{+})}\lesssim \sup_{s\ge a} e^{2t\tilde\theta(\bar \gamma^{+}_{2}(s))}\le e^{-2mt},\quad t>0.
\end{align}
On the other hand, using \eqref{eq-r-t}, \eqref{eq-hh-22}, \eqref{kk-mm-11-22} and \eqref{form-s-3}, for any $t>0$, we have 
\begin{equation}
\begin{aligned}\label{kk-ll-mm-22}
&\|I - S_{4}\|^{p}_{L^{p}(\bar\gamma_{2}^{+})}\lesssim \|I - S_{3}\|^{p}_{L^{p}(\bar\gamma_{2}^{+})}\lesssim \int_{a}^{\infty}|e^{2t\tilde\theta(\bar\gamma^{+}_{2}(s))}|^{p}|(\bar\gamma^{+}_{2})'(s)|\,ds\\
&\qquad \lesssim \int_{a}^{a_{1}}e^{2pt\mathrm{Re}\,\tilde\theta(\bar\gamma_{2}^{+}(s))}\,ds+
\int_{a_{1}}^{\infty}e^{2pt\mathrm{Re}\,\tilde\theta(\bar\gamma_{2}^{+}(s))}\,ds\\
&\qquad = \int_{a}^{a_{1}}e^{-2pc_{0}t}\,ds+
\int_{a_{1}}^{\infty}e^{-2pts}\,ds =  (a_{1}-a)e^{-2pc_{0}t} + (2pt)^{-1}e^{-2pta_{1}}.
\end{aligned}
\end{equation}
Combining \eqref{kk-ll-mm} and \eqref{kk-ll-mm-22} we deduce that 
\begin{align}\label{est-inf-1}
\|S_{4}-I\|_{(L^{p}\cap L^{\infty})(\bar\gamma_{2}^{+})} \lesssim e^{-ct},\quad t>0,
\end{align}
where $c>0$ is a constant. Let us assume that $\bar\gamma^{+}_{1}:[0,1]\to\C$ is a parametrization of the curve $\bar\gamma^{+}_{1}$. Applying the sign changing diagram from Figure \ref{d26} once again, we obtain the existence of $m_{1}>0$ such that the following inequality holds
\begin{align}\label{ineq-kk-hh-23}
\mathrm{Re}\,\tilde\theta(\bar\gamma^{+}_{1}(s))\ge m_{1},\quad s\in[0,1].
\end{align}
Using the form of the jump matrix $S_{3}$ on the curve $\tilde\gamma_{1}^{+}$ (see Figure \ref{d4b}), we obtain
\begin{align*}
\|S_{3}(\bar\gamma^{+}_{1}(s)) - I\| = |s_{3}||1-s_{1}s_{3}|^{-1}e^{-2t\mathrm{Re}\,\tilde\theta(\bar\gamma^{+}_{1}(s))},\quad 
s\in[0,1], \ t>0,
\end{align*}
which together with \eqref{eq-r-t} and \eqref{ineq-kk-hh-23} imply that
\begin{equation*}
\|I - S_{4}\|^{p}_{L^{p}(\bar\gamma_{1}^{+})}\lesssim \|I - S_{3}\|^{p}_{L^{p}(\bar\gamma_{1}^{+})}\lesssim \int_{0}^{1}|e^{2t\tilde\theta(\bar\gamma^{+}_{1}(s))}|^{p}|\bar\gamma^{+}_{1}(s)'|\,ds\lesssim e^{-2pm_1t}
\end{equation*}
for $t>0$ and furthermore
\begin{align*}
\|I - S_{4}\|_{L^{\infty}(\bar\gamma_{1}^{+})}\lesssim \|I - S_{3}\|_{L^{\infty}(\bar\gamma_{1}^{+})}\lesssim\sup_{s\ge 0} |e^{2t\tilde\theta(\bar\gamma^{+}_{1}(s))}|\le e^{-2m_1 t},\quad t>0.
\end{align*}
Consequently there is a constant $c>0$ such that 
\begin{align}\label{est-inf-2}
\|S_{4}-I\|_{(L^{p}\cap L^{\infty})(\bar\gamma_{1}^{+})} = e^{-ct},\quad t>0.
\end{align}
Proceeding as above way we can obtain the estimates \eqref{est-inf-1} and \eqref{est-inf-2} for the remaining components $\bar\gamma^{\pm}_{k}$, where $k=2,3,4$. This leads to the inequalities
\begin{align}\label{ee-5}
\|S_{4} - I\|_{L^{\infty}(\Sigma_{4}')}\lesssim e^{-ct},\quad \|S_{4} - I\|_{L^{p}(\Sigma_{4}')}\lesssim e^{-ct},\quad t> 0,
\end{align}
where $c>0$ is some constant. Combining \eqref{ee-m}, \eqref{ee-4-gg}, \eqref{ee-4} and \eqref{ee-5} yields the desired relation \eqref{ineq-11aa} and the proof of lemma is completed. \hfill $\square$\\

We proceed to the proof of the following proposition concerning asymptotic behavior of the function $R(0,t)$ as $t\to+\infty$.

\begin{proposition}\label{prop-final}
There is $t_{1}>0$ such that, for any $t>t_{1}$, the RH problem (H1)\,--\,(H3) admits a unique solution $R(z,t)$ with the property that 
\begin{align}\label{aaa1}
\|R(0,t) - I\| = O(t^{-1/2}),\quad t\to+\infty.
\end{align}
\end{proposition}
\proof From Lemma \ref{lem-est-3} it follows that there are $t_{0}>0$ such that 
\begin{align}\label{est-t9c}
\|S_{4} - I\|_{(L^{1}\cap L^{2}\cap L^{\infty})(\Sigma_{4})}\lesssim t^{-1/2}, \quad t>t_{0}.
\end{align}
Let us assume that $\RR:L^{2}_{I}(\Sigma_{4})\to L^{2}_{I}(\Sigma_{4})$ is a linear map given by the formula
\begin{equation*}
\RR(\rho):= \CC_{-}(\rho(S_{4} - I)),\quad \rho\in L^{2}_{I}(\Sigma_{4}),
\end{equation*}
where $\mathcal{C}_-$ represents the Cauchy operator on the contour $\Sigma_{4}$. 
Let us take $\rho\in L^{2}_{I}(\Sigma_{4})$ with $\rho = \rho_{0}+\rho_{\infty}$, where $\rho_{0}\in L^{2}(\Sigma_{4})$ and $\rho_{\infty}\in M_{2\times 2}(\C)$. Then, by the linearity of the Cauchy operator, we have
\begin{equation*}
\RR(\rho)= \CC_{-}((\rho-a)(S_{4} - I))+ \CC_{-}(\rho_{\infty}(S_{4} - I)).
\end{equation*}
Therefore the following estimates hold
\begin{equation}
\begin{aligned}\label{est-t5c}
\|\RR(\rho)\|_{L^{2}(\Sigma_{4})}&\lesssim \|\rho_{0}(S_{4} - I)\|_{L^{2}(\Sigma_{4})}+
\|\rho_{\infty}\|\|(S_{4} - I)\|_{L^{2}(\Sigma_{4})}\\
&\le \|S_{4}-I\|_{(L^{2}\cap L^{\infty})(\Sigma_{4})}\left(\|\rho_{0}\|_{L^{2}(\Sigma_{4})}+\|\rho_{\infty}\|\right)\\
&= \|\rho\|_{L^{2}_{I}(\Sigma_{4})}\|S_{4}-I\|_{(L^{2}\cap L^{\infty})(\Sigma_{4})}.
\end{aligned}
\end{equation}
Combining the inequalities \eqref{est-t9c} and \eqref{est-t5c}, we have
\begin{align}\label{ee10c}
\|\RR(\rho)\|_{L^{2}_{I}(\Sigma_{4})}\lesssim t^{-1/2}\|\rho\|_{L^{2}_{I}(\Sigma_{4})}, \quad t>t_{0},
\end{align}
which, in particular, implies that 
\begin{align}\label{ee11c}
\|\RR I\|_{L^{2}_{I}(\Sigma_{4})} \lesssim t^{-1/2}, \quad t>t_{0}.
\end{align}
Furthermore the inequality \eqref{ee10c} shows the existence of $t_{1}>t_{0}$ such that 
\begin{equation}\label{hh-kk-11}
\|\RR\|_{L^{2}_{I}(\Sigma_{4})}<1/2,\quad t>t_{1},
\end{equation}
which implies that the equation $\rho - \RR(\rho) = I$ has a unique solution $\rho\in L^{2}_{I}( \Sigma_{4})$, given by the series $\rho := \sum_{i=0}^{\infty} \RR^{i}I$, which is convergent in the space $L^{2}_{I}(\Sigma_{4})$. Therefore the inequalities \eqref{ee11c} and \eqref{hh-kk-11} yield
\begin{equation}
\begin{aligned}\label{ee15c}
\|\rho - I\|_{L^{2}_{I}(\Sigma_{4})}\le \|\RR I\|_{L^{2}( \Sigma_{4})}\sum_{i=0}^{\infty} \|\RR\|^{i}_{L^{2}_{I}( \Sigma_{4})}
\lesssim\|\RR I\|_{L^{2}(\Sigma_{4})} \lesssim t^{-\frac{1}{2}}, \quad t>t_{1}.
\end{aligned}
\end{equation}
It is known that the solutions of the RH problem (H1)\,--\,(H3) is represented by 
\begin{align*}
R(\lambda)= I +\frac{1}{2\pi i}\int_{\Sigma_{4}} \frac{\rho(\xi)(S_{4}(\xi)-I)}{\xi-\lambda}\,d\xi, \quad \lambda\not\in \Sigma_{4}.
\end{align*}
Therefore applying the H\"older inequality we obtain the estimates
\begin{align*}
&\|R(0) - I\| \le \frac{1}{2\pi}\int_{\Sigma_{4}} \frac{\|\rho(\xi)(S_{4}(\xi)-I)\|}{|\xi|}\,|d\xi|\le
\delta^{-1}\int_{\Sigma_{4}}\|\rho(\xi)(S_{4}(\xi)-I)\|\,|d\xi|\\
&\quad\lesssim \int_{\Sigma_{4}}\|(\rho(\xi)-I)(S_{4}(\xi)-I)\|\,|d\xi|+\int_{\Sigma_{4}}\|S_{4}(\xi)-I\|\,|d\xi| \\
&\quad\le \|\rho-I\|_{L^{2}( \Sigma_{4})}\|S_{4}-I\|_{L^{2}(\Sigma_{4})}+\|S_{4}-I\|_{L^{1}(\Sigma_{4})},\quad t> t_{1}.
\end{align*}
Combining this with \eqref{est-t9c} and \eqref{ee15c} gives
\begin{align*}
&\|R(0) - I\|\lesssim t^{-1/2}t^{-1/2}+ t^{-1/2}= t^{-1}+t^{-1/2}\lesssim t^{-1/2},\quad t>t_{1}
\end{align*}
and the proof of proposition is completed. \hfill $\square$\\

In the following proposition we derive the representation formula for the function $Z(\lambda,x)$, provided $x<0$ and $|x|$ is sufficiently large.
\begin{proposition}\label{prop-con}
There is $x_{-}<0$ such that, for any $x<x_{-}$, there is $\ve>0$ with the property that, for any $\lambda\in B(0,\ve)$ with $\mathrm{arg}\,\lambda=\pi/12$, we have
\begin{align}\label{eq-hh-kk-22}
Z(\lambda,x)=R(z)E(z)\hat\Phi^{0}(e^{2\pi i}it\eta(z))D K z^{-\alpha\sigma_3}e^{t\tilde\theta(z)\sigma_{3}}t^{-\frac{\alpha\sigma_3}{3}},
\end{align}
where $z=(-x)^{-1/2}\lambda$, $t=(-x)^{3/2}$ and the matrix $K$ is defined as
\begin{align}\label{eq-k}
K:=(1-s_1s_3)^{-\frac{1}{2}}E\hat S_2\begin{pmatrix}1& -s_{3}\\[10pt] -s_1& 1\end{pmatrix}E^{-1}.
\end{align}
\end{proposition}
\proof By Proposition \ref{prop-final}, there is $t_{1}>0$ such that, for any $t>t_{1}$, we have
\begin{align}\label{ee-dd-22}
T(z,t) = R(z,t)T_{0}(z,t),\quad |z|<\delta.
\end{align}
Let us define $x_{-}:=-t_{1}^{2/3}$ and fix $x<x_{-}$. Clearly $t=(-x)^{3/2}>t_{1}$. Let us take $\ve>0$ such that \eqref{rep-1} holds for $z=\lambda t^{-\frac{1}{3}}$, where $|\lambda|<\ve$ with $\mathrm{arg}\,\lambda = \pi/12$. As we have seen in the construction of the contour $\Sigma_{3}$ (see Section \ref{sect-constr-sigma-t}), the angle between the curves $\t\gamma^{+}_{1}$ and $\t\gamma^{+}_{0}$ at the point $O$ is equal to $\pi/4$. Since $\eta'(0)=1$ we can decrease $\ve>0$ if necessary such that
\begin{align}\label{et-omeg}
\eta(z) = \eta(\lambda t^{-\frac{1}{3}})\in \bar\Omega^{2}_{r}\cap \{z\in\C \ | \ \mathrm{Im}\, z>0\}, \quad |\lambda|\le\ve, \ \mathrm{arg}\,\lambda= \pi/12,
\end{align}
where the set $\bar\Omega^{2}_{r}$ is shown on Figure \ref{d5b}. Therefore, if $|\lambda|\le\ve$ and $\mathrm{arg}\,\lambda= \pi/12$, then
\begin{align*}
Z(\lambda,x)&=\Phi(\lambda,x)e^{\theta(\lambda)\sigma_{3}}\lambda^{-\alpha\sigma_3}=
\Phi^{2}(\lambda,x)E^{-1}\lambda^{-\alpha\sigma_3}e^{\theta(\lambda)\sigma_{3}}\\
&=\Phi^{2}(t^{1/3}z,-t^{2/3})E^{-1}z^{-\alpha\sigma_3} e^{t\tilde\theta(z)\sigma_{3}} t^{-\frac{\alpha\sigma_3}{3}}\\
& =U(z,t)e^{-t\tilde\theta(z)\sigma_{3}}E^{-1}z^{-\alpha\sigma_3} e^{t\tilde\theta(z)\sigma_{3}} t^{-\frac{\alpha\sigma_3}{3}}  \\
& =T(z,t)\begin{pmatrix}1& \hspace{-5pt}-\frac{s_3e^{-2t\tilde\theta(z)}}{1-s_1s_3}\\[5pt] 0&\hspace{-5pt} 1\end{pmatrix}
\begin{pmatrix}1&\hspace{-5pt}0\\[5pt] -s_1 e^{2t\tilde\theta(z)}&\hspace{-5pt}1   \end{pmatrix}e^{-t\tilde\theta(z)\sigma_{3}}E^{-1}z^{-\alpha\sigma_3}e^{t\tilde\theta(z)\sigma_{3}}t^{-\frac{\alpha\sigma_3}{3}}.
\end{align*}
It is not difficult to check that 
\begin{align*}
\begin{pmatrix}1&\hspace{-5pt} -\frac{s_3e^{-2t\tilde\theta}}{1-s_1s_3}\\[5pt] 0& \hspace{-5pt}1\end{pmatrix}\begin{pmatrix}1&\hspace{-5pt}0\\[5pt] -s_1 e^{2t\tilde\theta(z)}&\hspace{-5pt}1\end{pmatrix}e^{-t\tilde\theta(z)\sigma_3} = e^{-t\tilde\theta(z)\sigma_3}
\begin{pmatrix}1&\hspace{-5pt} -\frac{s_3}{1-s_1s_3}\\[5pt] 0& \hspace{-5pt}1\end{pmatrix}S_{1}^{-1},
\end{align*}
which together with \eqref{ee-dd-22}, give
\begin{align}\label{eq-c1}
Z(\lambda,x)=R(z)T_0(z)e^{-t\tilde\theta(z)\sigma_{3}}\begin{pmatrix}1&\hspace{-5pt} -\frac{s_3}{1-s_1s_3}\\[5pt] 0& \hspace{-5pt}1\end{pmatrix}S_{1}^{-1}E^{-1}z^{-\alpha\sigma_3}e^{t\tilde\theta(z)\sigma_{3}}t^{-\frac{\alpha\sigma_3}{3}}.
\end{align} 
In view of \eqref{et-ze}, \eqref{equa-p} and \eqref{et-omeg}, we obtain
\begin{equation*}
\begin{aligned}
&T_0(z)=E(z)\bar\Phi(t\eta(z))e^{-it\eta(z)\sigma_3} e^{-i\pi\nu\sigma_3}
=E(z)\tilde\Phi(t\eta(z))e^{-it\eta(z)\sigma_3} e^{-i\pi\nu\sigma_3}\\
&=E(z)\check\Phi(it\eta(z))e^{-it\eta(z)\sigma_3}e^{-i\pi\nu\sigma_3}
=E(z)\hat\Phi(it\eta(z))\sigma_{2}E\sigma_{2}e^{t\tilde\theta(z)\sigma_3}e^{-i\pi\nu\sigma_3}.
\end{aligned}
\end{equation*}
Combining this with Lemma \ref{lem-eq} and Remark \ref{rem-real} yields
\begin{equation}
\begin{aligned}\label{eq-c5}
T_0(z)&=E(z)\hat\Phi^{0}(e^{2\pi i}it\eta(z))DE\hat S_{2}e^{t\tilde\theta(z)\sigma_3} e^{-i\pi\nu\sigma_3}\\
&=E(z)\hat\Phi^{0}(e^{2\pi i}it\eta(z))DE\hat S_{2}e^{t\tilde\theta(z)\sigma_3} (1-s_{1}s_{3})^{\frac{\sigma_{3}}{2}}.
\end{aligned}
\end{equation}
Substituting \eqref{eq-c5} into \eqref{eq-c1}, we obtain
\begin{align}\label{eq-rr-11}
Z(\lambda,x)=R(z)E(z)\hat\Phi^{0}(e^{2\pi i}it\eta(z))DE\hat S_{2}WE^{-1}z^{-\alpha\sigma_3}e^{t\tilde\theta(z)\sigma_{3}}t^{-\frac{\alpha\sigma_3}{3}},
\end{align}
where we define
\begin{align*}
W:= e^{t\tilde\theta(z)\sigma_3}(1-s_{1}s_{3})^{\frac{\sigma_{3}}{2}}e^{-t\tilde\theta(z)\sigma_{3}}\begin{pmatrix}1&\hspace{-5pt} -\frac{s_3}{1-s_1s_3}\\[5pt] 0& \hspace{-5pt}1\end{pmatrix} S_{1}^{-1}.
\end{align*}
On the other hand the following equalities hold
\begin{equation}
\hspace{-10pt}\begin{aligned}\label{eq-c3}
W& =\begin{pmatrix}(1-s_1s_3)^{1/2}&  0 \\[5pt] 0& (1-s_1s_3)^{-1/2}\end{pmatrix} \begin{pmatrix}1&\hspace{-5pt} -\frac{s_3}{1-s_1s_3}\\[5pt] 0& \hspace{-5pt}1\end{pmatrix}
\begin{pmatrix}1& 0 \\[5pt] -s_1& 1\end{pmatrix}\\[7pt]
& =\begin{pmatrix}(1-s_1s_3)^{1/2}& \hspace{-5pt} \frac{-s_3}{(1-s_1s_3)^{1/2}}\\[5pt] 0& \hspace{-5pt}\frac{1}{(1-s_1s_3)^{1/2}}\end{pmatrix} 
\begin{pmatrix}1& \hspace{-5pt}0 \\[5pt] -s_1& \hspace{-5pt}1\end{pmatrix}
=\frac{1}{(1-s_1s_3)^{1/2}} \begin{pmatrix}1& \hspace{-5pt}-s_3\\[5pt] -s_1& \hspace{-5pt}1\end{pmatrix}.
\end{aligned}
\end{equation}
Therefore, combining \eqref{eq-c3} with \eqref{eq-rr-11}, we obtain \eqref{eq-hh-kk-22}, which completes the proof of proposition. \hfill $\square$\\

In the following proposition we perform some calculations to present the matrix $K$ in a simpler diagonal form. 

\begin{proposition}\label{prop-diag}
The matrix $K$ defined in \eqref{eq-k} is expressed by the formula
\begin{align*}
K=\frac{1}{(1-s_1s_3)^{1/2}} \begin{pmatrix}\cos(\pi\alpha)-k&0\\[5pt]0&\cos(\pi\alpha)+k\end{pmatrix} e^{-i\pi\alpha\sigma_3}.
\end{align*}
In particular the matrix $K$ is diagonal.
\end{proposition}
\proof Let us observe that, by the constraint condition \eqref{cc1} and \eqref{matrix-s}, we have
\begin{align*}
\hat S_2\begin{pmatrix}1&\hspace{-5pt} -s_{3}\\[5pt] -s_1& \hspace{-5pt}1\end{pmatrix}&=\begin{pmatrix}1& -s_3\\[5pt] -2\sin(\pi\alpha)-s_1& 1+2s_3\sin(\pi\alpha)\end{pmatrix}
=\begin{pmatrix}1& -s_3\\[5pt] s_3& 1+2s_3\sin(\pi\alpha)\end{pmatrix}.
\end{align*}
Since the matrix $E$ is given by \eqref{eq-e}, its inverse is of the form
\begin{align*}
E^{-1}=\frac{1}{-2i\cos(\pi\alpha)}\begin{pmatrix}-ie^{i\pi\alpha}& -ie^{-i\pi\alpha}\\[5pt] -1& 1\end{pmatrix}\begin{pmatrix}p&0\\[5pt] 0&q\end{pmatrix}^{-1}.
\end{align*}
Therefore, if we define 
\begin{align*}
A:=\begin{pmatrix}A_{11}&A_{12}\\[3pt]A_{21}&A_{22}\end{pmatrix}=\begin{pmatrix}1& ie^{-i\pi\alpha}\\[5pt] 1& -ie^{i\pi\alpha}\end{pmatrix}\begin{pmatrix}1& -s_3\\[5pt] s_3& 1+2s_3\sin(\pi\alpha)\end{pmatrix} \begin{pmatrix}-ie^{i\pi\alpha}& -ie^{-i\pi\alpha}\\[5pt] -1& 1\end{pmatrix}.
\end{align*}
then the matrix $K$ has the following form 
\begin{align}\label{matrix-kk}
K&= -\frac{(1-s_1s_3)^{-\frac{1}{2}}}{2i\cos(\pi\alpha)}\begin{pmatrix}p&0\\[3pt] 0& q\end{pmatrix}
\begin{pmatrix}A_{11}&A_{12}\\[3pt]A_{21}&A_{22}\end{pmatrix}\begin{pmatrix}p& 0\\[3pt] 0& q\end{pmatrix}^{-1}.
\end{align}
Observe that after multiplication, we obtain 
\begin{align*}
K_{0}&= \begin{pmatrix}1+is_3e^{-i\pi\alpha} & -s_3 + ie^{-i\pi\alpha}+2is_3e^{-i\pi\alpha}\sin(\pi\alpha)\\[5pt]
1-is_3e^{i\pi\alpha} & -s_3 - ie^{i\pi\alpha}-2is_3e^{i\pi\alpha}\sin(\pi\alpha)\end{pmatrix}
\begin{pmatrix}-ie^{i\pi\alpha}& -ie^{-i\pi\alpha}\\[5pt] -1& 1\end{pmatrix},
\end{align*}
which in turn implies that the entries are given by 
\begin{align*}
A_{11}&=-ie^{i\pi\alpha}(1+is_3e^{-i\pi\alpha})+s_3 - ie^{-i\pi\alpha}-2is_3e^{-i\pi\alpha}\sin(\pi\alpha),\\
A_{12}& = -ie^{-i\pi\alpha}(1+is_3e^{-i\pi\alpha})-s_3 + ie^{-i\pi\alpha}+2is_3e^{-i\pi\alpha}\sin(\pi\alpha),\\
A_{21}& = -ie^{i\pi\alpha}(1-is_3e^{i\pi\alpha})+s_3 + ie^{i\pi\alpha}+2is_3e^{i\pi\alpha}\sin(\pi\alpha),\\
A_{22}&=-ie^{-i\pi\alpha}(1-is_3e^{i\pi\alpha})-s_3 - ie^{i\pi\alpha}-2is_3e^{i\pi\alpha}\sin(\pi\alpha).
\end{align*}
Calculating the coefficient $A_{12}$ gives
\begin{align*}
A_{12}& = -ie^{-i\pi\alpha}+ s_3e^{-2i\pi\alpha} -s_3 + ie^{-i\pi\alpha}+2is_3e^{-i\pi\alpha}\sin(\pi\alpha) \\
&= s_3(\cos(\pi\alpha) - i\sin(\pi\alpha))^2 -s_3 + 2is_3(\cos(\pi\alpha) - i\sin(\pi\alpha))\sin(\pi\alpha) \\
&= s_3(\cos^2(\pi\alpha) -\sin^2(\pi\alpha)) -s_3  + 2s_3\sin^2(\pi\alpha)\\
&= s_3(\cos^2(\pi\alpha) +\sin^2(\pi\alpha)) -s_3 =0
\end{align*}
and similar computations for $A_{21}$ yields
\begin{align*}
A_{21}& = -ie^{i\pi\alpha}- s_3e^{2i\pi\alpha} + s_3 + ie^{i\pi\alpha}+2is_3e^{i\pi\alpha}\sin(\pi\alpha) \\
&= -s_3(\cos(\pi\alpha) + i\sin(\pi\alpha))^2 +s_3 + 2is_3(\cos(\pi\alpha) + i\sin(\pi\alpha))\sin(\pi\alpha) \\
&= -s_3(\cos^2(\pi\alpha) -\sin^2(\pi\alpha)) +s_3  - 2s_3\sin^2(\pi\alpha)\\
&= -s_3(\cos^2(\pi\alpha) +\sin^2(\pi\alpha)) +s_3 =0.
\end{align*}
For the coefficient $A_{11}$ we have
\begin{equation}
\begin{aligned}\label{e-c-a}
A_{11}&=-ie^{i\pi\alpha}+s_3+s_3 - ie^{-i\pi\alpha}-2is_3e^{-i\pi\alpha}\sin(\pi\alpha)\\
&=2s_3-2i\cos(\pi\alpha)-2is_3\cos(\pi\alpha)\sin(\pi\alpha)-2s_3\sin^2(\pi\alpha)\\
&=2s_3\cos^2(\pi\alpha) -2i\cos(\pi\alpha)-2is_3\cos(\pi\alpha)\sin(\pi\alpha)\\
&=2\cos(\pi\alpha)(s_3\cos(\pi\alpha) - i-is_3\sin(\pi\alpha)).
\end{aligned}
\end{equation}
Since the triple $(s_{1},s_{2},s_{3})$ is given by \eqref{stokes-11}, it follows that
\begin{align*}
&s_3\cos(\pi\alpha) -i- is_3\sin(\pi\alpha) \\
&\quad = (-\sin(\pi\alpha)+ik)\cos(\pi\alpha)-i-i(-\sin(\pi\alpha)+ik)\sin(\pi\alpha)\\
&\quad= -\sin(\pi\alpha)\cos(\pi\alpha)+ik\cos(\pi\alpha) +k\sin(\pi\alpha) - i + i\sin^2(\pi\alpha) \\
&\quad= -\sin(\pi\alpha)\cos(\pi\alpha)+ik\cos(\pi\alpha) +k\sin(\pi\alpha) - i\cos^2(\pi\alpha) \\
&\quad= -i\cos(\pi\alpha)(\cos(\pi\alpha) -i\sin(\pi\alpha))+ik(\cos(\pi\alpha)-i\sin(\pi\alpha)) \\
&\quad= (\cos(\pi\alpha)-i\sin(\pi\alpha))(-i\cos(\pi\alpha)+ik) = -ie^{-i\pi\alpha}(\cos(\pi\alpha)-k),
\end{align*}
which after substitution to \eqref{e-c-a} gives
\begin{align*}
A_{11}=-2i\cos(\pi\alpha)(\cos(\pi\alpha)-k)e^{-i\pi\alpha}. 
\end{align*}
It remains to calculate the coefficient $A_{22}$. To this end let us observe that
\begin{equation}
\begin{aligned}\label{e-c-b}
A_{22}&=-ie^{-i\pi\alpha}(1-is_3e^{i\pi\alpha})-s_3 - ie^{i\pi\alpha}-2is_3e^{i\pi\alpha}\sin(\pi\alpha)\\
&=-ie^{-i\pi\alpha}-s_3-s_3 - ie^{i\pi\alpha}-2is_3e^{i\pi\alpha}\sin(\pi\alpha)\\
&=-2s_3-2i\cos(\pi\alpha)-2is_3\cos(\pi\alpha)\sin(\pi\alpha)+2s_3\sin^2(\pi\alpha)\\
&=-2\cos(\pi\alpha)(s_3\cos(\pi\alpha) +i+is_3\sin(\pi\alpha)).
\end{aligned}
\end{equation}
Using \eqref{stokes-11} once again we obtain
\begin{align*}
&s_3\cos(\pi\alpha) +i+is_3\sin(\pi\alpha) \\
&\quad = (-\sin(\pi\alpha)+ik)\cos(\pi\alpha)+i+i(-\sin(\pi\alpha)+ik)\sin(\pi\alpha)\\
&\quad= -\sin(\pi\alpha)\cos(\pi\alpha)+ik\cos(\pi\alpha) -k\sin(\pi\alpha) + i - i\sin^2(\pi\alpha) \\
&\quad= -\sin(\pi\alpha)\cos(\pi\alpha)+ik\cos(\pi\alpha) -k\sin(\pi\alpha) + i\cos^2(\pi\alpha) \\
&\quad= i\cos(\pi\alpha)(\cos(\pi\alpha) +i\sin(\pi\alpha))+ik(\cos(\pi\alpha)+i\sin(\pi\alpha)) \\
&\quad= (\cos(\pi\alpha)+i\sin(\pi\alpha))(i\cos(\pi\alpha)+ik) = ie^{i\pi\alpha}(\cos(\pi\alpha)+k),
\end{align*}
which together with \eqref{e-c-b} provides
\begin{equation*}
A_{22}=-2i\cos(\pi\alpha)(\cos(\pi\alpha)+k)e^{i\pi\alpha}.
\end{equation*}
Therefore we have
\begin{equation}\label{ec-aa}
A = -2i\cos(\pi\alpha)\begin{pmatrix}(\cos(\pi\alpha)-k)e^{-i\pi\alpha}&0\\[5pt]
0&(\cos(\pi\alpha)+k)e^{i\pi\alpha}\end{pmatrix}.
\end{equation}
In view of \eqref{matrix-kk} and \eqref{ec-aa}, we have
\begin{align*}
K&= \frac{1}{(1-s_1s_3)^{1/2}}\begin{pmatrix}p&0\\[5pt] 0& q\end{pmatrix}
\begin{pmatrix}(\cos(\pi\alpha)-k)e^{-i\pi\alpha}&\hspace{-5pt}0\\[5pt]
0&\hspace{-5pt}(\cos(\pi\alpha)+k)e^{i\pi\alpha}\end{pmatrix}\begin{pmatrix}p& 0\\[5pt] 0& q\end{pmatrix}^{-1} \\
&= \frac{1}{(1-s_1s_3)^{1/2}}\begin{pmatrix}\cos(\pi\alpha)-k&0\\[5pt]
0&\cos(\pi\alpha)+k\end{pmatrix} e^{-i\pi\alpha\sigma_3}
\end{align*}
and the proof is completed. \hfill $\square$ \\

In the following proposition we obtain the representation formula for the matrix $P(x)$ provided $x<0$ and $|x|$ is sufficiently large.

\begin{proposition}\label{prop-p-l}
There is $x_{-}<0$ such that, for any $x<x_{-}$, we have
\begin{align*}
P(x)=\frac{1}{2}R(0,(-x)^{3/2})e^{-i\frac{\pi}{4}\sigma_3}\begin{pmatrix}1&1-2\alpha \\[5pt]-1&1-2\alpha\end{pmatrix}DK e^{2\pi i\alpha\sigma_{3}}(-ix)^{\alpha\sigma_3},
\end{align*}
where the function $R(z,t)$ satisfies the asymptotic condition \eqref{aaa1}.
\end{proposition}
\proof We write $\Lambda:=\{z\in\C \ | \ \mathrm{arg}\,z = \pi/12\}$. Proposition \ref{prop-con} implies the existence of $x_{-}<0$ with the property that, for any $x<x_{-}$, there is $\ve >0$ such that
\begin{align*}
Z(\lambda,x)=R(z)E(z)\hat\Phi^{0}(e^{2\pi i}it\eta(z))D K z^{-\alpha\sigma_3}e^{t\tilde\theta(z)\sigma_{3}}t^{-\frac{\alpha\sigma_3}{3}},\quad |\lambda|<\ve, \ \lambda\in\Lambda,
\end{align*}
where the change of variables is given by $\lambda = (-x)^{1/2}z$ and $t=(-x)^{3/2}$.
By Proposition \ref{prop-diag}, the matrix $DK$ is diagonal, which implies that
\begin{align}\label{eqq-l4}
Z(\lambda,x)=R(z)E(z)\hat\Phi^{0}(e^{2\pi i}it\eta(z))z^{-\alpha\sigma_3}DK e^{t\tilde\theta(z)\sigma_{3}}t^{-\frac{\alpha\sigma_3}{3}}, \ |\lambda|<\ve, \, \lambda\in\Lambda.
\end{align}
Lemma \ref{lem-lim-1} says that the function $z\mapsto \hat\Phi^{0}(z)z^{-\alpha\sigma_3}$ is holomorphic on the complex plane and the convergence holds
\begin{align}\label{eqq-l3}
\lim_{z\to 0}\hat\Phi^{0}(z)z^{-\alpha\sigma_3}=\frac{1}{2}e^{-i\frac{\pi}{4}\sigma_3}\begin{pmatrix}1&1-2\alpha \\[5pt]-1&1-2\alpha\end{pmatrix}.
\end{align}
On the other hand, we have 
\begin{align*}
\lim_{z\in\Lambda,z\to 0}(e^{2\pi i} it\eta(z))^{\alpha\sigma_3} z^{-\alpha\sigma_3} & 
= \lim_{z\in\Lambda,z\to 0} e^{2\pi i\alpha\sigma_{3}}[it\eta(z)/z]^{\alpha\sigma_3} = e^{2\pi i\alpha}(it)^{\alpha\sigma_3},
\end{align*}
which together with \eqref{eqq-l3} and the equality
\begin{align*}
\hat\Phi^{0}(e^{2\pi i}it\eta(z))z^{-\alpha\sigma_3} = \hat\Phi^{0}(e^{2\pi i}it\eta(z))(e^{2\pi i}it\eta(z))^{-\alpha\sigma_3} (e^{2\pi i}it\eta(z))^{\alpha\sigma_3} z^{-\alpha\sigma_3}
\end{align*}
give the following limit
\begin{align}\label{eqq-l5}
\lim_{z\in\Lambda, z\to 0}\hat\Phi^{0}(e^{2\pi i}it\eta(z))z^{-\alpha\sigma_3} = \frac{1}{2}e^{-i\frac{\pi}{4}\sigma_3}\begin{pmatrix}1&1-2\alpha \\[5pt]-1&1-2\alpha\end{pmatrix}e^{2\pi i\alpha\sigma_{3}}(it)^{\alpha\sigma_3}.
\end{align}
Let us observe that the function $E(z)$ is holomorphic in the neighborhood of the origin and $E(z)\to I$ as $z\to 0$. Combining this with \eqref{eqq-l5}, we pass in the equation \eqref{eqq-l4} to the limit with $z\to 0$ along the ray $\Lambda$ and deduce that
\begin{align*}
P(x)&=\lim_{\lambda\in\Lambda,\lambda\to 0}Z(\lambda,x) =\lim_{z\in\Lambda,z\to 0}R(z)E(z)\hat\Phi^{0}(e^{2\pi i}it\eta(z))z^{-\alpha\sigma_3}DK e^{t\tilde\theta(z)\sigma_{3}}t^{-\frac{\alpha\sigma_3}{3}}\\
& = \frac{1}{2}R(0,t)e^{-i\frac{\pi}{4}\sigma_3}\begin{pmatrix}1&1-2\alpha \\[5pt]-1&1-2\alpha\end{pmatrix}e^{2\pi i\alpha\sigma_{3}}(it)^{\alpha\sigma_3}DK t^{-\frac{\alpha\sigma_3}{3}}.
\end{align*}
In view of the fact that the matrix $DK$ is diagonal we obtain
\begin{align*}
P(x)& = \frac{1}{2}R(0,t)e^{-i\frac{\pi}{4}\sigma_3}\begin{pmatrix}1&1-2\alpha \\[5pt]-1&1-2\alpha\end{pmatrix}DK e^{2\pi i\alpha\sigma_{3}}(it)^{\alpha\sigma_3}t^{-\frac{\alpha\sigma_3}{3}}\\
& = \frac{1}{2}R(0,(-x)^{3/2})e^{-i\frac{\pi}{4}\sigma_3}\begin{pmatrix}1&1-2\alpha \\[5pt]-1&1-2\alpha\end{pmatrix}DK e^{2\pi i\alpha\sigma_{3}}(-ix)^{\alpha\sigma_3},
\end{align*}
and the proof of proposition is completed. \hfill $\square$

\section{Proof of Theorems \ref{th-total} and \ref{th-total-2}}

\noindent Propositions \ref{prop-p-g} and \ref{prop-p-l} say that, there is $x_{0}>0$ such that, for $x>x_{0}$, the functions $P(x)$ and $P(-x)$ have the following forms
\begin{align*}
P(x)&= \frac{1}{2}\chi^{1}(0,x)\chi^{2}(0,x) e^{-i\frac{\pi}{4}\sigma_3}\begin{pmatrix}1&1\\[5pt]-1&1\end{pmatrix}\begin{pmatrix}1&0\\[5pt] 0 &1-2\alpha\end{pmatrix}e^{2\pi i\alpha\sigma_{3}}(-ix)^{\alpha\sigma_3} D,\\
P(-x)&=\frac{1}{2}R(0,x^{3/2})e^{-i\frac{\pi}{4}\sigma_3}\begin{pmatrix}1&1\\[5pt]-1&1\end{pmatrix}\begin{pmatrix}1&0\\[5pt] 0 &1-2\alpha\end{pmatrix}DK e^{2\pi i\alpha\sigma_{3}}(ix)^{\alpha\sigma_3}.
\end{align*}
Let us observe that, for any $x>0$, we have
\begin{align*}
(ix)^{\alpha\sigma_3}(-ix)^{-\alpha\sigma_3} = \left(|x|e^{i\frac{\pi}{2}}\right)^{\alpha\sigma_3} \left(|x|e^{-i\frac{\pi}{2}}\right)^{-\alpha\sigma_3} = e^{i\frac{\pi}{2}\alpha\sigma_3}e^{i\frac{\pi}{2}\alpha\sigma_3} =
e^{i\pi\alpha\sigma_3},
\end{align*}
which together with the fact that matrices $D$ and $K$ are diagonal, imply that
\begin{align}\label{eq-p-p-inv}
P(-x)P(x)^{-1} = R(0,x^{3/2})H\chi^{2}(0,x)^{-1}\chi^{1}(0,x)^{-1}, \quad x>x_{0},
\end{align}
where the matrix $H$, is given by 
\begin{align*}
H:=\frac{1}{2}e^{-i\frac{\pi}{4}\sigma_3}\begin{pmatrix}1&1\\[5pt]-1&1\end{pmatrix}Ke^{i\pi\alpha\sigma_3} 
\begin{pmatrix}1&-1\\[5pt]1&1\end{pmatrix}e^{i\frac{\pi}{4}\sigma_3}.
\end{align*}
Writing $h_{\pm}:=(1-s_1s_3)^{-\frac{1}{2}}(\cos(\pi\alpha)\mp k)/2$ and applying Proposition \ref{prop-diag}, we obtain
\begin{align*}
H&=e^{-i\frac{\pi}{4}\sigma_3}\begin{pmatrix}1&1\\[5pt]-1&1\end{pmatrix} \begin{pmatrix}h_{+}&0\\[5pt]0&h_{-}\end{pmatrix}\begin{pmatrix}1&-1\\[5pt]1&1\end{pmatrix}e^{i\frac{\pi}{4}\sigma_3}\\
&=e^{-i\frac{\pi}{4}\sigma_3}\begin{pmatrix}h_{+}+h_{-}&h_{-}-h_{+}\\[5pt]h_{-}-h_{+}&h_{+}+h_{-}\end{pmatrix}e^{i\frac{\pi}{4}\sigma_3}=
\begin{pmatrix}h_{+}+h_{-}&i(h_{+}-h_{-})\\[5pt]i(h_{-}-h_{+})&h_{+}+h_{-}\end{pmatrix},
\end{align*}
which together with \eqref{est-kk}, \eqref{asy-11bb}, \eqref{aaa1} and \eqref{eq-p-p-inv} give
\begin{align*}
\lim_{x\to+\infty} P(x)P(-x)^{-1} = H^{-1} = \begin{pmatrix}h_{+}+h_{-}&i(h_{-}-h_{+})\\[5pt]i(h_{+}-h_{-})&h_{+}+h_{-}\end{pmatrix}.
\end{align*}
Let us observe that the formula \eqref{eq-p-2} implies that
$$\exp\left(\int_{-x}^x u(y;\alpha,k)\,dy\right) = [P(x)P(-x)^{-1}]_{11} + i [P(x)P(-x)^{-1}]_{21},\quad x>0$$
and consequently 
\begin{align*}
&\lim_{x\to+\infty} \exp\left(\int_{-x}^x u(y;\alpha,k)\,dy\right) \!= \!\lim_{x\to+\infty}[P(x)P(-x)^{-1}]_{11}\!+\! i\lim_{x\to+\infty}[P(x)P(-x)^{-1}]_{21}\\
& \quad = 2h_{-} =(1-s_1s_3)^{-1/2}(\cos(\pi\alpha)+k) = \left(\cos^2(\pi\alpha)-k^2\right)^{-1/2}(\cos(\pi\alpha)+k).
\end{align*}
which gives the desired total integral formula 
\begin{align}\label{t-int-f}
\lim_{x\to+\infty} \exp\left(\int_{-x}^x u(y;\alpha,k)\,dy\right) = \left(\cos^2(\pi\alpha)-k^2\right)^{-1/2}(\cos(\pi\alpha)+k).
\end{align}
in the case $u$ is the purely imaginary AS solution. If $u$ is the real AS solution, then $\mathrm{Im}\, u(x;\alpha,k) = 0$ for $x\in\R$, which enables us to take both side logarithm of the limit \eqref{t-int-f} to obtain \eqref{wz1}. This completes the proof of Theorems \ref{th-total} and \ref{th-total-2}. \hfill $\square$

\parindent = 0 pt

\end{document}